\documentclass[12pt]{article}

\usepackage{latexsym,amsfonts,amsmath,amssymb,graphicx,hyperref}

\newcommand{\comment}[1]{}

\newcommand{\Om}{\Omega}

\newcommand{\St}{\mathbb{S}}
\newcommand{\EE}{\mbox{\bf E}\,}
\newcommand{\PP}{\mbox{\bf P}\,}
\newcommand{\QQ}{\mbox{\bf Q}\,}

\newcommand{\R}{\mathbb{R}}
\newcommand{\C}{\mathbb{C}}
\newcommand{\Q}{\mathbb{Q}}
\newcommand{\HH}{\mathbb{H}}
\newcommand{\N}{\mathbb{N}}

\newcommand{\Z}{\mathbb{Z}}

\newcommand{\pa}{\partial}

\newcommand{\F}{{\cal F}}

\newcommand{\no}{\noindent}

\newcommand{\mr}{\mathring}
\newcommand{\BGE}{\begin{equation}}
\newcommand{\BGEN}{\begin{equation*}}
\newcommand{\EDE}{\end{equation}}
\newcommand{\EDEN}{\end{equation*}}

\def\eps{\varepsilon}
\def\til{\widetilde}
\def\ha{\widehat}
\def\sem{\setminus}

\def\lin{\overline}
\def\vphi{\varphi}

\def\del{\delta}

\DeclareMathOperator{\sign}{sign} 
\DeclareMathOperator{\dist}{dist} 
\DeclareMathOperator{\hcap}{hcap} \DeclareMathOperator{\id}{id}
\DeclareMathOperator{\Imm}{Im } \DeclareMathOperator{\Ree}{Re }
 
 \DeclareMathOperator{\HP}{HP}
\DeclareMathOperator{\scap}{scap} 
\DeclareMathOperator{\out}{out}

\def\h0{{\bf h}}

\newtheorem{Lemma}{Lemma}[section]
\newtheorem{Theorem}{Theorem}[section]
\newtheorem{Definition}{Definition}[section]
\newtheorem{Corollary}{Corollary}[section]
\newtheorem{Proposition}{Proposition}[section]
\newtheorem{Conjecture}{Conjecture}
\numberwithin{equation}{section}

\begin{document}
\title{\bf Duality of Chordal SLE}
\date{\today}
\author{Dapeng Zhan}
\maketitle
\begin{abstract}
We derive some geometric properties of chordal SLE$(\kappa;\vec{\rho})$
processes. Using these results and the method of coupling two SLE processes, we prove that the
outer boundary of the final hull of a chordal SLE$(\kappa;\vec{\rho})$
process has the same distribution as the
image of a chordal SLE$(\kappa';\vec{\rho'})$ trace, where
$\kappa>4$, $\kappa'=16/\kappa$, and the forces $\vec{\rho}$ and
$\vec{\rho'}$ are suitably chosen. We find that for $\kappa\ge 8$, the
boundary of a standard chordal SLE$(\kappa)$ hull stopped on swallowing
a fixed $x\in\R\sem\{0\}$ is the image of some SLE$(16/\kappa;\vec{\rho})$ trace
started from $x$. Then we obtain a new proof of the fact that chordal SLE$(\kappa)$
trace is not reversible for $\kappa>8$. We also prove that the reversal of
SLE$(4;\vec{\rho})$ trace has the same distribution as the
time-change of some SLE$(4;\vec{\rho'})$ trace for certain values of
$\vec{\rho}$ and $\vec{\rho'}$.
\end{abstract}

\section{Introduction}
The Schramm-Loewner evolution (SLE) has become a fast growing area
in Probability Theory since 1999 (\cite{S-SLE}). SLE describes some
random fractal curve, which is called an SLE trace, that grows in a plane
domain. The behavior of the trace depends on a real parameter
$\kappa>0$. We write SLE$(\kappa)$ to emphasize the parameter
$\kappa$. If $\kappa\in(0,4]$, the trace is a simple curve; if
$\kappa>4$, the trace is not simple; and if $\kappa\ge 8$, the trace
is space-filling. For basic properties of SLE, see \cite{LawSLE} and \cite{RS-basic}.

Many two-dimensional lattice models from statistical physics have
been proved to have SLE as their scaling limits when the mesh of the
grid tends to $0$, e.g., the convergence of critical percolation
on triangular lattice to SLE$(6)$
(\cite{SS-6}), loop-erased random walk (LERW) to SLE$(2)$
(\cite{LSW-2}\cite{LERW}), uniform spanning tree (UST) Peano curve to SLE$(8)$
(\cite{LSW-2}), Gaussian free field contour line to SLE$(4)$ (\cite{SS}), and
some Ising models to SLE$(3)$ and SLE$(16/3)$ (\cite{SS-3}). And
there are some promising conjectures, e.g., the convergence of
self-avoiding walk to SLE$(8/3)$ (\cite{LSW-8/3}), and double domino
tilling to SLE$(4)$ (\cite{RS-basic}).

For $\kappa>4$, people are also interested in the hulls that are
generated by the SLE$(\kappa)$ traces. Duplantier proposed a rough conjecture
about the duality between SLE$(\kappa)$ and SLE$(16/\kappa)$,
which says that when $\kappa>4$, the boundary of an SLE$(\kappa)$
hull looks locally like an SLE$(16/\kappa)$ trace.

For $\kappa\le 8$, the Hausdorff dimension of an SLE$(\kappa)$ trace
was proved to be $1+\kappa/8$ (\cite{dim-SLE}). If the duality
conjecture is true, then we may conclude that for $\kappa>4$, the
Hausdorff dimension of the boundary of an SLE$(\kappa)$ hull is
$1+2/\kappa$.

For some parameter $\kappa$, the duality is already known. The
duality between SLE$(8)$ and SLE$(2)$ follows from the convergence
of UST and LERW to SLE$(8)$ and SLE$(2)$, respectively, and the
Wilson's algorithm (\cite{Wilson}) that links UST with LERW. The
duality between SLE$(6)$ and SLE$(8/3)$ follows from the conformal
restriction property (\cite{LSW-8/3}). The duality between
SLE$(16/3)$ and SLE$(3)$ follows from the convergence of Ising
models.

In \cite{Julien-Duality}, J.\ Dub\'edat proposed some specific
conjectures about the duality of SLE, one of which says that for $\kappa>4$,
the right boundary of the final hull of a chordal
SLE$(\kappa;\kappa-4)$ process started from $(0,0^+)$ has the same law as a chordal
SLE$(\kappa';\frac12(\kappa'-4))$ trace started from $(0,0^-)$, where $\kappa'=16/\kappa$.
And he justified his conjecture by studying the distributions of the
sets obtained by adding Brownian loop soups to
SLE$(\kappa;\kappa-4)$ and SLE$(\kappa';\frac12(\kappa'-4))$,
respectively.

Recently, a new technique about constructing a coupling of two SLE
processes that grow in the same domain was introduced
(\cite{reversibility}) to prove the reversibility of chordal
SLE$(\kappa)$ trace when $\kappa\in(0,4]$. In this paper, we will
use this technique to prove some specific versions of the duality
conjecture, which are not exactly the same as those in
\cite{Julien-Duality}. For example, one of our results is that for
$\kappa>4$ and $\kappa'=16/\kappa$, the
right boundary of the final hull of a chordal
SLE$(\kappa;\kappa-4)$ process started from $(0,0^+)$ has the same law as
the image under the map $z\mapsto 1/\lin z$ of a chordal
SLE$(\kappa';\frac12(\kappa'-4))$ trace started from $(0,0^-)$.
If the degenerate chordal
SLE$(\kappa';\frac12(\kappa'-4))$ trace satisfies reversibility, which is
Conjecture \ref{conjec} of this paper, then
Dub\'edat's conjecture is proved.

This paper is organized in the following way. In Section
\ref{prelim}, we review the definitions of the chordal and strip
(i.e., dipolar) Loewner equations and
SLE$(\kappa;\vec{\rho})$ processes. The conformal invariance of
chordal and strip SLE$(\kappa;\vec{\rho})$ processes are introduced.
In Section \ref{geom}, we study the tail behavior of a chordal or
strip SLE$(\kappa;\vec{\rho})$ trace when the force points and
forces satisfy certain conditions. In Section \ref{coupl}, for
$\kappa\ge 4\ge\kappa'>0$ with $\kappa\kappa' =16$, some commutation result
of a chordal SLE$(\kappa;\vec{\rho})$ process with a chordal
SLE$(\kappa';\vec{\rho'})$ process is described in terms of a two-dimensional
martingale. This is closely
related with J.\ Dub\'edat's work in \cite{Julien-Comm}. Then
the technique in \cite{reversibility} is applied to get a coupling
of the above two SLE processes. In
Section \ref{appli}, we consider the coupling in the previous
section with some special choices of force points and forces, and
apply the geometry results from Section \ref{geom} to prove that in
this coupling, the chordal SLE$(\kappa';\vec{\rho'})$ trace becomes
the outer boundary of the chordal SLE$(\kappa;\vec{\rho})$ hull, and
so prove the duality conjecture. Then we derive the equation of the
boundary of a standard chordal SLE$(\kappa)$ hull, $\kappa\ge 8$,
at the time when a fixed $x\in\R\sem\{0\}$ is swallowed. Then we
give a new proof of the fact that for $\kappa>8$, the
chordal SLE$(\kappa)$ trace is not reversible. This result was claimed
in \cite{RS-basic}. At the end, we derive the
reversibility property of some chordal SLE$(4;\rho)$ traces.

\section{Preliminary} \label{prelim}
\subsection{Chordal SLE} If $H$ is a bounded and relatively  closed
subset of $\HH=\{z\in\C:\Imm z>0\}$, and $\HH\sem H$ is simply
connected, then we call $H$ a hull in $\HH$ w.r.t.\ $\infty$. For
such $H$, there is $\vphi_H$ that maps $\HH\sem H$ conformally onto
$\HH$, and satisfies $\vphi_H(z)=z+\frac{c}{z}+O(\frac 1{z^2})$ as
$z\to\infty$, where $c=\hcap(H)\ge 0$ is called the capacity of $H$
in $\HH$ w.r.t.\ $\infty$. If $H_1\subset H_2$ are hulls
in $\HH$ w.r.t.\ $\infty$, then $H_2/H_1:=\vphi_{H_1}(H_2\sem H_1)$
is also a hull in $\HH$ w.r.t.\ $\infty$, and
$\hcap(H_2/H_1)=\hcap(H_2)-\hcap(H_1)$. If $H_1\subset H_2\subset
H_3$ are three hulls in $\HH$ w.r.t.\ $\infty$, then $H_2/H_1\subset
H_3/H_1$ and $(H_3/H_1)/(H_2/H_1)=H_3/H_2$.

\begin{Proposition} Suppose $\Om$ is an open neighborhood of $x_0\in\R$
in $\HH$. Suppose $W$ maps $\Om$ conformally into $\HH$ such that
for some $r>0$, if $z\in\Om$ approaches $(x_0-r,x_0+r)$ then $W(z)$ approaches
$\R$. So
$W$ extends conformally across $(x_0-r,x_0+r)$ by Schwarz reflection
principle. Then for any $\eps>0$, there is some $\del>0$ such that
if a hull $H$ in $\HH$ w.r.t.\ $\infty$ is contained in
$\{z\in\HH:|z-x_0|<\del\}$, then $W(H)$ is also a hull in $\HH$
w.r.t.\ $\infty$, and
$$|\hcap(W(H))-W'(x_0)^2\hcap(H)|\le\eps|\hcap(H)|.$$
\label{hcap}
\end{Proposition}\vskip -9mm
{\bf Proof.} This is Lemma 2.8 in \cite{LSW1}. $\Box$

\vskip 4mm

For a real interval $I$, we use $C(I)$ to denote the space of real
continuous functions on $I$. For $T>0$ and $\xi\in C([0,T))$, the
chordal Loewner equation driven by $\xi$ is
$$\pa_t\vphi(t,z)=\frac{2}{\vphi(t,z)-\xi(t)},\quad \vphi(0,z)=z.$$
For $0\le t<T$, let $K(t)$ be the set of $z\in\HH$ such that the
solution $\vphi(s,z)$ blows up before or at time $t$. We call $K(t)$
and $\vphi(t,\cdot)$, $0\le t<T$, chordal Loewner hulls and maps,
respectively, driven by $\xi$.

\begin{Definition}
We call $(K(t),0\le t<T)$ a Loewner chain in $\HH$ w.r.t.\ $\infty$, if each $K(t)$ is a
hull in $\HH$ w.r.t.\ $\infty$; $K(0)=\emptyset$; $K(s)\subsetneqq
K(t)$ if $s<t$; and for each fixed $a\in(0,T)$ and compact $F\subset\HH\sem K(a)$,
 the extremal length (\cite{Ahl}) of the curves in $\HH\sem K(t+\eps)$
that disconnect $K(t+\eps)\sem K(t)$ from $F$ tends to $0$ as $\eps\to 0^+$, uniformly in
$t\in[0,a]$.
\end{Definition}

\begin{Proposition} (a) Suppose $K(t)$ and $\vphi(t,\cdot)$, $0\le t<T$,
are chordal Loewner hulls and maps, respectively, driven by $\xi\in
C([0,T))$. Then $(K(t),0\le t<T)$ is a Loewner chain in $\HH$
w.r.t.\ $\infty$, $\vphi_{K(t)}=\vphi(t,\cdot)$, and
$\hcap(K(t))=2t$ for any $0\le t<T$. Moreover, for every
$t\in[0,T)$, $$ \{\xi(t)\}=\bigcap_{\eps\in(0,T-t)}\lin{K(t+\eps)/
K(t)}.$$ (b) Let $(L(s),0\le s<S)$ be a
Loewner chain in $\HH$ w.r.t.\ $\infty$. Let $v(s)=\hcap(L(s))/2$,
$0\le s<S$. Then $v$ is a continuous (strictly) increasing function with
$u(0)=0$. Let $T=v(S)$ and $K(t)=L(v^{-1}(t))$, $0\le t<T$. Then
$K(t)$, $0\le t<T$, are chordal Loewner hulls driven by some $\xi\in
C([0,T))$. \label{chordal-Loewner-chain}
\end{Proposition}
{\bf Proof.} This is almost the same as Theorem 2.6 in \cite{LSW1}.
$\Box$

\vskip 4mm

Let $D$ be a domain and $K\subset D$. Let $p_1$ and $p_2$ be two
boundary points or prime ends of $D$. We say that $K$ does not
separate $p_1$ from $p_2$ in $D$ if there are neighborhoods $U_1$
and $U_2$ of $p_1$ and $p_2$, respectively, in $D$ such that $U_1$
and $U_2$ lie in the same pathwise connected component of $D\sem K$.
In our definition, $K$ may separates some $p$ from itself. Let $Q$
be a set of boundary points or prime ends of $D$. We say that $K$
does not divide $Q$ in $D$ if for any $p_1,p_2\in D$, $K$ does not
separate $p_1$ from $p_2$ in $D$.

Let $\vphi(t,\cdot)$ and $K(t)$ be as before. Let $x\in\R$. If at
time $t$, $\vphi(t,x)$ does not blow up, then
$K(t)$ does not separate $x$ from $\infty$ in $\HH$, and vice versa.
In fact, we have a slightly
stronger result: if $\vphi(s,x)$ blows up before or at
$s=t\in[0,T)$, then $\cup_{s<t} K(s)$ also separates $x$ from
$\infty$ in $\HH$. This follows from the property of a Loewner
chain.

Let $B(t)$, $0\le t<\infty$, be a (standard linear) Brownian motion.
Let $\kappa \ge 0$. Then $K(t)$ and $\vphi(t,\cdot)$, $0\le
t<\infty$, driven by $\xi(t)=\sqrt\kappa B(t)$, $0\le t<\infty$, are
called standard chordal SLE$(\kappa)$ hulls and maps, respectively.
It is known (\cite{RS-basic}\cite{LSW-2}) that almost surely for any
$t\in[0,\infty)$, \BGE\beta(t):=\lim_{\HH\ni
z\to\xi(t)}\vphi(t,\cdot)^{-1}(z)\label{trace}\EDE exists, and
$\beta(t)$, $0\le t<\infty$, is a continuous curve in $\lin{\HH}$.
Moreover, if $\kappa\in(0,4]$ then $\beta$ is a simple curve, which
intersects $\R$ only at the initial point, and for any $t\ge 0$,
$K(t)=\beta((0,t])$; if $\kappa>4$ then $\beta$ is not simple, and
intersects $\R$ at infinitely many points; and in general, $\HH\sem
K(t)$ is the unbounded component of $\HH\sem\beta((0,t])$ for any
$t\ge 0$. Such $\beta$ is called a standard chordal SLE$(\kappa)$
trace.

If $(\xi(t))$ is a semi-martingale, and $d\langle \xi(t)\rangle
=\kappa dt$ for some $\kappa>0$, then from Girsanov theorem and the
existence of standard chordal SLE$(\kappa)$ trace, almost surely for
any $t\in[0,T)$, $\beta(t)$ defined by (\ref{trace}) exists, and has
the same property as a standard chordal SLE$(\kappa)$ trace
(depending on the value of $\kappa$) as described in the last
paragraph.

Let $\kappa\ge 0$, $\rho_1,\dots,\rho_N\in\R$, $x\in\R$, and
$p_1,\dots,p_N\in\ha\R\sem\{x\}$, where $\ha\R=\R\cup\{\infty\}$ is
a circle. Let $\xi(t)$ and $p_k(t)$, $1\le k\le N$, be the solutions
to the SDE:
\begin{equation}\left\{\begin{array}{lll} d\xi(t) & = &
\sqrt\kappa d B(t)+\sum_{k=1}^N\frac{\rho_k\, dt}{\xi(t)-p_k(t)}\\
dp_k(t) & = & \frac{2dt}{p_k(t)-\xi(t)},\quad 1\le k\le
N,\end{array}\right.\label{kappa-rho}\end{equation} with initial
values $\xi(0)=x$ and $p_k(0)=p_k$, $1\le k\le N$. If
$\vphi(t,\cdot)$ are chordal Loewner maps driven by $\xi(t)$, then
$p_k(t)=\vphi(t,p_k)$. Here if some $p_k=\infty$ then
$p_k(t)=\infty$ and $\frac{\rho_k}{\xi(t)-p_k(t)}=0$ for all $t\ge
0$, so $p_k$ has no effect on the equation. Suppose $[0,T)$ is the
maximal interval of the solution. Let $K(t)$, $0\le t<T$, be chordal
Loewner hulls driven by $\xi$. Then we call $K(t)$, $0\le t<T$, a
(full) chordal SLE$(\kappa;\rho_1,\dots,\rho_N)$ process started
from $(x;p_1,\dots,p_N)$. Since $(\xi(t))$ is a semi-martingale, and
$d\langle \xi(t)\rangle =\kappa dt$, so the chordal Loewner trace
$\beta(t)$, $0\le t<T$, driven by $\xi$ exists, and is called a chordal
SLE$(\kappa;\rho_1,\dots,\rho_N)$ trace started from
$(x;p_1,\dots,p_N)$. If we let $\vec\rho$ and $\vec p$ to denote the
vectors $(\rho_1,\dots,\rho_N)$ and $(p_1,\dots,p_N)$, then we may
call $K(t)$ and $\beta(t)$, $0\le t<T$, chordal
SLE$(\kappa;\vec{\rho})$ process and trace, respectively, started
from $(x;\vec{p})$. If $S\in(0,T]$ is a stopping time, then
$K(t)$ and $\beta(t)$, $0\le t<S$, are called partial chordal
SLE$(\kappa;\vec{\rho})$ process and trace, respectively, started
from $(x;\vec{p})$.

These $p_k$'s and $\rho_k$'s are called force points and forces, respectively.
For $0\le t<T$ and $1\le k\le N$,
$\vphi(t,p_k)$ does not blow up, so $K(t)$ does not divide
$\{\infty,p_1,\dots,p_N\}$ in $\HH$. If $T<\infty$ then there must
exist some $p_k\in\R$ such that $\vphi(t,p_k)-\xi(t)\to 0$ as $t\to
T$, so $\cup_{t<T}K(t)$ separates $p_k$ from $\infty$ in $\HH$. If
$T=\infty$ then $\cup_{t<T}K(t)$ is unbounded, so $\cup_{t<T}K(t)$
separates $\infty$ from itself in $\HH$. Thus in any case,
$\cup_{t<T}K(t)$ divides $\{\infty,p_1,\dots,p_N\}$ in $\HH$.

The chordal SLE$(\kappa;\vec{\rho})$ defined above are of generic
cases. We now introduce degenerate SLE$(\kappa;\vec{\rho})$, where
one of the force points takes value $x^+$ or $x^-$, or two of the
force points take values $x^+$ and $x^-$, respectively, where
$x\in\R$ is the initial point of the trace. Let $\kappa \ge 0$;
$\rho_1,\dots,\rho_N\in\R$, and $\rho_1\ge \kappa/2-2$; $p_1=x^+$,
$p_2,\dots,p_N\in \ha\R\sem\{x\}$. Let $\xi(t)$ and $p_k(t)$, $1\le
k\le N$, $0<t<T$, be the maximal solution to (\ref{kappa-rho}) with
initial values $\xi(0)=p_1(0)=x$, and $p_k(0)=p_k$, $1\le k\le N$.
Moreover, we require that $p_1(t)>\xi(t)$ for any $0<t<T$. If $N=1$,
the existence of the solution follows from the Bessel Process (see
\cite{LSW-8/3}). The condition $\rho_1\ge \kappa/2-2$ is to
guarantee that $p_1$ is not immediately swallowed after time $0$. If
$N\ge 2$, the existence of the solution follows from the above
result and Girsanov Theorem. Then we obtain chordal
SLE$(\kappa;\rho_1,\dots,\rho_N)$ process and trace started from
$(x;x^+,p_2,\dots,p_N)$. If the condition $p_1(t)>\xi(t)$ is
replaced by $p_1(t)<\xi(t)$, then we get chordal
SLE$(\kappa;\rho_1,\dots,\rho_N)$ process and trace started from
$(x;x^-,p_2,\dots,p_N)$. Now suppose $N\ge 2$, $\rho_1,\rho_2\ge
\kappa/2-2$, $p_1=x^+$, and $p_2=x^-$. Let $\xi(t)$ and $p_k(t)$,
$1\le k\le N$, $0<t<T$, be the maximal solution to (\ref{kappa-rho})
with initial values $\xi(0)=p_1(0)=p_2(0)=x$, and $p_k(0)=p_k$,
$1\le k\le N$, such that $p_1(t)>\xi(t)>p_2(t)$ for all $0<t<T$.
Then we obtain chordal SLE$(\kappa;\rho_1,\dots,\rho_N)$ process and
trace started from $(x;x^+,x^-,p_3,\dots,p_N)$. The existence of the
solution to the equation equation follows from \cite{SS} and
Girsanov Theorem.

The force point $x^+$ or $x^-$ is called a degenerate force point.
Other force points are called generic force points. Let
$\vphi(t,\cdot)$ be the chordal Loewner maps driven by $\xi$. Since
for any generic force point $p_j$, we have $p_j(t)=\vphi(t,p_j)$, so
it is reasonable to write $\vphi(t,p_j)$ for $p_j(t)$ in the case
that $p_j$ is a degenerate force point. Suppose $\rho_j$ is the
force associated with some degenerate force point $p_j$. If we allow
that the process continues growing after $p_j$ is swallowed, the
condition that $\rho_j\ge \kappa/2-2$ may be weakened to $\rho_j>-2$
(\cite{LSW-8/3}).

>From the work in \cite{SW}, we get the conformal invariance of
chordal SLE$(\kappa;\vec{\rho})$ processes, which is the following
lemma.

\begin{Lemma} Suppose $\kappa\ge 0$ and
$\vec{\rho}=(\rho_1,\dots,\rho_N)$ with  $\sum_{m=1}^N\rho_m =\kappa-6$.
For $j=1,2$, let $K_j(t)$, $0\le
t<T_j$, be a generic or degenerate chordal SLE$(\kappa;\vec{\rho})$ process started from
$(x_{j};\vec{p}_{j})$, where $\vec{p}_{j}=(p_{j,1},\dots,p_{j,N})$,
$j=1,2$. Suppose $W$ is a conformal
or conjugate conformal map from $\HH$ onto $\HH$ such that $W(x_1)=x_2$ and
$W(p_{1,m})=p_{2,m}$, $1\le m\le N$. Let
$p_{1,\infty}=W^{-1}(\infty)$ and $p_{2,\infty}=W(\infty)$. For
$j=1,2$, let $S_j\in(0,T_j]$ be the largest number such that for $0\le
t<S_j$,  $K_j(t)$ does not separate $p_{j,\infty}$ from $\infty$ in
$\HH$. Then $(W(K_1(t)),0\le t<S_1)$ has the same law as
$(K_2(t),0\le t<S_2)$ up to a time-change. A similar result holds
for the traces. \label{coordinate}
\end{Lemma}
{\bf Proof.} Here we only consider the generic cases. The proof of
the degenerate cases is similar. Let
$Q_j=\{\infty,p_{j,1},\dots,p_{j,N},p_{j,\infty}\}$, $j=1,2$. Then
$W(Q_1)=Q_2$, and $S_j$ is the maximum number in $(0,T_j]$ such that
for $0\le t<S_j$, $K_j(t)$ does not divide $Q_j$ in $\HH$. For $0\le
t<S_1$, since $K_1(t)$ does not divide $Q_1$ in $\HH$, so
$W(K_1(t))$ does not divide $Q_2$ in $\HH$. From Theorem 3 in
\cite{SW}, after a time-change, $(W(K_1(t)),0\le t<S_1)$ is a partial
chordal SLE$(\kappa;\vec{\rho})$ process started from
$(x_2;\vec{p}_2)$. We now suffice to show that this chordal Loewner
chain can not be further extended without dividing $Q_2$ in $\HH$.
If this is not true, then $\cup_{0\le t<S_1} W(K_1(t))$ does not
divide $Q_2$ in $\HH$. So $\cup_{0\le t<S_1} K_1(t)$ does not divide
$Q_1$ in $\HH$, which contradicts the choice of $S_1$. $\Box$

\vskip 3mm

Note that if $\kappa\in(0,4]$ then $S_{j}=T_{j}$, $j=1,2$, so we
conclude that $(W(K_{1}(t)),0\le t<T_{1})$ has the same
distribution as $(K_{2}(t),0\le t<T_{2})$ up to a time-change.
In general, by adding $\infty$
to be a force point with suitable value of force, we may always make
the sum of forces equals to $\kappa-6$, so the lemma can be applied.

\subsection{Strip SLE}
Strip SLE is studied independently in \cite{thesis} and
\cite{Dipolar} (where it is called dipolar SLE). For $h>0$, let
$\St_h=\{z\in\C:0<\Imm z<h\}$ and $\R_h=ih+\R$. If $H$ is a bounded
closed subset of $\St_\pi$,  $\St_\pi\sem H$ is simply connected,
and has $\R_\pi$ as a boundary arc, then we call $H$ a hull in
$\St_\pi$ w.r.t.\ $\R_\pi$. For such $H$, there is a unique $\psi_H$
that maps $\St_\pi\sem H$ conformally onto $\St_\pi$, such that for
some $c\ge 0$, $\psi_H(z)=z\pm c +o(1)$ as $z\to\pm\infty$ in
$\St_\pi$. We call such $c$ the capacity of $H$ in $\St_\pi$ w.r.t.\
$\R_\pi$, and denote it by $\scap(H)$.

For $\xi\in C([0,T))$, the strip Loewner equation driven by $\xi$ is
\BGE \pa_t\psi(t,z)=\coth\big(\frac{\psi(t,z)-\xi(t)}2\big),\quad
\psi(0,z)=z.\label{strip}\EDE For $0\le t<T$, let $L(t)$ be the set
of $z\in\St_\pi$ such that the solution $\psi(s,z)$ blows up before
or at time $t$. We call $L(t)$ and $\psi(t,\cdot)$, $0\le t<T$,
strip Loewner hulls and maps, respectively, driven by $\xi$. It
turns out that $\psi(t,\cdot)=\psi_{L(t)}$ and $\scap(L(t))=t$ for
each $t$. From now on, we write $\coth_2(z)$, $\tanh_2(z)$, $\cosh_2(z)$, and $\sinh_2(z)$
for functions $\coth(z/2)$, $\tanh(z/2)$, $\cosh(z/2)$, and $\sinh(z/2)$, respectively.

Let $\kappa\ge 0$, $\rho_1,\dots,\rho_N\in\R$, $x\in\R$, and
$p_1,\dots,p_N\in \R\cup\R_\pi\cup\{+\infty,-\infty\}\sem\{x\}$. Let
$B(t)$ be a Brownian motion. Let $\xi(t)$ and $p_k(t)$, $1\le k\le
N$, be the solutions to the SDE:
\begin{equation}\left\{\begin{array}{lll} d\xi(t) & = &
\sqrt\kappa d B(t)+\sum_{k=1}^N\frac{\rho_k}2\coth_2({\xi(t)-p_k(t)})dt\\
dp_k(t) & = & \coth_2({p_k(t)-\xi(t)})dt,\quad 1\le k\le
N,\end{array}\right.\label{strip-kappa-rho}\end{equation} with
initial values $\xi(0)=x$ and $p_k(0)=p_k$, $1\le k\le N$. Here if
some $p_k=\pm\infty$ then $p_k(t)=\pm\infty$ and
$\coth_2({\xi(t)-p_k(t)})=\mp 1$ for all $t\ge 0$, so $p_k$ has
a constant effect on the equation. Suppose $[0,T)$ is the maximal
interval of the solution. Let $L(t)$, $0\le t<T$, be strip Loewner
hulls driven by $\xi$. Then we call $L(t)$ , $0\le t<T$, a (full)
strip SLE$(\kappa;\vec{\rho})$ process started from $(x;\vec{p})$,
where $\vec{\rho}=(\rho_1,\dots,\rho_N)$ and
$\vec{p}=(p_1,\dots,p_N)$.

The following two lemmas show that strip SLE$(\kappa;\vec{\rho})$
processes also satisfy conformal invariance, and are conformally
equivalent to the corresponding chordal SLE$(\kappa;\vec{\rho})$
processes. The proofs are similar to that of Lemma \ref{coordinate},
and use the result of Section 4 in \cite{SW}, so we omit the
proofs.

\begin{Lemma} Suppose $\kappa\ge 0$ and $\vec{\rho}=(\rho_1,\dots,\rho_N)$
with $\sum_{m=1}^N\rho_m =\kappa-6$. For $j=1,2$, let $L_j(t)$, $0\le
t<T_j$, be a strip SLE$(\kappa;\vec{\rho})$ process started from
$(x_{j};\vec{p}_{j})$, where $\vec{p}_{j}=(p_{j,1},\dots,p_{j,N})$.
Suppose $W$ is a conformal
or conjugate conformal map from $\St_\pi$ onto $\St_\pi$ such
that $W(x_1)=x_2$ and $W(p_{1,m})=p_{2,m}$, $1\le m\le N$. Let
$I_1=W^{-1}(\R_\pi)$ and $I_{2}=W(\R_\pi)$. For $j=1,2$, let
$S_j\in(0,T_j]$ be the largest number such that for $0\le t<S_j$,  $L_j(t)$
does not separate $I_j$ from $\R_\pi$ in $\St_\pi$. Then
$(W(L_1(t)),0\le t<S_1)$ has the same law as $(L_2(t),0\le t<S_2)$
up to a time-change. \label{coordinate*}
\end{Lemma}

\begin{Lemma} Suppose $\kappa\ge 0$ and $\vec{\rho}=(\rho_1,\dots,\rho_N)$
with $\sum_{m=1}^N\rho_m =\kappa-6$.   Let $K(t)$, $0\le t<T$, be a
chordal SLE$(\kappa;\vec{\rho})$ process started from
$(x;\vec{p})$, where $\vec{p}=(p_{1},\dots,p_{N})$. Let $L(t)$, $0\le t<S$, be a strip
SLE$(\kappa;\vec{\rho})$ process started from $(y;\vec{q})$, where
$\vec{q}=(q_{1},\dots,q_{N})$. Suppose $W$ is a conformal
or conjugate conformal map from $\HH$ onto $\St_\pi$ such that $W(x)=y$ and
$W(p_k)=q_k$, $1\le k\le N$. Let $I=W^{-1}(\R_\pi)$ and
$q_\infty=W(\infty)$. Let $T'\in(0,T]$ be the largest number such that for
$0\le t<T'$, $K(t)$ does not separate $I$ from $\infty$ in $\HH$.
Let $S'\in(0,S]$ be the largest number such that for $0\le t<S'$, $L(t)$
does not separate $q_\infty$ from $\R_\pi$. Then $(W(K(t)),0\le
t<T')$ has the same law as $(L(t),0\le t<S')$ up to a time-change.
\label{coordinate2}
\end{Lemma}

As usual, if $\kappa\in[0,4]$, then $S_j=T_j$, $j=1,2$, in Lemma
\ref{coordinate*}, and $T'=T$ and $S'=S$ in Lemma \ref{coordinate2}.
In general, for a strip SLE$(\kappa;\vec{\rho})$ process, by adding
$+\infty$ and $-\infty$ to be force points with suitable values of
forces, we may always make the sum of forces equals to $\kappa-6$, so the above two
lemmas can be applied. From Lemma \ref{coordinate2}, we have the
existence of the strip SLE$(\kappa;\vec{\rho})$ trace, and the above
two lemmas also hold for traces.

\section{Geometric Properties} \label{geom}
Suppose $\beta(t)$, $0\le t<T$, is a chordal
SLE$(\kappa;\vec{\rho})$ trace. In this section, we will study the
existence and property of the limit or subsequential limit of $\beta(t)$
as $t\to T$ in certain cases.
The three lemmas in the last section
will be frequently used.

\subsection{Many force points}
Let $\kappa>0$, $\vec{p}_\pm=(p_{\pm 1},\dots,p_{\pm N_\pm})$,
$\vec{\rho}_\pm=(\rho_{\pm 1},\dots,\rho_{\pm N_\pm})$, where
$0<p_1<\dots<p_{N_+}$,
$0>p_{-1}>\dots>p_{-N-}$, and $\rho_{\pm j}\in\R$, $j=1,\dots,N_\pm$.
Let
$\beta(t)$, $0\le t<T$, be a chordal
SLE$(\kappa;\vec{\rho}_+,\vec{\rho}_-)$ trace started from
$(0;\vec{p}_+,\vec{p}_-)$. Let $\vphi(t,\cdot)$ and  $\xi(t)$, $0\le
t<T$, be the chordal Loewner maps and driving function,
respectively, for the trace $\gamma$.

\begin{Theorem} Suppose for any $1\le k\le N_\pm$,
$\sum_{j=1}^k\rho_{\pm j}\ge \kappa/2-2$. \\ (i) Almost surely $T=\infty$,
so $\infty$ is a subsequential limit of $\beta(t)$ as
$t\to T$.
\\(ii) If in addition, $\kappa\in(0,4]$, then almost
surely $\lin{\beta((0,\infty))}\cap(\R\sem\{0\})=\emptyset$.
\label{more force} \end{Theorem} {\bf Proof.} Let
$\rho_\infty=\kappa-6-\sum_{j=1}^{N_+}\rho_j-\sum_{j=1}^{N_-}\rho_{-j}$.
Let $\chi_0=0$. For $1\le k\le N_\pm$, let $\chi_{\pm
k}=\sum_{j=1}^k\rho_{\pm j}\ge \kappa/2-2$. Let
$\chi^\pm_{\max}=\max\{\chi_{\pm k}:1\le k\le N_\pm\}$. \vskip 3mm

\no (i) If $T=\infty$, then the diameter of $\beta((0,t])$ tends to
$\infty$ as $t\to \infty$, so $\infty$ is a subsequential limit of
$\beta(t)$ as $t\to T$. So we suffice to prove that $T=\infty$ a.s..
If $T<\infty$, then for $x=p_1$ or $p_{-1}$, $\vphi(t,x)-\xi(t)\to
0$ as $t\to T$, where $\xi(t)$ and $\vphi(t,\cdot)$ are the driving
function and chordal Loewner map. For any $t\in[0,T)$,
$\vphi(t,p_{-1})<\xi(t)<\vphi(t,p_1)$, so
$\pa_t\vphi(t,p_1)=2/(\vphi(t,p_1)-\xi(t))>0$ and
$\pa_t\vphi(t,p_{-1})=2/(\vphi(t,p_{-1})-\xi(t))<0$. Thus
$\vphi(t,p_1)-\vphi(t,p_{-1})$ increases. If $\vphi(t,p_1)-\xi(t)\to
0$, then $(\vphi(t,p_1)-\xi(t))/(\vphi(t,p_1)-\vphi(t,p_{-1}))
\to 0$, so $(\vphi(t,p_1)-\xi(t))/(\xi(t)-\vphi(t,p_{-1}))
\to 0$. Similarly, if $\xi(t)-\vphi(t,p_{-1})\to
0$, then $(\xi(t)-\vphi(t,p_{-1}))/(\vphi(t,p_1)-\xi(t))
\to 0$. Thus if $T<\infty$, then
$\ln(\xi(t)-\vphi(t,p_{-1}))-\ln(\vphi(t,p_1)-\xi(t))$ tends to
$+\infty$ or $-\infty$ as $t\to T$.

Suppose $W$ maps $\HH$ conformally onto $\St_\pi$ such that $W(0)=0$
and $W(p_{\pm 1})=\pm\infty$. Let $q_\infty=W(\infty)$ and $q_{\pm
j}=W(p_{\pm j})$, $1\le j\le N_\pm$. Let $\gamma(t)=\beta(u^{-1}(t))$
for $0\le t<S=u(T)$, where $u$ is a continuous increasing function
on $[0,T)$ such that $\scap(L(t))=t$ for any $t$, and $L(t)$ is
the hull in $\St_\pi$ w.r.t.\ $\R_\pi$ generated by $\gamma((0,t])$.
>From Lemma \ref{coordinate2}, $\gamma(t)$, $0\le t<S$, is a strip
 SLE$(\kappa;\rho_\infty,\vec{\rho}_+,\vec{\rho}_-)$
 trace started from
 $(0;q_\infty,\vec{q}_+,\vec{q}_-)$, where $\vec{q}_\pm=(q_{\pm 1},\dots,q_{\pm N_\pm})$.
 Since all $q_j$'s are either $\pm\infty$ or lie on $\R_\pi$, which will
 never be swallowed, so $S=\infty$.

Let $\psi(t,\cdot)$ and $\eta(t)$, $0\le t<\infty$, be the strip
Loewner maps and driving function, respectively, for the trace
$\gamma$. Let $X_\infty(t)=\Ree\psi(t,q_\infty)-\eta(t)$ and $X_{\pm
j}(t)=\Ree\psi(t,q_{\pm j})-\eta(t)$, $1\le j\le N_\pm$. Then
$X_{-2}(t)<\dots<X_{-N_-}(t)<X_\infty(t)<X_{N_+}(t)<\dots <X_2(t)$.
And for some Brownian motion $B(t)$, $\eta(t)$ satisfies the SDE:
$$d\eta(t)=\sqrt\kappa
dB(t)-\frac{\rho_\infty}2
\tanh_2({X_\infty(t)})dt$$
$$-\sum_{j=1}^{N_+}\frac{\rho_j}2 \tanh_2({X_j(t)})dt-
\sum_{j=1}^{N_-}\frac{\rho_{-j}}2
\tanh_2({X_{-j}(t)})dt.$$ For $0\le t<T$, let
$W_t=\psi(u(t),\cdot)\circ W\circ\vphi(t,\cdot)^{-1}$. Then $W_t$
maps $\HH$ conformally onto $\St_\pi$,  $W_t(\xi(t))=\eta(u(t))$,
$W_t(\infty)=\psi(u(t),q_\infty)$, and $W_t(\vphi(t,p_{\pm
1}))=\pm\infty$. Thus
$$\ln(\xi(t)-\vphi(t,p_{-1}))-\ln(\vphi(t,p_1)-\xi(t))=
\Ree\psi(u(t),q_\infty)-\eta(u(t))=X_\infty(u(t)).$$ Thus if
$T<\infty$ then $X_\infty(t)$ tends to $+\infty$ or $-\infty$ as
$t\to \infty$. So now we suffice to show that a.s.\
$\limsup_{t\to\infty} X_\infty(t)=+\infty$ and $\liminf_{t\to\infty}
X_\infty(t)=-\infty$. We will prove that a.s.\ $\limsup_{t\to\infty}
X_\infty(t)=+\infty$. The other statement follows from symmetry.

Let
$X_{N_++1}(t)=X_{-N_--1}(t)=X_\infty(t)$. Then $X_\infty(t)$
satisfies the SDE:
$$dX_\infty(t)=-\sqrt\kappa
dB(t)+(\frac\kappa 2-2-\frac{\chi_{N_+}+\chi_{-N_-}}2)
\tanh_2({X_\infty(t)})dt$$
$$+\sum_{j=1}^{N_+}\frac{\chi_j-\chi_{j-1}}2 \tanh_2({X_j(t)})dt+
\sum_{j=1}^{N_-}\frac{\chi_{-j}-\chi_{-j+1}}2
\tanh_2({X_{-j}(t)})dt$$
$$=-\sqrt\kappa dB(t)+(\frac\kappa
2-2)\tanh_2({X_\infty(t)})dt$$$$
+\sum_{j=1}^{N_+}\frac{\chi_j}2
\big(\tanh_2({X_j(t)})-\tanh_2({X_{j+1}(t)})\big)dt$$
$$ +\sum_{j=1}^{N_-}\frac{\chi_{-j}}2
\big(\tanh_2({X_{-j}(t)})-\tanh_2({X_{-j-1}(t)})\big)dt.$$
Note that for $1\le j\le N_{\pm}$, $\pm (\tanh_2(X_{\pm
j}(t))-\tanh_2(X_{\pm j\pm 1}(t)))>0$. Since $\kappa/2-2\le
\chi_j\le \chi^\pm_{\max}$ for $1\le j\le N_{\pm}$, so for
some adapted process $A(t)\ge 0$,
$$dX_\infty(t)= -\sqrt\kappa dB(t)+A(t)dt+(\frac\kappa
2-2)\tanh_2({X_\infty(t)})dt$$$$+(\frac\kappa
4-1)\sum_{j=1}^{N_+}
\big(\tanh_2({X_j(t)})-\tanh_2({X_{j+1}(t)})\big)dt$$
$$+\frac{\chi^-_{\max}}2 \sum_{j=1}^{N_-}
\big(\tanh_2({X_{-j}(t)})-\tanh_2({X_{-j-1}(t)})\big)dt$$
$$=-\sqrt\kappa dB(t)+A(t)dt+(\frac\kappa
4-1-\frac{\chi^-_{\max}}2)\big(1+
\tanh_2({X_\infty(t)})\big)dt.$$ Note that $\tanh_2(X_{\pm 1}(t))=\pm 1$.
 Define $f$ on $\R$ such that for any
$x\in \R$, $f'(x)=(e^x+1)^{\frac
2\kappa({\chi^-_{\max}}+2-\kappa/2)}$. Since $\chi^-_{\max}\ge
\kappa/2-2$, so $f'(x)\ge 1$ for any $x\in\R$. Thus
  $f$ maps $\R$ onto $\R$. Let $Y(t)=f(X_\infty(t))$,
and $\til A(t)=f'(X_\infty(t))A(t)\ge 0$. From Ito's formula, we
have
$$dY(t)=-\sqrt\kappa f'(X_\infty(t))dB(t)+\til A(t)dt.$$
Let $M(t)=Y(t)-\int_0^t \til A(s)ds$. Then $Y(t)\ge M(t)$ and
$dM(t)=-\sqrt\kappa f'(X_\infty(t))dB(t)$. Let $v(t)=\int_0^t\kappa
f'(X_\infty(s))^2ds$. Then $v$ is a continuous increasing function
on $[0,\infty)$, and maps $[0,\infty)$ onto $[0,\infty)$. And
$M(v^{-1}(t))$, $0\le t<\infty$, is a Brownian motion. Thus a.s.\
$\limsup_{t\to\infty} M(t)=+\infty$. Since $Y(t)\ge M(t)$ for any
$t$, so a.s.\ $\limsup_{t\to\infty} Y(t)=+\infty$. Since
$X_\infty(t)=f^{-1}(Y(t))$, so a.s.\ $\limsup_{t\to\infty}
X_\infty(t)=+\infty$, as desired.

\vskip 3mm

\no (ii) From symmetry, we suffice to show that a.s.\
$\lin{\beta((0,\infty))}\cap(-\infty,0)=\emptyset$. Fix any $r_+\in
(-\infty,p_{-N_-})\cap\Q$ and $r_-\in (p_{-1},0)\cap\Q$. We suffice
to show that a.s.\ $\lin{\beta((0,\infty))}\cap(r_+,r_-)=\emptyset$.
Choose $W$ that maps $\HH$ conformally onto $\St_\pi$ such that
$W(0)=0$ and $W(r_\pm)=\pm\infty$. Let $q_{\pm j}=W(p_{\pm j})$,
$1\le j\le N_\pm$, and
$\vec{q}_\pm=(q_{\pm 1},\dots,q_{\pm N_\pm})$. Let
$q_\infty=W(\infty)\in (0,\infty)$. Let $\gamma(t)=\beta(u^{-1}(t))$
for $0\le t<S=u(T)$, where $u$ is a continuous increasing function
on $[0,T)$ such that $\scap(\gamma((0,t]))=t$ for any $t\in[0,S)$.
>From Lemma \ref{coordinate2}, $\gamma(t)$, $0\le t<S$, is a strip
 SLE$(\kappa;\rho_\infty,\vec{\rho}_+,\vec{\rho}_-)$
 trace started from
 $(0;q_\infty,\vec{q}_+,\vec{q}_-)$.

Let $\psi(t,\cdot)$ and $\eta(t)$, $0\le t<S$, be the strip Loewner
maps and driving function, respectively, for the trace $\gamma$. Let
$X_\infty(t)=\psi(t,q_\infty)-\eta(t)$,
$q_{N_++1}=q_{-N_--1}=q_\infty$, and $X_{\pm j}(t)=\psi(t,q_{\pm
j})-\eta(t)$, $1\le j\le N_{\pm}+1$. Then there is a Brownian motion
$B(t)$ such that $X_\infty(t)$ satisfies:
$$dX_\infty(t)=-\sqrt\kappa
dB(t)+\big(1+\frac{\rho_\infty}2\big)
\coth_2({X_\infty(t)})dt$$
$$+\sum_{j=1}^{N_+}\frac{\rho_j}2 \coth_2({X_j(t)})dt+
\sum_{j=1}^{N_-}\frac{\rho_{-j}}2
\coth_2({X_{-j}(t)})dt$$
$$=-\sqrt\kappa
dB(t)+\big(\frac\kappa 2-2\big)
\coth_2({X_\infty(t)})dt$$$$
+\sum_{j=1}^{N_+}\frac{\chi_j}2
\big(\coth_2({X_j(t)})-\coth_2({X_{j+1}(t)})\big)dt$$
$$ +\sum_{j=1}^{N_-}\frac{\chi_{-j}}2
\big(\coth_2({X_{-j}(t)})-\coth_2({X_{-j-1}(t)})\big)dt.$$
Since $X_j(t)$, $1\le j\le N_++1$, lie on the boundary of $\St_\pi$
in the counterclockwise direction; and $X_{-j}(t)$, $1\le j\le
N_-+1$, lie on the boundary of $\St_\pi$ in the clockwise direction,
so we have $\pm(\coth_2(X_{\pm j}(t))-\coth_2(X_{\pm(j+1)}))>0$ for
$1\le j\le N_{\pm}$. Since $\chi_{-j}\ge \kappa/2-2$, $1\le j\le
N_-$, and $\chi_j\le \chi^+_{\max}$, $1\le j\le N_+$, so for
some adapted process $A_1(t)\ge 0$,
$$dX_\infty(t)=-\sqrt\kappa
dB(t)-A_1(t)dt+\big(\frac\kappa 2-2\big)
\coth_2({X_\infty(t)})dt$$
$$+(\frac\kappa
4-1)\big(\coth_2({X_{-1}(t)})
-\coth_2({X_{-N_--1}(t)})\big)dt$$$$+\frac{\chi^+_{\max}}2
\big(\coth_2({X_{1}(t)})
-\coth_2({X_{N_++1}(t)})\big)dt$$
$$=-\sqrt\kappa
dB(t)-A_1(t)dt+(\frac\kappa
4-1)\big(\coth_2({X_{-1}(t)})
+\coth_2({X_{\infty}(t)})\big)dt$$
$$+\frac{\chi^+_{\max}}2
\big(\coth_2({X_{1}(t)})
-\coth_2({X_{\infty}(t)})\big)dt$$ Note that
$X_{-1}(t)\in\R_\pi$ and $X_\infty(t)\in(0,\infty)$, so
$\coth_2(X_{-1}(t))+\coth_2(X_{\infty}(t))>0$. Since
$\kappa\in(0,4]$, so $\kappa/4-1\le 0$. Thus for some adapted process
$A_2(t)\ge A_1(t)\ge 0$,
$$dX_\infty(t)=-\sqrt\kappa
dB(t)-A_2(t)dt+\frac{\chi^+_{\max}}2
\big(\coth_2({X_{1}(t)})
-\coth_2({X_{\infty}(t)})\big)dt.$$ For $0\le t<S$,
since $X_\infty(t)>0$, so
$$\sqrt\kappa B(t)\le\frac{\chi^+_{\max}}2\int_0^t
\big(\coth_2({X_{1}(s)})
-\coth_2({X_{\infty}(s)})\big)ds.$$ Since
$0<X_1(s)<X_\infty(s)$ for $0\le s<S$, so the integrand is positive.
Thus if $\chi^+_{\max}\le 0$, then $B(t)\le 0$ for $0\le t<S$. Now
suppose $\chi^+_{\max}> 0$. Let $q_1(t)=\psi(t,q_1)$ and
$q_\infty(t)=\psi(t,q_\infty)$. From the strip Loewner equation, for
$0\le t<S$,
$$\sqrt\kappa B(t)\le\frac{\chi^+_{\max}}2(q_1(s)-q_\infty(s))|^{s=t}_{s=0}\le
-\frac{\chi^+_{\max}}2(q_1(s)-q_\infty(s))|_{s=0}=
\frac{\chi^+_{\max}}2(q_\infty-q_1),$$ where the second ``$\le$''
follows from the fact that $q_1(t)<q_\infty(t)$. Thus in any case,
$B(t)$ is uniformly bounded above on $[0,S)$. So we have $S<\infty$
a.s..

For a hull $H$ in $\St_\pi$ w.r.t.\ $\R_\pi$, if $\scap(H)=s$ then
the height of $H$ is no more than $2\cos^{-1}(e^{-s/2})$, and the
equality is attained when $H$ is some vertical line segment. Now for
$0\le t<S$, $\scap(\gamma((0,t]))=t<S$, so the distance between
$\gamma((0,t])$ and $\R_\pi$ is bigger than
$\pi-2\cos^{-1}(e^{-S/2})$. Since a.s.\ $S<\infty$, so
$\gamma((0,S))$ is bounded away from $\R_\pi$. From the property of
$W$ and the definition of $\gamma$, we conclude that a.s.\
$\beta((0,\infty))$ is bounded away from $(r_+,r_-)$. So we are
done. $\Box$

\begin{Theorem} Suppose $x\in\R$, $\kappa\in(0,4]$,
$\rho_1,\rho_2\ge\kappa/2-2$, and $\beta(t)$, $0\le t<\infty$, is a
chordal SLE$(\kappa;\rho_1,\rho_2)$ trace started from
$(x;p_1,p_2)$.\\
(i) If $p_1=x^-$ and $p_2=x^+$, then a.s.\ $\lim_{t\to \infty}\beta(t)=\infty$.\\
(ii) If $p_1\in(-\infty,x)$ and $p_2\in(x,+\infty)$, then a.s.\
$\lim_{t\to \infty}\beta(t)=\infty$.
 \label{lim=infty}
\end{Theorem}
{\bf Proof.} We may suppose $x=0$. We first consider the case that
$p_1=x^-=0^-$ and $p_2=x^+=0^+$. Let $Z$ denote the set of
subsequential limits in $\lin{\HH}$ of $\beta(t)$ as $t\to \infty$.
We suffice to show that $Z=\emptyset$ a.s.. From Lemma
\ref{coordinate}, for any $a>0$, $a^2\beta(t)$, $0\le t<\infty$, has
the same distribution as $\beta(at)$, $0\le t<\infty$, which implies
that $a^2Z$ has the same distribution as $Z$. Thus we suffice to
show that a.s.\ $0\not\in Z$.

Let $\vphi(t,\cdot)$ and $\xi(t)$ be the chordal
Loewner maps and driving function for the trace $\beta$.
Choose $W_t$ that maps $(\HH;\xi(t),\vphi(t,0^+),\vphi(t,0^-))$
conformally onto $(\St_\pi;0,+\infty,-\infty)$, and let $X_\infty(t)=\Ree W_t(\infty)$.
Then $X_\infty(t)=\ln(\vphi(t,0^+)-\xi(t))-\ln(\xi(t)-\vphi(t,0^-))$. From
the proof of Theorem \ref{more force} (i), we see that a.s.\ $\limsup X_\infty(t)=+\infty$
and $\liminf X_\infty(t)=-\infty$. Thus a.s.\ there is $t\ge 1$ such that $X_\infty(t)=0$, i.e.,
$\vphi(t,0^+)-\xi(t)=\xi(t)-\vphi(t,0^-)$. Let $T$ denote the first $t$ with this property.
So $T$ is a finite stopping time.

 Let $g(z)=(\vphi(T,z)-\xi(T))/(\vphi(T,0^+)-\xi(T))$
and $f=g^{-1}$. Then $g$ maps $\HH\sem\beta((0,T])$ conformally onto
$\HH$, $g(\beta(T))=0$; and $f$ extends continuously to $\HH\cup\R$
such that $f^{-1}(0)=\{-1,1\}$. Let $\gamma(t)=g(\beta(T+t))$, $t\ge
0$. Then after a time-change, $\gamma(t)$, $0\le t<\infty$, has the
same distribution as a chordal SLE$(\kappa;\rho_1,\rho_2)$ trace
started from $(0;-1,1)$. From Theorem \ref{more force} (ii),
$\gamma((0,\infty))$ is bounded away from $\{-1,1\}$ a.s.. Thus
a.s.\ $\beta([T,\infty))$ is bounded away from $0$, which implies
that  $0\not\in Z$. So we proved (i).

(ii) Suppose $p_1\in(-\infty,x)$ and $p_2\in(x,\infty)$. Let
$r=(p_2-x)/(x-p_1)$. Let $\beta_0(t)$ be a chordal
SLE$(\kappa;\rho_1,\rho_2)$ trace started from $(0;0^-,0^+)$. Let
$\vphi(t,\cdot)$ and $\xi(t)$, $0\le t<\infty$, be the chordal
Loewner maps and driving function for the trace $\beta_0$. Let $X_\infty(t)$
be defined as in the last paragraph with $\beta$ replaced by $\beta_0$. Then
there is
a.s.\ $t\ge 1$ such that $X_\infty(t)=\ln(r)$, i.e.,
$(\vphi(t,0^+)-\xi(t))/(\xi(t)-\vphi(t,0^-))=r$. Let $T_r$ denote this
time. Since $(X_\infty(t))$ is recurrent, $T$ is a finite stopping time. Let
$$g(z)=x+\frac{(p_2-p_1)(\vphi(T_r,z)-\xi(T_r))}{\vphi(T_r,0^+)-\vphi(T_r,0^-)}.$$
Then $g$ maps $(\HH\sem\beta_0((0,T_r]);\beta_0(T_r),0^-,0^+)$ conformally onto $(\HH;x,p_1,p_2)$.
So after a time-change, $(g(\beta_0(T_r+t)),0\le t<\infty)$, has the same
distribution as $(\beta(t),0\le t<\infty)$. From (i), a.s.\
$\lim_{t\to\infty}\beta_0(t)=\infty$, so we have a.s.\
$\lim_{t\to\infty}\beta(t)=\infty$. $\Box$

\begin{Conjecture} (Reversibility)  Suppose $\kappa\in(0,4)$,
$\rho_-,\rho_+\ge\kappa/2-2$, and $\beta(t)$, $0\le t<\infty$, is a
chordal SLE$(\kappa;\rho_-,\rho_+)$ trace started from
$(0;0^-,0^+)$. Let $W(z)=1/\lin{z}$. Then after a time-change, the
reversal of $(W(\beta(t)))$ has the same distribution as $(\beta(t))$.\label{conjec}
\end{Conjecture}

If $\kappa=0$, the conjecture is trivial because the trace is a half line. If
$\rho_+=\rho_-=0$, i.e., $\beta$ is a standard chordal SLE$(\kappa)$ trace, the
reversibility is known in \cite{reversibility}. If $\kappa=4$, the reversibility is
a result of the convergence of discrete Gaussian free field contour line in \cite{SS};
and is also a special case of Theorem \ref{kappa=4} in this paper.
To prove this conjecture using the technique in \cite{reversibility} and this paper, one may
need to know the conditional distribution of $\beta(t)$, $T_1\le t<T_2$, given its
initial segment $\beta([0,T_1])$  and final segment $\beta([T_2,\infty))$,
where $T_1$ is a stopping time, $T_2$ is a ``backward'' stopping time, and $T_1<T_2$.
In the case that $\beta$ is a standard chordal SLE$(\kappa)$ trace, we find that
$\beta(t)$, $T_1\le t<T_2$, is a chordal SLE$(\kappa)$ trace in
$\HH\sem(\beta((0,T_1]\cup[T_2,\infty)))$ from $\beta_1(T_1)$ to $\beta_2(T_2)$, up to
a time-change. If $\kappa=4$, we will see in the proof of Theorem \ref{kappa=4}
 that after a time-change, $\beta(t)$, $T_1\le t<T_2$, is a generic
SLE$(\kappa;\rho_-,\rho_+)$ trace in $\HH\sem(\beta((0,T_1]\cup[T_2,\infty)))$
In general, this conditional distribution may not be an
SLE$(\kappa;\vec{\rho})$ trace.

\subsection{Two force points} We now study a strip SLE
process with two force points at $\infty$ and $-\infty$. Let $\kappa>0$ and
$\rho_+,\rho_-\in\R$. Suppose $\beta(t)$, $0\le t<T$, is a strip
SLE$(\kappa;\rho_+,\rho_-)$ trace started from
$(0;+\infty,-\infty)$. Let $\sigma=(\rho_--\rho_+)/2$. Then
$T=\infty$ and the driving function is
$\xi(t)=\sqrt{\kappa}B(t)+\sigma t$, $0\le t<\infty$, for some
Brownian motion $B(t)$. Let $L(t)$ and $\psi(t,\cdot)$, $0\le
t<\infty$, be the strip Loewner hulls and maps, respectively, driven
by $\xi$.

We first consider the case that $|\sigma|<1$. Then $\xi(t)$
satisfies
\begin{equation} |\xi(t)|\le A(\omega)+\sigma' t,\,\,\,\,\forall t\ge 0,\label{bound0}
\end{equation}where $\sigma':=(1+|\sigma|)/2<1$ and
$A(\omega)>0$ is a random number.

\begin{Lemma} If $|\sigma|<1$, then $L(\infty)$ is bounded. \label{bdd}
\end{Lemma}
{\bf Proof.} Let $\sigma''=(1+|\sigma'|)/2$. We may choose $R>0$ such
that $\Ree\coth_2(z)>\sigma''$ when $z\in\St_\pi$ and $\Ree z\ge R$.
>From (\ref{bound0}) there is $a=a(\omega)\ge R+1$ such that
$R+1+\xi(t)-\sigma'' t\le a$ for all $t\ge 0$. Consider a point
$z\in\St_\pi$ with $\Ree z\ge a$. Suppose there is $t$ such that
$\Ree\psi(t,z)-\xi(t)<R$. Since $\psi(0,z)=z$, so $\Ree
\psi(0,z)=\Ree z\ge a>R$. Since $\xi(0)=0$, so
$\Ree\psi(0,z)-\xi(0)\ge a>R$. Thus there is a first $t_0$ such that
$\Ree\psi(t_0,z)-\xi(t_0)=R$. For $t\in[0,t_0]$, we have
$\Ree\psi(t,z)-\xi(t)\ge R$, and so
$$\pa_t\Ree\psi(t,z)=\Ree\coth_2(\psi(t_0,z)-\xi(t_0))\ge \sigma''.$$
Integrating the above inequality w.r.t.\ $t$ from $0$ to $t_0$, we
get
$$R=\Ree\psi(t_0,z)-\xi(t_0)\ge\Ree\psi(0,z)+\sigma'' t_0-\xi(t_0)
\ge a+\sigma'' t_0-\xi(t_0)\ge R+1,$$ where the last inequality uses
the property of $a$. So we get a contradiction. Therefore
$\Ree\psi(t,z)-\xi(t)\ge R$ for all $t\ge 0$. So $\psi(t,z)$ will
never blow up, which means that $z\not\in L(t)$ for all $t\ge 0$,
and so $z\not\in L(\infty)$. Similarly, there is $a'=a'(\omega)>0$
such that if $z\in\St_\pi$ and $\Ree z\le -a'$ then $z\not\in
L(\infty)$. Thus $L(\infty)$ is contained in $\{x+iy:-a'< x<
a,0<y<\pi\}$, and so is bounded. $\Box$

\vskip 3mm

Let
\begin{equation}f_{\kappa,\sigma}(x)=\int_{-\infty}^x
\exp(s/2)^{\frac{4\sigma}\kappa}\cosh_2(
s)^{-\frac{4}{\kappa}}ds. \label{f}
\end{equation}
 Since $|\sigma|<1$, so $f_{\kappa,\sigma}$ maps $\R$
onto the interval $(0,A_{\kappa,\sigma})$ for some
$A_{\kappa,\sigma}<\infty$.

Let $$X_t(z)=\Ree\psi(t,z)-\xi(t).$$ Now fix $z_0=x_0+\pi
i\in\R_\pi$. Then $\psi(t,z_0)\in\R_\pi$ for all $t$. Let $X_t$
denote $X_t(z_0)$ temporarily. Then $dX_t=\tanh_2(X_t)dt-d\xi(t)$.
>From Ito's formula, we have
$$df_{\kappa,\sigma}(X_t)=-\exp({
X_t}/{2})^{\frac{4\sigma}\kappa}\cosh_2(
{X_t})^{-\frac{4}{\kappa}}\sqrt{\kappa}dB(t).$$ Thus
$f_{\kappa,\sigma}(X_t)$ is a local martingale.

Let $u(0)=0$ and $u'(t)=[\exp({
X_s}/{2})^{\frac{4\sigma}\kappa}\cosh_2(
{X_s})^{-\frac{4}{\kappa}}\sqrt{\kappa}]^2$.
Then $u$ is a continuous increasing function. Let $T=u(\infty)\in(0,+\infty]$, and
$v=u^{-1}$. Then
$(f_{\kappa,\sigma}(X_{v(t)}),0\le t<T)$ has the same distribution
as $(B(t), 0\le t<T)$. Since $ f_{\kappa,\sigma}(X_{v(t)})$ stays
inside $(0,A_{\kappa,\sigma})$, so from the property of Brownian
motion, we have a.s.\ $T<\infty$ and $\lim_{t\to T}
f_{\kappa,\sigma}(X_{v(t)})$ exists. If  $\lim_{t\to T}
f_{\kappa,\sigma}(X_{v(t)})$ is neither $0$ nor $A_{\kappa,\sigma}$,
then $f_{\kappa,\sigma}(X_{v(t)})$ is uniformly bounded away from
$0$ and $A_{\kappa,\sigma}$ on $[0,T)$, so $X_t$ is uniformly
bounded on $[0,\infty)$, which implies that $u'(t)$ is uniformly
bounded below, and so $T=u(\infty)=\infty$.
Since $T<\infty$ a.s., so $\lim_{t\to T}
f_{\kappa,\sigma}(X_{v(t)})\in\{0,A_{\kappa,\sigma}\}$ a.s.. Thus
$\lim_{t\to\infty} X_t\in\{\pm \infty\}$ a.s.. Moreover, the
probability that $X_t\to +\infty$ is equal to
$f_{\kappa,\sigma}(x_0)/A_{\kappa,\sigma}$ by the Markov property.

Define
$$J_+=\inf\{x\in\R:\lim_{t\to\infty}X_t(x+\pi i)=+\infty\};$$
$$J_-=\sup\{x\in\R:\lim_{t\to\infty}X_t(x+\pi i)=-\infty\}.$$
Since $x_1<x_2$ implies $\Ree\psi(t,x_1+\pi i)<\Ree\psi(t,x_2+\pi
i)$ for all $t$, so we have $J_-\le J_+$; and for $x<J_-$,
$X_t(x+\pi i)\to -\infty$, for $x>J_+$, $X_t(x+\pi i)\to +\infty$ as
$t\to\infty$. Hence $\PP\{J_+<x\}\le
f_{\kappa,\sigma}(x)/A_{\kappa,\sigma}\le\PP\{J_-\le x\}$ for all
$x\in\R$. Since $f_{\kappa,\sigma}$ is strictly increasing, so
$J_-=J_+$ a.s.. By discarding an event of probability $0$, we may
assume that $J_+=J_-$, and let it be denoted by $J$. The density of
$J$ is
$\exp(x/2)^{\frac{4\sigma}\kappa}\cosh_2(x)^{-{4}/{\kappa}}/A_{\kappa,\sigma}$.

\begin{Lemma} $\lin{L(\infty)}\cap\R_\pi=\{J+\pi i\}$. \label{onept}
\end{Lemma}
{\bf Proof.} If $J+\pi i\not\in\lin{L(\infty)}$, then there are
$b,c>0$ such that $\dist(x+\pi i,L(\infty))>c$ for all
$x\in[J-b,J+b]$. From the definition of $J$, $X_t(J\pm b+\pi i)\to
\pm\infty$ as $t\to\infty$. Thus $\Ree\psi(t,J+b+\pi
i)-\Ree\psi(t,J-b+\pi i)\to+\infty$ as $t\to\infty$. By mean value
theorem, for each $t$, there is $x_t\in[J-b,J+b]$ such that
$|\pa_z\psi(t,x_t+\pi i)|\to\infty$ as $t\to\infty$. From Koebe's
$1/4$ theorem, we conclude that $\dist(x_t+\pi i, L(t))\to 0$, which
is a contradiction. Thus $J+\pi i\in\lin{L(\infty)}$.

Suppose $x_0>J$. Then $X_t(x_0+\pi i)\to+\infty$ as $t\to\infty$.
Thus $\pa_t\psi(t,x_0+\pi i)\to 1$ as $t\to\infty$. Recall that
$0<\sigma'<1$, and $|\xi(t)|\le A(\omega)+\sigma' t$ for all $t\ge
0$. So there is $H>0$ such that when $t\ge H$, $X_t(x_0+\pi
i)=\Ree\psi(t,x_0+\pi i)-\xi(t)
>\frac{1-\sigma'}2\,t$. So $X_t(x+\pi i)
>\frac{1-\sigma'}2\,t$ for all $x\ge x_0$ and $t\ge H$.

Differentiate equation (\ref{strip}) w.r.t.\ $z$, then we get
\begin{equation*}\pa_t\pa_z\psi(t,z)={-1/2 \cdot \pa_z\psi(t,z)}\cdot
{\sinh_2(\psi(t,z)-\xi(t))^{-2}}.
\end{equation*}
Thus \BGE\pa_t\ln|\pa_z\psi(t,z_0)|=\Ree(-1/2\cdot
{\sinh_2(\psi(t,z_0)-\xi(t))^{-2}}).\label{diffstrip}\EDE It
follows that for all $x\ge x_0$,
$$|\pa_z\psi(t,x+\pi i)|=\exp\left(\int_0^t\Ree \left(\frac{-1/2}
{\sinh_2(\psi(s,x+\pi i)-\xi(s))^2}\right)ds\right)$$
$$=\exp\left(\int_0^t\Ree \left(\frac{-1/2}
{\sinh_2(X_s(x+\pi i)+\pi i)^2}\right)ds\right)$$
$$=\exp\left(\int_0^t\frac{1/2}{\cosh_2(X_s(x+\pi i))^2}ds\right).$$
$$\le\exp\left(\int_0^H \frac{ds}2+\int_H^\infty
\frac{1}{2\cosh_2(\frac{1-\sigma'}2s)^2}ds\right)<+\infty.$$ Then by Koebe's
$1/4$ theorem, for all $x\ge x_0$, $x+\pi i$ is bounded away from
$L(\infty)$ uniformly. Thus $\lin{L(\infty)}$ is disjoint from
$[x_0+\pi i,+\infty)$ for all $x_0>J$. So $\lin{L(\infty)}$ is
disjoint from $(J+\pi i,+\infty)$. Similarly, $\lin{L(\infty)}$ is
disjoint from $(-\infty,J+\pi i)$. Thus $\lin{L(\infty)}$ intersects
$\R_\pi$ only at $J+\pi i$. $\Box$

\begin{Theorem} If $\kappa\in(0,4]$ and $|\sigma|<1$, then a.s.\
$\lim_{t\to\infty}\beta(t)\in\R_\pi$. \label{lim-1p}
\end{Theorem}
{\bf Proof}. Let $Q=\cap_{0\le t<\infty}\lin{\beta[t,\infty)}$. By
Lemma \ref{bdd}, $Q$ is nonempty and compact. Suppose $\til\xi$ has
the same law as $\xi$, and is independent of $\xi$. Let
$\til\beta(t)$ and $\til\psi(t,\cdot)$, $0\le t<\infty$, be the
strip Loewner trace and maps driven by $\til\xi$, respectively. Let
$(\til\F_t)$ be the filtration generated by $\til\xi$. For
$h\in(0,1)$, let $T_h$ be the first $t$ such that
$\Imm\til\beta(t)=\pi-h$. From Lemma \ref{onept}, $T_h$ is a finite
$(\til\F_t)$-stopping time. Let $\xi_*(t)=\til\xi(t)$ for $0\le t\le
T_h$; $\xi_*(t)=\til\xi(T_h)+\xi(t-T_h)$ for $t\ge T_h$. Then
$\xi_*$ has the same distribution as $\xi$. Let $\beta_*(t)$ be the
strip Loewner trace driven by $\xi_*$. Then
$\beta_*(t)=W_{T_h}(\beta(t-T_h))$ for $t\ge T_h$, where
$W_{T_h}(z):=\til\psi(T_h,\cdot)^{-1}(\til\xi(T_h)+z)$. Since
$\beta_*$ has the same distribution as $\beta$, so $W_{T_h}(Q)$ has
the same law as $Q$.

 Let $\Lambda_-$ denote the set of curves in
$\St_\pi\sem\til\beta((0,T_h])$ that connecting $(-\infty,0)$ with
the union of $[0,\infty)$ and the righthand side of
$\til\beta((0,T_h])$. Let $p=\Ree\til\beta(T_h)+\pi i$, and
$A=\{z\in\St_\pi:h<|z-p|<\pi\}$. Then every curve in $\Lambda_-$
crosses $A$. Thus the extremal length (\cite{Ahl}) of $\Lambda_-$ is
at least $(\ln(\pi)-\ln(h))/\pi$. From the property of $\psi_{T_h}$,
$W_{T_h}$ maps $\St_\pi$ conformally onto
$\St_\pi\sem\til\beta((0,T_h])$. There are $c_h<0<d_h$ such that
$W_{T_h}((-\infty,c_h])=(-\infty,0]$ and
$W_{T_h}([d_h,\infty))=[0,\infty)$. Since $\kappa\in(0,4]$, so
$W_{T_h}((c_h,d_h)=\beta((0,T_h])\subset\St_\pi$. Moreover,
$W_{T_h}(0)=\til\beta(T_h)$, and $W_{T_h}$ maps $[0,d_h)$ to the
righthand side of $\til\beta((0,T_h])$. From conformal invariance of
extremal length, the extremal distance between $(-\infty,c_h)$ and
$[0,\infty)$ in $\St_\pi$ is not less than $(\ln(\pi)-\ln(h))/\pi$.
Thus $c_h\to -\infty$ uniformly as $h\to 0$. Similarly,
$d_h\to+\infty$ uniformly as $h\to 0$.

For any $z\in\lin\St_\pi$, we have $\Imm W_{T_h}(z)\ge \Imm z$; and
the strict inequality holds when $z\in\St_\pi$ or $z\in(c_h,d_h)$.
Thus $\min\{\Imm W_{T_h} (Q)\}\ge \min\{\Imm Q\}$. Since
$W_{T_h}(Q)$ has the same law as $Q$, so a.s.\ $\min\{\Imm W_{T_h}
(Q)\}=\min\{\Imm Q\}$. Suppose now $Q\not\subset \R_\pi$ holds with
a positive probability. Since $Q$ is a bounded set, there is $R>0$
such that $\PP[{\cal E}_R]>0$, where ${\cal E}_R$ denotes the event
that $Q\subset\{z:|\Ree z|<R\}$ and $Q\not\subset\R_\pi$ both hold.
If $h$ is small enough, we have $|c_h|,|d_h|>R$. Assume that ${\cal
E}_R$ occurs. For any $z\in Q\sem\R_\pi$, either $z\in\St_\pi$ or
$z\in (c_h,d_h)$. In both cases, we have $\Imm W_{T_h}(z)>\Imm z$.
Thus $\min\{\Imm W_{T_h} (Q)\}> \min\{\Imm Q\}$ on ${\cal E}_R$,
which is a contradiction. Thus a.s.\ $Q\subset\R_\pi$. From Lemma
\ref{onept}, we have a.s.\ $Q=\{J+\pi i\}$, which means that
$\lim_{t\to\infty}\beta(t)=J+\pi i$. $\Box$

\vskip 3mm

Now we consider the case that $|\sigma|\ge 1$.

\begin{Theorem} If $\kappa\in(0,4]$ and $\pm\sigma\ge 1$, then almost surely
$\lim_{t\to\infty}\beta(t)=\pm\infty$. \label{lim-1p*}
\end{Theorem}
{\bf Proof.} Let $\sigma\ge 1$. Let $W(z)=e^z-1$. Then $W$ maps
$(\St_\pi;0,+\infty,-\infty)$ conformally onto $(\HH;0,\infty,-1)$. From Lemma
\ref{coordinate2}, after a time-change, $W(\beta(t))$, $0\le
t<\infty$, has the same distribution as a chordal
SLE$(\kappa;\frac\kappa2-3+\sigma)$ trace started from $(0;-1)$,
which is also a chordal SLE$(\kappa;\frac\kappa2-3+\sigma,0)$ trace
started from $(0;-1,1)$. Since $\sigma\ge 1$, so
$\frac\kappa2-3+\sigma\ge \frac\kappa2-2$. Since $\kappa\in(0,4]$, so $0\ge
\frac\kappa2-2$. Thus from Theorem \ref{lim=infty} (ii), a.s.\
$\lim_{t\to\infty} W(\beta(t))=\infty$, which implies that
$\lim_{t\to\infty}\beta(t)=+\infty$. The case $\sigma\le -1$ is
similar. $\Box$

\vskip 3mm

\no{\bf Remark.} Theorem \ref{lim-1p} and Theorem \ref{lim-1p*} should hold true in
the case $\kappa>4$. For example, the only part that the condition
$\kappa\in(0,4]$ is used in the proof of Theorem \ref{lim-1p} is
that $\Imm W_{T_h}(x)>0=\Imm x$ for $c_h<x<d_h$. If this is not true
for any $\kappa>4$, then we get some cut point of the hull that lies
on the real line, which does not seem to be possible. If $\kappa>4$ in
Theorem \ref{lim-1p*}, we can prove that if $\sigma\ge
1$ (resp.\ $\sigma\le -1$), then $L(\infty)$ is bounded from left
(resp.\ right) and unbounded from right (resp.\ left), and
$\lin{L(\infty)}\cap\R_\pi=\emptyset$.

\subsection{Three or four force points} First, we consider a strip Loewner
process with three force points. Let $\kappa>0$ and
$\rho_++\rho_-+\rho=\kappa-6$. Suppose $\beta(t)$, $0\le t<T$, is a
strip SLE$(\kappa;\rho_+,\rho_-,\rho)$ trace started from
$(0;+\infty,-\infty,p)$ for some $p\in\R_\pi$. Then $T=\infty$. Let
$\bar p=\Ree p$. Then the driving function $\xi(t)$, $0\le
t<\infty$, is the solution to the SDE:
\begin{equation}\left\{\begin{array}{lll} d\xi(t) & = &
\sqrt\kappa d B(t)+\frac{\rho_--\rho_+}2\, dt-\frac{\rho}2\tanh_2({\bar p(t)-\xi(t)})dt;\\
d\bar p(t) & = & \tanh_2({\bar
p(t)-\xi(t)})dt.\end{array}\right.\label{3pt}\end{equation}
Here $\bar p(t)\in\R$ and $\bar p(t)+\pi i=\psi(t,p)$ for any $t\ge
0$, where $\psi(t,\cdot)$, $0\le t<\infty$, are strip Loewner maps
driven by $\xi$. Let $X(t)=\bar p(t)-\xi(t)$. Then $X(t)$ satisfies
the SDE: \BGE dX(t)=-\sqrt\kappa dB(t)-\frac{\rho_--\rho_+}2\,
dt+(1+\frac{\rho}2)\tanh_2({X(t)})dt.\label{X}\EDE
Suppose $f$ is a real valued function on $\R$, and for any $x\in\R$,
$$f'(x)=\exp(x/2)^{\frac{4}\kappa\cdot\frac{\rho_--\rho_+}2}\cosh_2(x)^{-\frac 4\kappa(1+\frac{\rho}2)}.$$
>From Ito's formula, $f(X(t))$ is a local martingale.

Let $I=f(\R)$. Recall that $\rho=\kappa-6-\rho_+-\rho_-$. If
$\rho_+\ge \kappa/2-2$ and $\rho_-\ge \kappa/2-2$, then $I=\R$, so
a.s.\ $\limsup X(t)=+\infty$ and $\liminf X(t)=-\infty$. If
$\rho_+<\kappa/2-2$ and $\rho_-\ge \kappa/2-2$, then $I=(a,\infty)$
for some $a\in\R$, so a.s.\ $\lim X(t)=-\infty$. If
$\rho_+\ge\kappa/2-2$ and $\rho_-<\kappa/2-2$, then $I=(-\infty,b)$
for some $b\in\R$, so a.s.\ $\lim X(t)=+\infty$. If $\rho_+<
\kappa/2-2$ and $\rho_-< \kappa/2-2$, then $I=(a,b)$ for some
$a,b\in\R$, so with some probability $P\in(0,1)$, $\lim
X(t)=-\infty$; and with probability $1-P$, $\lim X(t)=+\infty$.

Let $I_1=[\kappa/2-2,\infty)$, $I_2=(\kappa/2-4,\kappa/2-2)$, and
$I_3=(-\infty,\kappa/2-4]$. Let Case (jk) denote the case that
 $\rho_+\in I_j$ and $\rho_-\in I_k$. We
use $(p,+\infty)$ or $(-\infty,p)$ to denote the open subarc of
$\R_\pi$ between $p$ and $+\infty$ or between $p$ and $-\infty$,
respectively.

\begin{Theorem} Suppose $\kappa\in(0,4]$. In Case (11),
a.s.\ $\lim_{t\to\infty}\beta(t)=p$. In Case (12), a.s.\ $\lim_{t\to
\infty}\beta(t)\in (-\infty,p)$. In Case (21), a.s.\ $\lim_{t\to
\infty}\beta(t)\in(p,+\infty)$. In Case (13), a.s.\ $\lim_{t\to
\infty}\beta(t)=-\infty$. In Case (31), a.s.\ $\lim_{t\to
\infty}\beta(t)=+\infty$. In Case (22), a.s.\  $\lim_{t\to
\infty}\beta(t)\in(-\infty,p)$ or $\in(p,+\infty)$. In Case (23),
a.s.\ $\lim_{t\to \infty}\beta(t)=-\infty$ or $\in(p,+\infty)$. In
Case (32), a.s.\ $\lim_{t\to \infty}\beta(t)\in(-\infty,p)$ or
$=+\infty$. In Case (33), a.s.\ $\lim_{t\to \infty}\beta(t)=-\infty$
or $=+\infty$. And in each of the last four cases, both events
happen with some positive probability. \label{3pt**}
\end{Theorem}
{\bf Proof.} The result in Case (11) follows from Theorem
\ref{lim=infty} and Lemma \ref{coordinate2}. Now consider Case (12).
We have a.s.\ $\lim X(t)=+\infty$. Let $Y(t)=X(t)+\sqrt\kappa B(t)$.
>From (\ref{X}), a.s.\
$$Y'(t)=-\frac{\rho_--\rho_+}2+(1+\frac{\rho}2)\tanh_2({X(t)})
\to \frac{\rho+\rho_+-\rho_-}2+1=\frac\kappa2-2-\rho_-$$ as
$t\to \infty$. Thus a.s.\ \BGE \lim_{t\to\infty}
X(t)/t=\lim_{t\to\infty} Y(t)/t= \kappa/2-2-\rho_->0.\label{X/t}\EDE
>From (\ref{3pt}), we see that as $t\to\infty$, the SDE for $\xi(t)$
tends to $d\xi(t)=\sqrt\kappa dB(t)+\sigma dt$, where
$\sigma:=\frac{\rho_--\rho_+}2-\rho/2=\rho_--(\kappa/2-3)\in(-1,1)$.
>From Theorem \ref{lim-1p}, it is reasonable to guess that a.s.\
$\lim_{t\to\infty} \beta(t)\in\R_\pi$. This will be rigorously
proved below.

Let \BGE
a(t)=\frac{\rho/2}{\sqrt\kappa}\big(1-\tanh_2({X(t)})\big);\label{Dt}\EDE
\BGE M(t)=\exp\big(-\int_0^t a(s)dB(s)-\frac 12 \int_0^t
a(s)^2ds\big).\label{Mt}\EDE From (\ref{X/t}), a.s.\ $\int_0^\infty
a(t)^2dt<\infty$, so a.s.\ $\lim_{t\to\infty} M(t)\in (0,\infty)$.
>From Ito's formula, $M(t)$ is a positive local martingale, and
$dM(t)/M(t)=-a(t)dB(t)$. For $N\in\N$, let $T_N\in[0,\infty]$ be the
largest number
 such that $M(t)\in(1/N,N)$ for $0\le t<T_N$. Then
$T_N$ is a stopping time, $M(t\wedge T_N)$ is a bounded martingale,
and $\PP[\{T_N=\infty\}]\to 1$ as $N\to \infty$. Define $\QQ$ such
that $d\QQ=M(T_N)d\PP$, where $M(\infty):=\lim_{t\to\infty}M(t)$.
Then $\QQ$ is also a probability measure. For $t\ge 0$, let $\til
B(t)=B(t)+\int_0^t a(s)ds$. From (\ref{3pt}), we have
$$\xi(t)=\xi(0)+\sqrt\kappa\til
B(t)+\sigma t.$$ From Girsanov Theorem, $\til B(t)$, $0\le t<T_N$,
is a partial $\QQ$-Brownian motion. Since $\kappa\in(0,4]$ and
$|\sigma|<1$, so from Theorem \ref{lim-1p},
$\QQ$-a.s.\ $\lim_{t\to T_N}\beta(t)\in\R_\pi$ on $\{T_N=\infty\}$. Since $1/N\le
d\QQ/d\PP\le N$, so $\QQ$ is equivalent to $\PP$. Thus ($\PP$-)a.s.\
$\lim_{t\to T_N}\beta(t)\in\R_\pi$ on $\{T_N=\infty\}$. For any $\eps>0$, there is $N$
such that $\PP[\{T_N=\infty\}]>1-\eps$. Thus with probability
greater than $1-\eps$, $\lim_{t\to\infty}\beta(t)\in\R_\pi$. Since
$\eps>0$ is arbitrary, so a.s.\ $\lim_{t\to\infty}\beta(t)\in\R_\pi$.
Now for any $x\in\R$ and $x\ge \bar p$, $\psi(t,x+\pi i)\in\R_\pi$
and
 $\Ree\psi(t,x+\pi i)\ge \Ree\psi(t,\bar p+\pi i)$ for any $t\ge 0$. Thus
 $\Ree\psi(t,x+\pi i)-\xi(t)\to\infty$ as $t\to \infty$.
>From an argument in the proof of Lemma \ref{onept}, we have
$\dist(x+\pi i,\beta((0,\infty)))>0$. Thus
$\lim_{t\to\infty}\beta(t)\not\in [p,+\infty)$, so a.s.\
$\lim_{t\to\infty}\beta(t)\in (-\infty,p)$.

Now consider Case (13). The argument is similar to that in Case (12)
except that now $\sigma=\rho_--(\kappa/2-3)\le -1$, so from Theorem \ref{lim-1p*},
we have a.s.\
$\lim_{t\to \infty}\beta(t)=-\infty$.  Case (21) and (31) are
symmetric to the above two cases. In Case (22), a.s.\
$\lim_{t\to\infty}X(t)=+\infty$ or $=-\infty$. If
$\lim_{t\to\infty}X(t)=+\infty$, then as $t\to\infty$, the SDE for
$\xi(t)$ tends to $d\xi(t)=\sqrt\kappa dB(t)+\sigma dt$, where
$\sigma =\rho_--(\kappa/2-3)\in(-1,1)$. Using the argument in Case
(12), we get a.s.\ $\lim_{t\to\infty}\beta(t)\in(-\infty,p)$
whenever $\lim_{t\to\infty}X(t)=+\infty$. Similarly, a.s.\
$\lim_{t\to\infty}\beta(t)\in(p,+\infty)$ whenever
$\lim_{t\to\infty}X(t)=-\infty$. The arguments in the other three
cases are similar to that in Case (22). $\Box$

\vskip 4mm

Next, we consider a strip Loewner process with four force points.
Let $\kappa>0$ and $\rho_++\rho_-+\rho_1+\rho_2=\kappa-6$. Suppose
$\beta(t)$, $0\le t<T$, is a strip
SLE$(\kappa;\rho_+,\rho_-,\rho_1,\rho_2)$ trace started from
$(0;+\infty,-\infty,p_1,p_2)$ for some $p_1,p_2\in\R$ with
$p_1>0>p_2$. Then the driving function $\xi(t)$, $0\le t<T$, is the
maximal solution to the SDE:
\begin{equation}\left\{\begin{array}{lll} d\xi(t) & = &
\sqrt\kappa d B(t)+\frac{\rho_--\rho_+}2
dt-\sum_{j=1}^2\frac{\rho_j}2\coth_2({p_j(t)-\xi(t)})dt;\\
dp_j(t) & = & \coth_2({p_j(t)-\xi(t)})dt,\quad j=1,2.
\end{array}\right.\label{4pt*}\end{equation}
Here $p_j(t)=\psi(t,p_j)\in\R$, $0\le t<T$, $j=1,2$, where
$\psi(t,\cdot)$, $0\le t<T$, are strip Loewner maps driven by $\xi$.

\begin{Theorem} Suppose $\kappa\in(0,4]$,
$\rho_1\ge \kappa/2-2$, $\rho_2\ge \kappa/2-2$,
$|(\rho_1+\rho_+)-(\rho_2+\rho_-)|<2$, and $\min\{\rho_1,\rho_2\}\le 0$.
Then a.s\ $T=\infty$ and $\lim_{t\to
\infty}\beta(t)\in \R_\pi$. \label{lim-1p3}
\end{Theorem}
{\bf Proof.} We only consider the case that $\rho_2\le 0$. The
case $\rho_1\le 0$ is symmetric. Let $X_j(t)=p_j(t)-\xi(t)$, $j=1,2$. Then
$X_1(t)>0>X_2(t)$, $0\le t<T$. And we have
$$ dX_1(t)=-\sqrt\kappa d
B(t)-\frac{\rho_--\rho_+}2dt+\big(1+\frac{\rho_1}2\big)\coth_2
({X_1(t)})dt+\frac{\rho_2}2\coth_2({X_2(t)})dt.$$
Define $f$ on $(0,\infty)$ such that for any $x>0$,
$$f'(x)=\exp(x/2)^{\frac
4\kappa\cdot\frac{\rho_--\rho_++\rho_2}2}\sinh_2(x)^{-\frac
4\kappa(1+\frac{\rho_1}2)}.$$ Then for any $x>0$,
$$\frac\kappa 2f''(x)=f'(x)\big(\frac{\rho_--\rho_++\rho_2}2-\big(1+\frac{\rho_1}2\big)\coth_2(
x)\big).$$ Let $Y(t)=f(X_1(t))$ for any
$t\in[0,T)$. From Ito's formula, we have
$$dY(t)=-\sqrt\kappa
f'(X_1(t))dB(t)+\frac{\rho_2}2f'(X_1(t))(1+\coth_2({X_2(t)}))dt.$$
>From the conditions of $\rho_j$'s, $f$ maps $(0,\infty)$ onto
$(-\infty,b)$ for some $b\in\R$. Since $\rho_2\le 0$ and $X_2(t)<0$,
so the drift is non-negative. Thus
a.s.\ $\lim_{t\to T}Y(t)=b$, which implies that $\lim_{t\to
T}X_1(t)=+\infty$. Let $Z(t)=X_1(t)+\sqrt\kappa
B(t)$. Since $\coth_2(X_2(t))<-1$ and $\rho_2\le 0$, so if $T=\infty$, then as
$t\to\infty$,
$$Z'(t)\ge \frac{\rho_+-\rho_--\rho_2}2+\big(1+\frac{\rho_1}2\big)\coth_2
({X_1(t)})\to 1+\frac{\rho_++\rho_1-\rho_--\rho_2}2.$$
Then $\liminf_{t\to \infty} X_1(t)/t=\liminf_{t\to \infty}
Z(t)/t\ge\sigma:=1+({\rho_++\rho_1-\rho_--\rho_2})/2>0$.

Let $a(t)$ and $M(t)$ be defined by (\ref{Dt}) and (\ref{Mt}) except
that $\rho$ and $\tanh_2(X(t))$ in (\ref{Dt}) are replaced by
$\rho_1$ and $\coth_2(X_1(t))$, respectively. If $T=\infty$, since $\liminf_{t\to
\infty} X_1(t)/t\ge\sigma>0$, so a.s.\ $\lim_{t\to\infty} M(t)\in
(0,\infty)$. This is clearly true if $T<\infty$ because $a(s)$ is bounded.
 Let $\til B(t)=B(t)+\int_0^t a(s)ds$, $0\le t<T$. From (\ref{4pt*})
we have
$$d\xi(t)=\sqrt\kappa d\til B(t)+\frac{\rho_--\rho_+-\rho_1}2
dt-\frac{\rho_2}2\coth_2({X_2(t)})dt.$$ If under some probability measure
$\QQ$, $(\til B(t))$ is a partial
Brownian motion, then $\beta(t)$, $0\le t<T$, is a partial strip
SLE$(\kappa;\rho'_+,\rho_-,\rho_2)$ process started from
$(0;+\infty,-\infty,p_2)$, where $\rho'_+=\rho_++\rho_1$. Since
$\rho_+'+\rho_-+\rho_2=\kappa-6$,
$\rho_+'\in(\kappa/2-4,\kappa/2-2)$ and $\rho_2\ge \kappa/2-2$, so
from Lemma \ref{coordinate*} and Theorem \ref{lim-1p3}, we have $\QQ$-a.s.\
$\lim_{t\to T}\beta(t)\in\R_\pi\cup \St_\pi$. From the proof
in Case (12) of Theorem \ref{lim-1p3}, we have a.s.\
$\lim_{t\to T}\beta(t)\in\R_\pi\cup \St_\pi$.
Since $\beta$ is a full  trace, it separates either $p_1$ or $p_2$
from $\R_\pi$ in $\St_\pi$, so $\lim_{t\to T}\beta(t)\in\St_\pi$ is
not possible. Thus $\lim_{t\to T}\beta(t)\in\R_\pi$ a.s.. This implies
that $T=\lim_{t\to T}\scap(\beta((0,t])=\infty$.
$\Box$

\section{Coupling of Two SLE Processes} \label{coupl}
Let $\kappa_1,\kappa_2>0$; $\kappa_1\kappa_2=16$; $\rho_{j,m}\in\R$,
$1\le m\le N$, $j=1,2$, $N\in\N$;
$\rho_{2,m}=-\kappa_2\rho_{1,m}/4$, $1\le m\le N$;
$x_1,x_2,p_1,\dots,p_N\in\R$ are distinct points. Let
$\vec{\rho}_j=(\rho_{j,1},\dots,\rho_{j,N})$, $j=1,2$, and
$\vec{p}=(p_1,\dots,p_N)$. Note that if $\kappa_1=\kappa_2=4$, then
$\vec{\rho}_1+\vec{\rho}_2=\vec{0}$; if $\kappa_1,\kappa_2\ne 4$,
then $\vec{\rho}_1/(\kappa_1-4)=\vec{\rho}_2/(\kappa_2-4)$. The goal
of this section is to prove the following theorem.

\begin{Theorem} There is a coupling of $K_1(t)$,
$0\le t<T_1$, and $K_2(t)$, $0\le t<T_2$, such that (i) for $j=1,2$,
$K_j(t)$, $0\le t<T_j$, is a chordal
SLE$(\kappa_j;-\frac{\kappa_j}2,\vec{\rho}_j)$ process started from
$(x_j;x_{3-j},\vec{p})$; and (ii) for $j\ne k\in\{1,2\}$, if $\bar
t_k$ is an $(\F^k_t)$-stopping time with $\bar t_k<T_k$, then
conditioned on $\F^k_{\bar t_k}$, $\vphi_k(\bar t_k,K_j(t))$, $0\le
t\le T_j(\bar t_k)$, has the same distribution as a time-change of a
partial chordal SLE$(\kappa_j;-\frac{\kappa_j}2,\vec{\rho}_j)$
process started from $(\vphi_k(\bar t_k,x_j);\xi_k(\bar
t_k),\vphi_k(\bar t_k,\vec{p}))$, where
$\vphi_k(t,\vec{p})=(\vphi_k(t,p_1),\dots,\vphi_k(t,p_N))$,
$\vphi_k(t,\cdot)=\vphi_{K_k(t)}$, $T_j(\bar t_k)\in(0,T_j]$ is the
largest number such that $\lin{K_j(t)}\cap\lin{K_k(\bar t_k)}=\emptyset$
for $0\le t<T_j(\bar t_k)$, and $(\F^j_t)$ is the filtration
generated by $(K_j(t))$, $j=1,2$.
 \label{bartk}
\end{Theorem}

In many cases we can prove that $\vphi_k(\bar t_k,K_j(t))$, $0\le
t\le T_j(\bar t_k)$, has the same distribution as a time-change of a
{\it full} chordal SLE$(\kappa_j;-\frac{\kappa_j}2,\vec{\rho}_j,)$
process started from $(\vphi_k(\bar t_k,x_j);\xi_k(\bar
t_k),\vphi_k(\bar t_k,\vec{p}))$. From the property of $T_j(\bar
t_k)$, ${\cup_{0\le t<T_j(\bar t_k)} K_j(t)}$ touches ${K_k(\bar
t_k)}$, so $\cup_{0\le t<T_j(\bar t_k)} \vphi_k(\bar t_k,K_j(t))$
touches $\R$. So the chain can not be further extended while staying
bounded away from the boundary. Thus if $\kappa_j\le 4$,  it is a
full process. Another case is when there is some force point $p_k$
that lies between $x_1$ and $x_2$. Then ${\cup_{0\le t<T_j(\bar
t_k)} K_j(t)}$ separates $\vphi_k(\bar t_k,p_k)$ from $\infty$. So
again we get a full process.

\subsection{Ensembles} Let's review the results in Section 3 of
\cite{reversibility}. For $j=1,2$, let $K_j(t)$ and
$\vphi_j(t,\cdot)$, $0\le t<S_j$, be chordal Loewner hulls and maps
driven by $\xi_j\in C([0,S_j))$. Suppose $K_1(t_1)\cap
K_2(t_2)=\emptyset$ for any $t_1\in[0,S_1)$ and $t_2\in[0,S_2)$. For
$j\ne k\in\{1,2\}$, $t_0\in[0,S_k)$ and $t\in[0,S_j)$, let
\BGE K_{j,t_0}(t)=(K_j(t)\cup K_k(t_0))/K_k(t_0),\quad
\vphi_{j,t_0}(t,\cdot)=\vphi_{K_{j,t_0}(t)}.\label{Ensemble}\EDE
Then for any $t_1\in[0,S_1)$ and $t_2\in[0,S_2)$,
\BGE\vphi_{K_1(t_1)\cup
K_2(t_2)}=\vphi_{1,t_2}(t_1,\cdot)\circ\vphi_2(t_2,\cdot)=
\vphi_{2,t_1}(t_2,\cdot)\circ\vphi_1(t_1,\cdot).\label{circ=1}\EDE
We use $\pa_1$ and $\pa_z$ to denote the partial derivatives of
$\vphi_j(\cdot,\cdot)$ and $\vphi_{j,t_0}(\cdot,\cdot)$ w.r.t.\ the
first (real) and second (complex) variables, respectively, inside
the bracket; and use $\pa_0$ to denote the partial derivative of
$\vphi_{j,t_0}(\cdot,\cdot)$ w.r.t.\ the subscript $t_0$. From (3.10$\sim$3.14) in
Section 3 of \cite{reversibility}, we have \BGE
\pa_0\vphi_{k,t}(s,\xi_j(t))=-3\pa_z^2\vphi_{k,t}(s,\xi_j(t));\label{-3}\EDE
\BGE
\frac{\pa_0\pa_z\vphi_{k,t}(s,\xi_j(t))}{\pa_z\vphi_{k,t}(s,\xi_j(t))}
=\frac
12\cdot\left(\frac{\pa_z^2\vphi_{k,t}(s,\xi_j(t))}{\pa_z\vphi_{k,t}(s,\xi_j(t))}
\right)^2-\frac
43\cdot\frac{\pa_z^3\vphi_{k,t}(s,\xi_j(t))}{\pa_z\vphi_{k,t}(s,\xi_j(t))};
\label{1/2-4/3}\EDE \label{-3-1/2-4/3}
\BGE\pa_1\vphi_{j,t_0}(t,z)=\frac{2\pa_z\vphi_{k,t}(t_0,\xi_j(t))^2}
{\vphi_{j,t_0}(t,z)-\vphi_{k,t}(t_0,\xi_j(t))};\label{chordal*}\EDE
\BGE \frac{\pa_1\pa_z\vphi_{j,s}(t,z)}{\pa_z\vphi_{j,s}(t,z)}
=\frac{-2\pa_z\vphi_{k,t}(s,\xi_j(t))^2}
{(\vphi_{j,s}(t,z)-\vphi_{k,t}(s,\xi_j(t)))^2};\label{chordal**}\EDE
\BGE
\pa_1\Big(\frac{\pa_z^2\vphi_{j,s}(t,z)}{\pa_z\vphi_{j,s}(t,z)}\Big)
=\frac{4\pa_z\vphi_{k,t}(s,\xi_j(t))^2\pa_z\vphi_{j,s}(t,z)}
{(\vphi_{j,s}(t,z)-\vphi_{k,t}(s,\xi_j(t)))^3};\label{chordal*2}\EDE
\BGE
\pa_1\pa_z\Big(\frac{\pa_z^2\vphi_{j,s}(t,z)}{\pa_z\vphi_{j,s}(t,z)}\Big)
=\frac{4\pa_z\vphi_{k,t}(s,\xi_j(t))^2\pa_z^2\vphi_{j,s}(t,z)}
{(\vphi_{j,s}(t,z)-\vphi_{k,t}(s,\xi_j(t)))^3}-\frac{12\pa_z\vphi_{k,t}(s,\xi_j(t))^2\pa_z\vphi_{j,s}(t,z)^2}
{(\vphi_{j,s}(t,z)-\vphi_{k,t}(s,\xi_j(t)))^4}.\label{chordal*3}\EDE

\subsection{Martingales} \label{subsec-mart}
 Suppose $x_1,x_2,p_1,\dots,p_N$ are distinct
points on $\R$. Let $\xi_j\in C([0,T_j))$, $j=1,2$, be two independent
semi-martingales that satisfy $d\langle \xi_j(t)\rangle
=\kappa_jdt$, where $\kappa_1,\kappa_2>0$. Let $\vphi(t,\cdot)$ and
$K_j(t)$, $0\le t<\infty$, be chordal Loewner maps and hulls driven
by $\xi_j$, $j=1,2$. Let \BGE {\cal D}:=\{(t_1,t_2): \lin{K_1(t_1)}\cap
\lin{K_2(t_2)}=\emptyset, \vphi(t_j,p_m)\mbox{ does not blow up,
}1\le m\le N,j=1,2\}.\label{calD}\EDE  For $(t_1,t_2)\in {\cal D}$,
$j\ne k\in\{1,2\}$, and $h\in\Z_{\ge 0}$, let
$A_{j,h}(t_1,t_2)=\pa_z^h\vphi_{k,t_j}(t_k,\xi_j(t_j))$.
For $(t_1,t_2)\in {\cal D}$, $1\le m\le N$,
and $h\in\Z_{\ge 0}$, let
$B_{m,h}(t_1,t_2)=\pa_z^h\vphi_{K_1(t_1)\cup K_2(t_2)}(p_m)$. For $j=1$, $k=2$, and $1\le m\le N$,
we have the following SDEs:
\BGE \pa_j
A_{j,0}= A_{j,1} \pa \xi_j(t_j)+(\frac{\kappa_j} 2-3)A_{j,2}\pa
t_j;\label{pa-A_1-0}\EDE \BGE \frac{\pa_j A_{j,1}}{A_{j,1}}=
\frac{A_{j,2}} {A_{j,1}}\,\pa \xi_j(t_j)+\Big(\frac
12\cdot\frac{A_{j,2}^2} {A_{j,1}^2}+\Big(\frac{\kappa_j} 2-\frac
43\Big)\cdot\frac{A_{j,3}} {A_{j,1}}\Big)\,\pa
t_j;\label{pa-A_1-1}\EDE \BGE \pa_j A_{k,0}=\frac{2
A_{j,1}^2}{A_{k,0}-A_{j,0}}\,\pa t_j,\quad \frac{\pa_j
A_{k,1}}{A_{k,1}}=\frac{-2 A_{j,1}^2}{(A_{k,0}-A_{j,0})^2}\,\pa
t_j;\label{pa-A_2}\EDE \BGE \pa_j B_{m,0}=\frac{2
A_{j,1}^2}{B_{m,0}-A_{j,0}}\,\pa t_j,\quad \frac{\pa_j B_{m,1}}{
B_{m,1}}=\frac{-2 A_{j,1}^2}{(B_{m,0}-A_{j,0})^2}\,\pa
t_j.\label{Bm}\EDE
Here $\pa_j$ means the partial derivative w.r.t.\ $t_j$.
Note that (\ref{pa-A_1-0}) and (\ref{pa-A_1-1}) are
(4.10) and (4.11) in \cite{reversibility}; (\ref{pa-A_2}) follows
from (\ref{chordal*}) and (\ref{chordal**}) here; and (\ref{Bm}) follows from
(\ref{chordal*}), (\ref{chordal**}), and (\ref{circ=1}).
By symmetry, (\ref{pa-A_1-0}$\sim$
\ref{Bm}) also hold for $j=2$ and $k=1$.

For $j\ne k\in\{1,2\}$, let $E_{j,0}=A_{j,0}-A_{k,0}=-E_{k,0}$;
$E_{j,m}=A_{j,0}-B_{m,0}$, $1\le m\le N$; and
$C_{m_1,m_2}=B_{m_1,0}-B_{m_2,0}$, $1\le m_1<m_2\le N$. From
(\ref{pa-A_1-0}), (\ref{pa-A_2}), and (\ref{Bm}), for
$0\le m \le N$, \BGE \frac{\pa_j
E_{j,m}}{E_{j,m}}=\frac{A_{j,1}}{E_{j,m}}\,\pa\xi_j(t_j)+\Big(
(\frac{\kappa_j}2-3)\cdot\frac{A_{j,2}}{E_{j,m}}+2\,\frac{A_{j,1}^2}{E_{j,m}^2}\Big)\,\pa
t_j.\label{Djm}\EDE From (\ref{pa-A_2}) and (\ref{Bm}), for
$1\le m\le N$ and $1\le m_1< m_2\le N$ \BGE \frac{\pa_j
E_{k,m}}{E_{k,m}}=\frac{-2A_{j,1}^2}{E_{j,0}E_{j,m}}\,\pa t_j,\qquad
\frac{\pa_j
C_{m_1,m_2}}{C_{m_1,m_2}}=\frac{-2A_{j,1}^2}{E_{j,m_1}E_{j,m_2}}\,\pa
t_j.\label{Dkm}\EDE

Now suppose $\kappa_1\kappa_2=16$. For $j=1,2$, let \BGE \alpha_j
=\frac{6-\kappa_j}{2\kappa_j},\qquad
\lambda_j=\frac{(8-3\kappa_j)(6-\kappa_j)}{2\kappa_j}.\label{alpha-lambda}\EDE
Then $\lambda_1=\lambda_2$. Let it be denoted by $\lambda$. From (\ref{pa-A_1-1}) and
(\ref{pa-A_2}), we have \BGE \frac{\pa_j
A_{j,1}^{\alpha_j}}{A_{j,1}^{\alpha_j}}=\frac{6-\kappa_j}{2\kappa_j}\cdot
\frac{A_{j,2}} {A_{j,1}}\,\pa \xi_j(t_j)+\lambda\,\Big(\frac
14\cdot\frac{A_{j,2}^2} {A_{j,1}^2}-\frac 16\cdot\frac{A_{j,3}}
{A_{j,1}}\Big)\,\pa t_j;\label{pa-A-1-alpha}\EDE \BGE \frac{\pa_j
A_{k,1}^{\alpha_k}}{A_{k,1}^{\alpha_k}}=-2\alpha_k\,
\frac{A_{j,1}^2}{E_{j,0}^2}\,\pa t_j=-\frac{3\kappa_j-8}{8}\,
\frac{A_{j,1}^2}{E_{j,0}^2} \,\pa t_j.\label{pa-A-2-alpha}\EDE

Suppose
$\vec{\rho}_j=(\rho_{j,1},\dots,\rho_{j,N})\in\R^N$, $j=1,2$, and
$\vec\rho_{2}=-\frac{\kappa_2}4\vec{\rho}_{1}$. Let
$\rho^*_{j,m}=\rho_{j,m}/\kappa_j$, $1\le m\le N$, $j=1,2$. Then
$\rho^*_{2,m}=-\kappa_1\rho^*_{1,m}/4$ and
$\rho^*_{1,m}=-\kappa_2\rho^*_{2,m}/4$ for $1\le m\le N$. From
(\ref{Djm}) and (\ref{Dkm}), for $j\ne k\in\{1,2\}$ and $1\le m\le N$,
we have \BGEN \frac{\pa_j
|E_{j,m}|^{\rho^*_{j,m}}}{|E_{j,m}|^{\rho^*_{j,m}}}=\rho^*_{j,m}\,
\frac{A_{j,1}}{E_{j,m}}\,\pa\xi_j(t_j)+\rho^*_{j,m}\cdot\frac{\kappa_j-6}2
\cdot\frac{A_{j,2}}{E_{j,m}}\,\pa t_j+\EDEN \BGE
+\Big(\frac{\kappa_j}2\,\rho^*_{j,m}(\rho^*_{j,m}-1)
+2\rho^*_{j,m}\Big)\, \frac{A_{j,1}^2}{E_{j,m}^2}\,\pa
t_j;\label{Djm-beta}\EDE \BGE \frac{\pa_j
|E_{k,m}|^{\rho^*_{k,m}}}{|E_{k,m}|^{\rho^*_{k,m}}}=
-2\rho^*_{k,m}\,\frac{A_{j,1}^2}{E_{j,0}E_{j,m}}\,\pa
t_j=\frac{\kappa_j\rho^*_{j,m}}{2}\cdot\frac{A_{j,1}^2}{E_{j,0}E_{j,m}}\,\pa
t_j.\label{Dkm-beta}\EDE Let $E=|E_{1,0}|=|E_{2,0}|$. From
(\ref{Djm}), for $j=1,2$, \BGE \frac{\pa_j
E^{-1/2}}{E^{-1/2}}=-\frac12\cdot
\frac{A_{j,1}}{E_{j,0}}\,\pa\xi_j(t_j)+\Big(\frac{6-\kappa_j}4
\cdot\frac{A_{j,2}}{E_{j,0}} +\frac{3\kappa_j -8}{8}\,
\frac{A_{j,1}^2}{E_{j,0}^2} \Big)\,\pa t_j.\label{Djm-beta-0}\EDE
For $1\le m\le N$, let
$$\gamma_m=\frac{\kappa_1}4\,\rho^*_{1,m}(\rho^*_{1,m}-1)
+\rho^*_{1,m}=\frac{\kappa_2}4\,\rho^*_{2,m}(\rho^*_{2,m}-1)
+\rho^*_{2,m}.$$ For $j=1,2$, from (\ref{Bm}) we have \BGE
\frac{\pa_j
B_{m,1}^{\gamma_m}}{B_{m,1}^{\gamma_m}}=-(\frac{\kappa_j}2\,
\rho^*_{j,m}(\rho^*_{j,m}-1)
+2\rho^*_{j,m})\,\frac{A_{j,1}^2}{E_{j,m}^2} \,\pa
t_j.\label{B-gamma}\EDE For $1\le m_1< m_2\le N$, let
$$\delta_{m_1,m_2}=\frac{\kappa_1}2\rho^*_{1,m_1}\rho^*_{1,m_2}
=\frac{\kappa_2}2\rho^*_{2,m_1}\rho^*_{2,m_2}.$$ From (\ref{Dkm}),
for $j=1,2$, \BGE \frac{\pa_j
|C_{m_1,m_2}|^{\delta_{m_1,m_2}}}{|C_{m_1,m_2}|^{\delta_{m_1,m_2}}}=
-{\kappa_j}\rho^*_{j,m_1}\rho^*_{j,m_2}\,\frac{A_{j,1}^2}
{E_{j,m_1}E_{j,m_2}}\,\pa t_j.\label{m1m2}\EDE For $(t_1,t_2)\in\cal
D$, let \BGE F(t_1,t_2)=\exp\Big(\int_0^{t_2}\!\int_0^{t_1}
\frac{2A_{1,1}(s_1,s_2)A_{2,1}(s_1,s_2)}{E(s_1,s_2)^2}
\,ds_1ds_2\Big).\label{def-E}\EDE From (4.15) in
\cite{reversibility}, for $j=1,2$, \BGE \frac{\pa_j F^{-\lambda}}{F^{-\lambda}}=-\lambda\Big(\frac
14\cdot \frac{A_{j,2}^2}{A_{j,1}^2}-\frac 16\cdot\frac{A_{j,3}}
{A_{j,1}}\Big)\,\pa t_j.\label{E}\EDE Let \BGE \til
M=\frac{A_{1,1}^{\alpha_1}A_{2,1}^{\alpha_2}}{E^{1/2}}\,F^{-\lambda}\,\prod_{m=1}^N
\Big(B_{m,1}^{\gamma_m}\prod_{j=1}^2|E_{j,m}|^{\rho^*_{j,m}}\Big)
\prod_{1\le m_1<m_2\le N}
|C_{m_1,m_2}|^{\delta_{m_1,m_2}}.\label{def-tilM}\EDE

Now we compute the SDE for $\pa_j\til M/\til M$ in terms of $\pa\xi_j(t_j)$
and $\pa t_j$. The coefficient of the $\pa\xi_j(t_j)$
term should be the sum of the coefficients of the $\pa\xi_j(t_j)$
terms in (\ref{pa-A-1-alpha}$\sim$\ref{E}).
The SDEs in (\ref{pa-A-1-alpha}$\sim$\ref{E}) that contain stochastic terms are (\ref{pa-A-1-alpha}),
(\ref{Djm-beta}), (\ref{Djm-beta-0}). So the sum is equal to
\BGE \frac{6-\kappa_j}{2\kappa_j}\cdot \frac{A_{j,2}}
{A_{j,1}}-\frac12\cdot
\frac{A_{j,1}}{E_{j,0}}+\sum_{m=1}^N\rho^*_{j,m}\,
\frac{A_{j,1}}{E_{j,m}}.\label{MS}\EDE
The coefficient of the $\pa t_j$ term equals to the sum
of the coefficients of the $\pa t_j$
terms in (\ref{pa-A-1-alpha}$\sim$\ref{E})  plus the sum of the coefficients of the drift terms
coming out of products.
The drift term in the SDE for $\pa_j\til M/\til M$ contributed
by the products of (\ref{pa-A-1-alpha}) and SDEs in (\ref{Djm-beta}) is
\BGE \kappa_j\cdot\sum_{m=1}^N \frac{6-\kappa_j}{2\kappa_j}\cdot
\frac{A_{j,2}} {A_{j,1}}\cdot\rho^*_{j,m}\cdot
\frac{A_{j,1}}{E_{j,m}}=-\sum_{m=1}^N \rho^*_{j,m}\cdot \frac{\kappa_j-6}{2}\cdot
\frac{A_{j,2}}{E_{j,m}}\label{*1}\EDE
The drift term contributed
by the products of (\ref{Djm-beta-0}) and SDEs in (\ref{Djm-beta}) is
\BGE \kappa_j\cdot \sum_{m=1}^N \rho^*_{j,m}\cdot
\frac{A_{j,1}}{E_{j,m}}\cdot\Big(- \frac12\cdot
\frac{A_{j,1}}{E_{j,0}}\Big)=-\sum_{m=1}^N \frac{\kappa_j\rho^*_{j,m}}{2}
\cdot\frac{A_{j,1}^2}{E_{j,0}E_{j,m}}.\label{*2}\EDE
The drift term contributed
by the product of (\ref{pa-A-1-alpha}) and (\ref{Djm-beta-0}) is
\BGE \kappa_j\cdot \frac{6-\kappa_j}{2\kappa_j}\cdot\frac{A_{j,2}} {A_{j,1}}\cdot\Big(- \frac12\cdot
\frac{A_{j,1}}{E_{j,0}}\Big)=-\frac{6-\kappa_j}4
\cdot\frac{A_{j,2}}{E_{j,0}}.\label{*3}\EDE
The drift term contributed
by the products of pairs of SDEs in (\ref{Djm-beta}) is
\BGE \kappa_j\cdot \sum_{1\le m_1< m_2\le N} \rho^*_{j,m_1}\cdot
\frac{A_{j,1}}{E_{j,m_1}} \cdot \rho^*_{j,m_2}\cdot
\frac{A_{j,1}}{E_{j,m_2}}=\sum_{1\le m_1< m_2\le N} {\kappa_j}\rho^*_{j,m_1}\rho^*_{j,m_2}\,\frac{A_{j,1}^2}
{E_{j,m_1}E_{j,m_2}}.\label{*4}\EDE
The sum of the coefficients of the $\pa t_j$
terms in (\ref{pa-A-1-alpha}$\sim$\ref{E}) is equal to
$$\lambda\,\Big(\frac
14\cdot\frac{A_{j,2}^2} {A_{j,1}^2}-\frac 16\cdot\frac{A_{j,3}}
{A_{j,1}}\Big)-\frac{3\kappa_j-8}{8}\,
\frac{A_{j,1}^2}{E_{j,0}^2}+\sum_{m=1}^N\rho^*_{j,m}\cdot\frac{\kappa_j-6}2
\cdot\frac{A_{j,2}}{E_{j,m}}$$$$+\sum_{m=1}^N\Big(\frac{\kappa_j}2\,\rho^*_{j,m}(\rho^*_{j,m}-1)
+2\rho^*_{j,m}\Big)\, \frac{A_{j,1}^2}{E_{j,m}^2}
+\sum_{m=1}^N \frac{\kappa_j\rho^*_{j,m}}{2}\cdot\frac{A_{j,1}^2}{E_{j,0}E_{j,m}}$$
$$+\Big(\frac{6-\kappa_j}4
\cdot\frac{A_{j,2}}{E_{j,0}} +\frac{3\kappa_j -8}{8}\,
\frac{A_{j,1}^2}{E_{j,0}^2} \Big)-\sum_{m=1}^N(\frac{\kappa_j}2\,
\rho^*_{j,m}(\rho^*_{j,m}-1)
+2\rho^*_{j,m})\,\frac{A_{j,1}^2}{E_{j,m}^2} $$
$$-\sum_{1\le m_1< m_2\le N}  {\kappa_j}\rho^*_{j,m_1}\rho^*_{j,m_2}\,\frac{A_{j,1}^2}
{E_{j,m_1}E_{j,m_2}}-\lambda\Big(\frac
14\cdot \frac{A_{j,2}^2}{A_{j,1}^2}-\frac 16\cdot\frac{A_{j,3}}
{A_{j,1}}\Big)$$
$$=\sum_{m=1}^N\rho^*_{j,m}\cdot\frac{\kappa_j-6}2
\cdot\frac{A_{j,2}}{E_{j,m}}
+\sum_{m=1}^N \frac{\kappa_j\rho^*_{j,m}}{2}\cdot\frac{A_{j,1}^2}{E_{j,0}E_{j,m}}$$
\BGE +\frac{6-\kappa_j}4
\cdot\frac{A_{j,2}}{E_{j,0}}
-\sum_{1\le m_1< m_2\le N} {\kappa_j}\rho^*_{j,m_1}\rho^*_{j,m_2}\,\frac{A_{j,1}^2}
{E_{j,m_1}E_{j,m_2}}.\label{*0}\EDE
>From (\ref{*1}$\sim$\ref{*0}), the SDE for $\pa_j\til M/\til M$ has no $\pa t_j$ terms. Thus
from (\ref{MS}), for $j=1,2$, we have
\BGE\frac{\pa_j\til M}{\til
M}=\Big(\frac{6-\kappa_j}{2\kappa_j}\cdot \frac{A_{j,2}}
{A_{j,1}}-\frac12\cdot
\frac{A_{j,1}}{E_{j,0}}+\sum_{m=1}^N\rho^*_{j,m}\,
\frac{A_{j,1}}{E_{j,m}}\Big)\,\pa\xi_j(t_j).\label{tilM}\EDE For
$(t_1,t_2)\in\cal D$, let \BGE M(t_1,t_2)=\frac{\til M(t_1,t_2)\til
M(0,0)}{\til M(t_1,0)\til M(0,t_2)}.\label{def-M}\EDE Then
$M(t_1,0)=M(0,t_2)=1$ for $t_j\in[0,T_j)$, $j=1,2$.

The process $\til M$ turns out to be the local Radon-Nikodym derivative
of the coupling measure in Theorem \ref{bartk} w.r.t.\ the product measure of
two standard chordal SLE$(\kappa)$ processes. Fix $j\ne k\in\{1,2\}$. Such $\til M$
must satisfy SDE (\ref{tilM}). So there are factors $A_{j,1}^{\alpha_j}$,
$\prod_m |E_{j,m}|^{\rho^*_{j,m}}$, and $E^{-1/2}$ in (\ref{def-tilM}).
Other factors in (\ref{def-tilM}) make $\til M$ a local martingale in $t_j$, for any fixed $t_k$.
Moreover, if $t_j$ is fixed, $\til M$ should also be a local martingale in $t_k$.
And we expect some symmetry between $j$ and $k$ in the definition of $\til M$.
This gives
restrictions on the values of $\kappa_j$  and $\rho_{j,m}$, $j=1,2$, $1\le m\le N$. The process
$M$ then becomes the local Radon-Nikodym derivative
of the coupling measure in Theorem \ref{bartk} w.r.t.\ the product of
its marginal measures. The property of $M$ will be checked later.

\vskip 3mm

 Let $B_1(t)$ and $B_2(t)$ be independent Brownian
motions. Let $(\F^j_t)$ be the filtration generated by $B_j(t)$,
$j=1,2$. Fix $j\ne k\in \{1,2\}$. Suppose $\xi_j(t)$, $0\le t<T_j$,
is the maximal solution to the SDE: \BGE d\xi_j(t)=\sqrt{\kappa_j}d
B_j(t)+\Big(
\frac{-\kappa_j/2}{\xi_j(t)-\vphi_j(t_j,x_k)}+\sum_{m=1}^N\,
\frac{\rho_{j,m}}{\xi_j(t)-\vphi_j(t_j,p_m)}\Big)\,dt,\label{SDE-xi}\EDE
with $\xi_j(0)=x_j$. Then $(K_j(t),0\le t<T_j)$ is an
SLE$(\kappa_j;-\frac{\kappa_j}2,\vec{\rho}_{j})$ process started
from $(x_1;x_2,\vec{p})$. Since $\vphi_{k,t}(t,\cdot)=\id$, so at
$t_j=t$ and $t_k=0$, $A_{j,1}=1$, $A_{j,2}=0$,
$E_{j,0}=\xi(j)-\vphi_j(t,\xi_k)$, and
$E_{j,m}=\xi(j)-\vphi_j(t,p_m)$, $1\le m\le N$. Thus
$$\Big(\frac{6-\kappa_j}{2\kappa_j}\cdot \frac{A_{j,2}}
{A_{j,1}}-\frac12\cdot
\frac{A_{j,1}}{E_{j,0}}+\sum_{m=1}^N\rho^*_{j,m}\,
\frac{A_{j,1}}{E_{j,m}}\Big)\Big|_{t_j=t,t_k=0}$$ \BGE
=\frac{-1/2}{\xi_j(t)-\vphi_j(t_j,x_k)}+\sum_{m=1}^N\,
\frac{\rho^*_{j,m}}{\xi_j(t)-\vphi_j(t_j,p_m)}.\label{(t,0)}\EDE For
$j\ne k\in\{1,2\}$ and $t_k\in [0,T_k)$, let $T_j(t_k)\in(0,T_j]$ be
the largest number such that for $0\le t<T_j(t_k)$, $\lin{K_j(t)}\cap
\lin{K_k(t_k)}=\emptyset$.

\begin{Theorem} Fix $j\ne k\in\{1,2\}$. Let $\bar t_k$ be an $(\F^k_t)$-stopping
time. Then the process $t\mapsto M|_{t_j=t,t_k=\bar t_k}$, $0\le
t<T_j(\bar t_k)$, is an $(\F^j_t\times\F^k_{\bar t_k})_{t\ge
0}$-local martingale, and
 $$\frac{\pa_j  M}{
M}\Big|_{t_j=t,t_k=\bar
t_k}=\Big[\Big(\frac{6-\kappa_j}{2\kappa_j}\cdot \frac{A_{j,2}}
{A_{j,1}}-\frac12\cdot
\frac{A_{j,1}}{E_{j,0}}+\sum_{m=1}^N\rho^*_{j,m}\,
\frac{A_{j,1}}{E_{j,m}}\Big)\Big|_{t_j=t,t_k=\bar t_k}$$\BGE -\Big(
\frac{-1/2}{\xi_j(t)-\vphi_j(t,x_k)}+\sum_{m=1}^N\,
\frac{\rho^*_{j,m}}{\xi_j(t)-\vphi_j(t,p_m)}\Big)\Big]
\,\sqrt{\kappa_j}\pa B_j(t).\label{M-SDE}\EDE \label{M}
\end{Theorem}
{\bf Proof.} This follows from (\ref{tilM}$\sim$\ref{(t,0)}), where
all functions are valued at $t_j=t$ and $t_k=\bar t_k$, and all SDE
are $(\F^j_t\times\F^k_{\bar t_k})$-adapted. $\Box$

\vskip 3mm

Now we make some improvement over the above theorem. Let $\bar t_2$
be an $(\F^2_t)$-stopping time with $\bar t_2<T_2$. Suppose $R$ is
an $(\F^1_{t}\times\F^2_{\bar t_2})_{t\ge 0}$-stopping time with $R<
T_1(\bar t_2)$. Let $\F_{R,\bar t_2}$ denote the $\sigma$-algebra
obtained from the filtration $(\F^1_t\times\F^2_{\bar t_2})_{t\ge
0}$ and its stopping time $R$, i.e., ${\cal E}\in \F_{R,\bar t_2}$
iff for any $t\ge 0$, ${\cal E}\cap\{R\le t\}\in
\F^1_t\times\F^2_{\bar t_2}$. For every $t\ge 0$, $R+t$ is also an
$(\F^1_{t}\times\F^2_{\bar t_2})_{t\ge 0}$-stopping time. So we have
a filtration $(\F_{R+t,\bar t_2})_{t\ge 0}$.

\begin{Theorem} Let $\bar t_2$ and $R$ be as above. Let
$I\in[0,\bar t_2]$ be $\F_{R,\bar t_2}$-measurable. Then $(M(R+t,I),
0\le t<T_1(I)-R)$ is a continuous $(\F_{R+t,\bar t_2})_{t\ge
0}$-local martingale.
 \label{M-I}
\end{Theorem}
{\bf Proof.} (i)
Let $B_1^R(t)=B_1(R+t)-B_1(R)$, $0\le t<\infty$. Since $B_1(t)$ is
an $(\F^1_{t}\times\F^2_{\bar t_2})_{t\ge 0}$-Brownian motion, so
$B_1^R(t)$ is an $(\F_{R+t,\bar t_2})_{t\ge 0}$-Brownian motion. Since $\vphi_1(R+t,\cdot)$
is $(\F^1_{t}\times\F^2_{\bar t_2})_{t\ge 0}$-adapted, so $\xi_1(R+t)$ satisfies
the $(\F^1_{t}\times\F^2_{\bar t_2})_{t\ge 0}$-adapted SDE
$$d\xi_1(R+t)= \sqrt\kappa
dB^R_1(t)d\xi_j(t)+
\frac{-\kappa_1/2}{\xi_1(R+t)-\vphi_1(R+t,x_2)}\, dt$$$$+\sum_{m=1}^N\,
\frac{\rho_{1,m}}{\xi_1(R+t)-\vphi_1(R+t,p_m)}\,dt.$$
Now we show that
$\vphi_{2}(I,\cdot)$ is $\F_{R,\bar t_2}$-measurable. Fix $n\in\N$.
Let $I_n=\lfloor nI\rfloor/n$. For $m\in\N\cup\{0\}$, let ${\cal
E}_n(m)=\{m/n\le I<(m+1)/n\}$. Then ${\cal E}_n(m)$ is $\F_{R,\bar
t_2}$-measurable, and $I_n=m/n$ on ${\cal E}_n(m)$. Since $m/n\le
\bar t_2$ and $I_n=m/n$ on ${\cal E}_n(m)$, so $I_n$ agrees with
$(m/n)\wedge \bar t_2$ on ${\cal E}_n(m)$. Now $(m/n)\wedge \bar
t_2$ is an $(\F^2_t)$-stopping time, and $\F^2_{(m/n)\wedge \bar
t_2}\subset\F^2_{\bar t_2}\subset \F_{R,\bar t_2} $. So
$\vphi_2((m/n)\wedge \bar t_2,\cdot)$ is $\F_{R,\bar
t_2}$-measurable. Since $\vphi_2(I_n,\cdot)=\vphi_2((m/n)\wedge \bar
t_2,\cdot)$ on ${\cal E}_n(m)$, and  ${\cal E}_n(m)$ is $\F_{R,\bar
t_2}$-measurable for each $m\in\N\cup\{0\}$, so $\vphi_2(I_n,\cdot)$
is $\F_{R,\bar t_2}$-measurable. Since
$\vphi_2(I_n,\cdot)\to\vphi_2(I,\cdot)$ as $n\to\infty$, so
$\vphi_2(I,\cdot)$ is also $\F_{R,\bar t_2}$-measurable. Thus $K_2(I)$ is
$\F_{R,\bar t_2}$-measurable. Hence for any $t\ge 0$, $\vphi_{K_1(R+t)\cup K_2(I)}$ is
$\F_{R+t,\bar t_2}$-measurable.
>From (\ref{circ=1}), $\vphi_{1,I}(R+t,\cdot)$ and $\vphi_{2,R+t}(I,\cdot)$ are
both $\F_{R+t,\bar t_2}$-measurable. If the $t_j$ and $t_k$ in
(\ref{pa-A-1-alpha}$\sim$\ref{E}) are replaced by $R+t$ and $I$, respectively,
then all these SDEs are $\F_{R+t,\bar t_2}$-adapted. From the same computation, we conclude that
$(M(R+t,I),0\le t<T_1(I)-R)$ is a continuous $(\F_{R+t,\bar t_2})_{t\ge
0}$-local martingale.
$\Box$

\vskip 3mm

Let $\HP$ denote the set of $(H_1,H_2)$ such that $H_j$ is a hull in
$\HH$ w.r.t.\ $\infty$ that contains some neighborhood of $x_j$ in
$\HH$, $j=1,2$, $\lin{H_1}\cap\lin{H_2}=\emptyset$, and $p_m\not\in
\lin{H_1}\cup\lin{H_2}$, $1\le m\le N$. For $(H_1,H_2)\in\HP$, let
$T_j(H_j)$ be the first time that $\lin{K_j(t)}\cap\lin{\HH\sem
H_j}\ne\emptyset$, $j=1,2$. An argument that is similar to the proof
of Theorem 5.1 in \cite{reversibility} gives the following.

\begin{Theorem} For any $(H_1,H_2)\in\HP$, there are $C_2>C_1>0$
depending only on $H_1$ and $H_2$ such that $M(t_1,t_2)\in[C_1,C_2]$
for any $(t_1,t_2)\in[0,T_1(H_1)]\times [0,T_2(H_2)]$. \label{bound}
\end{Theorem}

Fix $(H_1,H_2)\in\HP$. Let $\mu$ denote the joint distribution of
$(\xi_1(t):0\le t\le T_1)$ and $(\xi_2(t):0\le t\le T_2)$. From
Theorem \ref{M} and Theorem \ref{bound}, we have $$\int
M(T_1(H_1),T_2(H_2))d\mu=\EE_\mu[M(T_1(H_1),T_2(H_2))]=M(0,0)=1.$$
Note that $M(T_1(H_1),T_2(H_2))>0$. Suppose $\nu$ is a measure on
$\F^1_{T_1(H_1)}\times\F^2_{T_2(H_2)}$ such that
$d\nu/d\mu=M(T_1(H_1),T_2(H_2))$. Then $\nu$ is a probability
measure. Now suppose the joint distribution of $(\xi_1(t),0\le t\le
T_1(H_1))$ and $(\xi_2(t),0\le t\le T_2(H_2))$ is $\nu$ instead of
$\mu$. Fix an $(\F^2_t)$-stopping time $\bar t_2$ with $\bar t_2\le
T_2(H_2)$. From (\ref{SDE-xi}$\sim$\ref{M-SDE}) and Girsanov
theorem, there is an $(\F^1_t\times\F^2_{\bar t_2})$-Brownian motion
$\bar B_1(t)$  such that $\xi_1(t)$ satisfies the
$(\F^1_{t}\times\F^2_{\bar t_2})$-adapted SDE for $0\le t\le
T_1(H_1)$:
$$d\xi_1(t)=\sqrt\kappa_1 d\bar B_1(t)+\Big(\frac{6-\kappa_1}{2}\cdot \frac{A_{1,2}}
{A_{1,1}}-\frac{\kappa_1}2\cdot
\frac{A_{1,1}}{E_{1,0}}+\sum_{m=1}^N\rho_{1,m}\,
\frac{A_{1,1}}{E_{1,m}}\Big)\Big|_{(t,\bar t_2)}\,dt.$$ Let
$\xi_{1,\bar t_2}(t)=A_{1,0}(t,\bar t_2)=\vphi_{2,t}(\bar
t_2,\xi_1(t))$, $0\le t\le T_1(H_1)$. From Ito's formula and
(\ref{-3}), $\xi_{1,\bar t_2}(t)$ satisfies \BGE d\xi_{1,\bar
t_2}(t)=A_{1,1}(t,\bar t_2)\sqrt\kappa_1 d\bar
B_1(t)+\Big(-\frac{\kappa_1}2\cdot
\frac{A_{1,1}^2}{E_{1,0}}+\sum_{m=1}^N\rho_{1,m}\,
\frac{A_{1,1}^2}{E_{1,m}}\Big)\Big|_{(t,\bar
t_2)}\,dt.\label{xi-1-bar}\EDE Since $\vphi_2(\bar t_2,\cdot)$ is a conformal map, and
from (\ref{Ensemble}), for $0\le t_1<T_1(\bar t_2)$,
$$K_{1,\bar t_2}(t)=(K_1(t)\cup K_2(\bar t_2))/K_2(\bar t_2)=
\vphi_2(\bar t_2,K_1(t)),$$
so $(K_{1,\bar t_2}(t))$ is a Loewner chain. Let $v(t)=\hcap (K_{1,\bar t_2}(t))/2$.
>From Proposition \ref{chordal-Loewner-chain}, $v(t)$ is a continuous increasing function with $v(0)=0$, and
$(\til K(t)=K_{1,\bar t_2}(v^{-1}(t)))$ are chordal Loewner hulls driven
by some real continuous function, say $\til\xi(t)$, and the chordal Loewner maps
are $\til\vphi(t,\cdot)=\vphi_{K_{1,\bar t_2}(v^{-1}(t))}=\vphi_{1,\bar t_2}(v^{-1}(t),\cdot)$.
Moreover,
$$\{\til\xi(v(t))\}=\cap_{\eps>0}\lin{\til K(v(t+\eps))/ \til K(v(t))};\quad
\{\xi_1(t)\}=\cap_{\eps>0}\lin{K_1(t+\eps)/ K_1(t)}.$$
Let $W_t=\vphi_{2,t}(\bar t_2,\cdot)$. From (\ref{circ=1}), for $\eps>0$, we have
$$W_t(K_1(t+\eps)/K_1(t))=
W_t\circ \vphi_1(t,\cdot)(K_1(t+\eps)\sem K_1(t))
$$$$=\vphi_{K_1(t)\cup K_2(\bar t_2)}(K_1(t+\eps)\sem K_1(t))$$$$=
\vphi_{K_1(t)\cup K_2(\bar t_2)}((K_1(t+\eps)\cup K_2(\bar t_2))\sem (K_1(t)\cup K_2(\bar t_2)))$$
$$=(K_1(t+\eps)\cup K_2(\bar t_2))/ (K_1(t)\cup K_2(\bar t_2))$$
$$=[(K_1(t+\eps)\cup K_2(\bar t_2))/K_2(\bar t_2)]/ [(K_1(t)\cup K_2(\bar t_2))/K_2(\bar t_2)]$$
$$=K_{1,\bar t_2}(t+\eps)/K_{1,\bar t_2}(t)=\til K(v(t+\eps))/\til K(v(t)).$$
Thus $\til\xi(v(t))=W_t(\xi_1(t))=\vphi_{2,t}(\bar t_2,\xi_1(t))=\xi_{1,\bar t_2}(t)$. Since
 $\hcap(K_1(t+\eps)/K_1(t))=2\eps$ and
 $\hcap(\til K(v(t+\eps))/\til K(v(t)))=2v(t+\eps)-2v(t)$,
so from Proposition \ref{hcap}, $v'(t)=W_t'(\xi_1(t))^2=
\pa_z\vphi_{2,t}(\bar t_2,\xi_1(t))^2=A_{1,1}(t,\bar t_2)^2$.

>From the definitions of $E_{1,0}$ and $E_{1,m}$, and (\ref{circ=1}), we have \BGE
E_{1,0}(v^{-1}(t),\bar t_2)=\til\xi(t)-\til\vphi(t,\xi_2(\bar
t_2));\label{D-1-0}\EDE \BGE E_{1,m}(v^{-1}(t),\bar
t_2)=\til\xi(t)-\til\vphi(t,\vphi_2(\bar t_2,p_m)).\label{D-1-m}\EDE
>From (\ref{xi-1-bar}$\sim$\ref{D-1-m}), and the properties of $v(t)$
and $\til\xi(t)$, there is a Brownian motion $\til B_1(t)$ such that
$\til\xi(t)$, $0\le t<v(T_1(H_1))$, satisfies the SDE:
$$d\til\xi(t)=\sqrt\kappa_1 d\til B_1(t)+\Big(
\frac{-{\kappa_1}/2}{\til\xi(t)-\til\vphi(t,\xi_2(\bar
t_2))}+\sum_{m=1}^N \frac{\rho_{1,m}}{\til\xi(t)-\til\vphi(t,
\vphi_2(\bar t_2,p_m))}\Big)\,dt.$$ Note that
$\til\xi(0)=\xi_{1,\bar t_2}(0)=\vphi_2(\bar t_2,x_1)$. Thus
conditioned on $\F^2_{\bar t_2}$, $(\til K(t),0\le
t<v(T_1(H_1))$, is a partial chordal
SLE$(\kappa_1;-\frac{\kappa_1}2,\vec{\rho}_{1})$ process started
from $(\vphi_2(\bar t_2,x_1);\xi_2(\bar t_2),\vphi_2(\bar
t_2,\vec{p}))$. By symmetry, we may
exchange the subscripts $1$ and $2$ in the above statement.

\begin{Theorem} Suppose $n\in\N$ and $(H_1^m,H_2^m)\in\HP$, $1\le m\le n$.
There is a continuous function $M_*(t_1,t_2)$ defined on
$[0,\infty]^2$ that satisfies the following properties: (i) $M_*=M$
on $[0,T_1(H_1^m)]\times[0,T_2(H_2^m)]$ for $1\le m\le n$; (ii)
$M_*(t,0)=M_*(0,t)=1$ for any $t\ge 0$; (iii) $M(t_1,t_2)\in
[C_1,C_2]$ for any $t_1,t_2\ge 0$, where $C_2>C_1>0$ are constants
depending only on $H_j^m$, $j=1,2$, $1\le m\le n$; (iv) for any
$(\F^2_t)$-stopping time $\bar t_2$, $(M_*(t_1,\bar t_2),t_1\ge 0)$
is a bounded continuous $(\F^1_{t_1}\times\F^2_{\bar t_2})_{t_1\ge
0}$-martingale; and (v) for any $(\F^1_t)$-stopping time $\bar t_1$,
$(M_*(\bar t_1, t_2),t_2\ge 0)$ is a bounded continuous $(\F^1_{\bar
t_1}\times\F^2_{t_2})_{t_2\ge 0}$-martingale.
 \label{martg}
\end{Theorem}
{\bf Proof.} This is Theorem 6.1 in \cite{reversibility}. For reader's
convenience, we include the proof here.
The first quadrant $[0,\infty]^2$ will be divided by the
 vertical or horizontal lines $\{x_j=T_j(H_j^m)\}$, $1\le m\le n$,
$j=1,2$, into small rectangles, and $M_*$ will be piecewise defined on
these rectangles. Theorem \ref{bound} will be used to prove the boundedness, and
Theorem \ref{M-I} will be used to prove the martingale properties.

Let $\N_n:=\{k\in\N:k\le n\}$. Write $T_j^k$ for $T_j(H^k_j)$,
$k\in\N_n$, $j=1,2$. Let $S\subset\N_n$ be such that $\cup_{k\in
S}[0,T_1^k]\times[0,T_2^k]=\cup_{k=1}^n [0,T_1^k]\times[0,T_2^k]$,
and $\sum_{k\in S} k\le \sum_{k\in S'}k$ if $S'\subset\N_n$ also
satisfies this property. Such $S$ is a random nonempty set, and
$|S|\in\N_n$ is a random number. Define a partial order
``$\preceq$'' on $[0,\infty]^2$ such that
$(s_1,s_2)\preceq(t_1,t_2)$ iff $s_1\le t_1$ and $s_2\le t_2$. If
$(s_1,s_2)\preceq(t_1,t_2)$ and $(s_t,s_2)\ne (t_1,t_2)$,  we write
$(s_1,s_2)\prec(t_1,t_2)$. Then for each $m\in\N_n$ there is $k\in
S$ such that $(T_1^m,T_2^m)\preceq(T_1^k,T_2^k)$; and for each $k\in
S$ there is no $m\in\N_n$ such that
$(T_1^k,T_2^k)\prec(T_1^m,T_2^m)$.

There is a map $\sigma$ from $\{1,\dots,|S|\}$ onto $S$ such that if
$1\le k_1<k_2\le |S|$, then \BGE
T_1^{\sigma(k_1)}<T_1^{\sigma(k_2)}, \quad
T_2^{\sigma(k_1)}>T_2^{\sigma(k_2)}.\label{<}\EDE Define
$T_1^{\sigma(0)}=T_2^{\sigma(|S|+1)}=0$ and
$T_1^{\sigma(|S|+1)}=T_2^{\sigma(0)}=\infty$. Then (\ref{<}) still
holds for $0\le k_1<k_2\le |S|+1$.

Extend the definition of $M$ to
$[0,\infty]\times\{0\}\cup\{0\}\times[0,\infty]$ such that
$M(t,0)=M(0,t)=1$ for $t\ge 0$. Fix $(t_1,t_2)\in[0,\infty]^2$.
There are $k_1\in \N_{|S|+1}$ and $k_2\in\N_{|S|}\cup\{0\}$ such
that \BGE T^{\sigma(k_1-1)}_1\le t_1\le T^{\sigma(k_1)}_1,\quad
T^{\sigma(k_2+1)}_2\le t_2\le T^{\sigma(k_2)}_2.\label{k1k2}\EDE If
$k_1\le k_2$, let \BGE M_*(t_1,t_2)=M(t_1,t_2).\label{M*M}\EDE It
$k_1\ge k_2+1$, let \BGE
M_*(t_1,t_2)=\frac{M(T_1^{\sigma(k_2)},t_2)M(T_1^{\sigma(k_2+1)},
T_2^{\sigma(k_2+1)})\cdots
M(T_1^{\sigma(k_1-1)},T_2^{\sigma(k_1-1)})M(t_1,T_2^{\sigma(k_1)}) }
{M(T_1^{\sigma(k_2)},T_2^{\sigma(k_2+1)})\cdots
M(T_1^{\sigma(k_1-2)},T_2^{\sigma(k_1-1)})
M(T_1^{\sigma(k_1-1)},T_2^{\sigma(k_1)}) }\label{M*}\EDE In the
above formula, there are $k_1-k_2+1$ terms in the numerator, and
$k_1-k_2$ terms in the denominator. For example, if $k_1-k_2=1$,
then
$$M_*(t_1,t_2)=M(T_1^{\sigma(k_2)},t_2)M(t_1,T_2^{\sigma(k_1)})
/M(T_1^{\sigma(k_2)},T_2^{\sigma(k_1)}).$$

We need to show that $M_*(t_1,t_2)$ is well-defined. First, we show
that the $M(\cdot,\cdot)$ in (\ref{M*M}) and (\ref{M*}) are defined.
Note that $M$ is defined on
$$Z:=\bigcup_{k=0}^{|S|+1}[0,T_1^{\sigma(k)}]\times[0,T_2^{\sigma(k)}].$$
If $k_1\le k_2$ then $t_1\le T_1^{\sigma(k_1)}\le T_1^{\sigma(k_2)}$
and $t_2\le T_2^{\sigma(k_2)}$, so $(t_1,t_2)\in Z$. Thus
$M(t_1,t_2)$ in (\ref{M*M}) is defined. Now suppose $k_1\ge k_2+1$.
Since $t_2\le T_2^{\sigma(k_2)}$ and $t_1\le T_1^{\sigma(k_1)}$, so
$(T_1^{\sigma(k_2)},t_2),(t_1,T_2^{\sigma(k_1)})\in Z$. It is clear
that $(T_1^{\sigma(k)}, T_2^{\sigma(k)})\in Z$ for $k_2+1\le k\le
k_1-1$. Thus the $M(\cdot,\cdot)$ in the numerator of (\ref{M*}) are
defined. For $k_2\le k\le k_1-1$, $T_1^{\sigma(k)}\le
T_1^{\sigma(k+1)}$, so $(T_1^{\sigma(k)}, T_2^{\sigma(k+1)})\in Z$.
Thus the $M(\cdot,\cdot)$ in the denominator of (\ref{M*}) are
defined.

Second, we show that the value of $M_*(t_1,t_2)$ does not depend on
the choice of $(k_1,k_2)$ that satisfies (\ref{k1k2}). Suppose
(\ref{k1k2}) holds with $(k_1,k_2)$ replaced by $(k_1',k_2)$, and
$k_1'\ne k_1$. Then $|k_1'-k_1|=1$. We may assume $k_1'=k_1+1$. Then
$t_1=T_1^{\sigma(k_1)}$. Let $M_*'(t_1,t_2)$ denote the
$M_*(t_1,t_2)$ defined using $(k_1',k_2)$. There are three cases.
Case 1. $k_1<k_1'\le k_2$. Then from (\ref{M*M}),
$M_*'(t_1,t_2)=M(t_1,t_2)=M_*(t_1,t_2)$. Case 2. $k_1=k_2$ and
$k_1'-k_2=1$. Then $T_1^{\sigma(k_2)}=T_1^{\sigma(k_1)}=t_1$. So
from (\ref{M*M}) and (\ref{M*}),
$$M_*'(t_1,t_2)={M(T_1^{\sigma(k_2)},t_2)M(t_1,T_2^{\sigma(k_1)})}/
{M(T_1^{\sigma(k_2)},T_2^{\sigma(k_1)})}=M(t_1,t_2)=M_*(t_1,t_2).$$
Case 3. $k_1'>k_1>k_2$. From (\ref{M*}) and that
$T_1^{\sigma(k_1)}=t_1$, we have
$$M_*'(t_1,t_2)=\frac{M(T_1^{\sigma(k_2)},t_2)M(T_1^{\sigma(k_2+1)},T_2^{\sigma(k_2+1)})
\cdots
M(T_1^{\sigma(k_1)},T_2^{\sigma(k_1)})M(t_1,T_2^{\sigma(k_1+1)})}
{M(T_1^{\sigma(k_2)},T_2^{\sigma(k_2+1)})\cdots
M(T_1^{\sigma(k_1-1)},T_2^{\sigma(k_1)})
M(T_1^{\sigma(k_1)},T_2^{\sigma(k_1+1)})}$$
$$=\frac{M(T_1^{\sigma(k_2)},t_2)M(T_1^{\sigma(k_2+1)},T_2^{\sigma(k_2+1)})
\cdots M(t_1,T_2^{\sigma(k_1)})}
{M(T_1^{\sigma(k_2)},T_2^{\sigma(k_2+1)})\cdots
M(T_1^{\sigma(k_1-1)},T_2^{\sigma(k_1)})}=M_*(t_1,t_2).$$ Similarly,
if (\ref{k1k2}) holds with $(k_1,k_2)$ replaced by $(k_1,k_2')$,
then $M_*(t_1,t_2)$ defined using $(k_1,k_2')$ has the same value as
$M(t_1,t_2)$. Thus $M_*$ is well-defined.

>From the definition, it is clear that for each $k_1\in \N_{|S|+1}$
and $k_2\in\N_{|S|}\cup\{0\}$, $M_*$ is continuous on
$[T_1^{\sigma(k_1-1)},T_1^{\sigma(k_1)}]
\times[T_2^{\sigma(k_2+1)},T_1^{\sigma(k_2)}]$. Thus $M_*$ is
continuous on $[0,\infty]^2$. Let $(t_1,t_2)\in[0,\infty]^2$.
Suppose $(t_1,t_2)\in[0,T_1^m]\times[0,T_2^{m}]$ for some
$m\in\N_n$. There is $k\in \N_{|S|}$ such that
$(T_1^m,T_2^m)\preceq(T_1^{\sigma(k)},T_2^{\sigma(k)})$. Then we may
choose $k_1\le k$ and $k_2\ge k$ such that (\ref{k1k2}) holds, so
$M_*(t_1,t_2)=M(t_1,t_2)$. Thus (i) is satisfied. If $t_1=0$, we may
choose $k_1=1$ in (\ref{k1k2}). Then either $k_1\le k_2$ or $k_2=0$.
If $k_1\le k_2$ then $M_*(t_1,t_2)=M(t_1,t_2)=1$ because $t_1=0$. If
$k_2=0$, then
$$M_*(t_1,t_2)=M(T_1^{\sigma(0)},t_2)M(t_1,T_2^{\sigma(1)})
/M(T_1^{\sigma(0)},T_2^{\sigma(1)})=1$$ because
$T_1^{\sigma(0)}=t_1=0$. Similarly, $M_*(t_1,t_2)=0$ if $t_2=0$. So
(ii) is also satisfied. And (iii) follows from Theorem \ref{bound} and
the definition of $M_*$.

\vskip 3mm

Now we prove (iv). Suppose $(t_1,t_2)\in[0,\infty]^2$ and $t_2\ge
\vee_{m=1}^n T_2^m=T_2^{\sigma(1)}$. Then (\ref{k1k2}) holds with
$k_2=0$ and some $k_1\in\{1,\dots,|S|+1\}$. So $k_1\ge k_2+1$. Since
$T_1^{\sigma(k_2)}=0$ and $M(0,t)=1$ for any $t\ge 0$, so from
(\ref{M*}) we have
$$M_*(t_1,t_2)=\frac{M(T_1^{\sigma(k_2+1)},
T_2^{\sigma(k_2+1)})\cdots
M(T_1^{\sigma(k_1-1)},T_2^{\sigma(k_1-1)})M(t_1,T_2^{\sigma(k_1)}) }
{M(T_1^{\sigma(k_2+1)},T_2^{\sigma(k_2+2)})\cdots
M(T_1^{\sigma(k_1-1)},T_2^{\sigma(k_1)}) }.$$ The right-hand side of
the above equality has no $t_2$. So
$M_*(t_1,t_2)=M_*(t_1,\vee_{m=1}^n T_2^m)$ for any $t_2\ge
\vee_{m=1}^n T_2^m$. Similarly, $M_*(t_1,t_2)=M_*(\vee_{m=1}^n
T_1^m,t_2)$ for any $t_1\ge \vee_{m=1}^n T_1^m$.

 Fix an $(\F^2_t)$-stopping time $\bar t_2$. Since
$M_*(\cdot,\bar t_2)=M_*(\cdot,\bar t_2\wedge(\vee_{m=1}^n T_2^m))$,
and $\bar t_2\wedge(\vee_{m=1}^n T_2^m)$ is also an
$(\F^2_t)$-stopping time, so we may assume that $\bar t_2\le
\vee_{m=1}^n T_2^m$. Let $I_0=\bar t_2$. For $s\in \N\cup\{0\}$,
define \BGE R_s=\sup\{T_1^m:m\in\N_n,T_2^m\ge I_s\};\quad
I_{s+1}=\sup\{T_2^m:m\in\N_n,T_2^m<I_s,T_1^m>R_s\}.\label{RI}\EDE
Here we set $\sup(\emptyset)=0$. Then we have a non-decreasing
sequence $(R_s)$ and a non-increasing sequence $(I_s)$. Let $S$ and
$\sigma(k)$, $0\le k\le |S|+1$, be as in the definition of $M_*$.
>From the property of $S$, for any $s\in\N\cup\{0\}$, \BGE
R_s=\sup\{T_1^k:k\in S,T_2^k\ge I_s\}.\label{R}\EDE Suppose for some
$s\in\N\cup\{0\}$, there is $m\in\N_n$ that satisfies $T_2^m<I_s$
and $T_1^m>R_s$. Then there is $k\in S$ such that $T_j^k\ge T_j^m$,
$j=1,2$. If $T_2^k\ge I_s$, then from (\ref{R}) we have $R_s\ge
T_1^k\ge T_1^m$, which contradicts that $T_1^m>R_s$. Thus
$T_2^k<I_s$. Now $T_2^k<I_s$, $T_1^k\ge T_1^m>R_s$, and $T_2^k\ge
T_2^m$. Thus for any $s\in\N\cup\{0\}$, \BGE
I_{s+1}=\sup\{T_2^k:k\in S,T_2^k<I_s,T_1^k>R_s\}.\label{I}\EDE

First suppose $\bar t_2>0$. Since $\bar t_2\le \vee_{m=1}^n
T_2^m=T_2^{\sigma(0)}$, so there is a unique $k_2\in \N_{|S|}$ such
that $T_2^{\sigma(k_2)}\ge \bar t_2>T_2^{\sigma(k_2+1)}$. From
(\ref{R}) and (\ref{I}), we have $R_s=T_1^{\sigma(k_2+s)}$ for $0\le
s\le |S|-k_2$; $R_s=T_1^{\sigma(|S|)}$ for $s\ge |S|-k_2$;
$I_s=T_2^{\sigma(k_2+s)}$ for $1\le s\le |S|-k_2$; and $I_s=0$ for
$s\ge |S|-k_2+1$. Since $R_0= T_1^{\sigma(k_2)}$ and $\bar t_2\le
T_2^{\sigma(k_2)}$, so from (i),
 \BGE M_*(t_1,\bar
t_2)=M(t_1,\bar t_2),\quad\mbox{for } t_1\in[0,R_0].\label{0,R0}\EDE
Suppose $t_1\in [R_{s-1},R_{s}]$ for some $s\in\N_{|S|-k_2}$. Let
$k_1=k_2+s$. Then $T_1^{\sigma(k_1-1)}\le t_1\le T_1^{\sigma(k_1)}$.
Since $I_s=T_2^{\sigma(k_2+s)}=T_2^{\sigma(k_1)}$, so from
(\ref{M*}), \BGE M_*(t_1,\bar t_2)/M_*(R_{s-1},\bar
t_2)=M(t_1,I_{s})/M(R_{s-1}, I_{s}),\quad\mbox{for }t_1\in
[R_{s-1},R_{s}].\label{M/M-1}\EDE Note that if $s\ge |S|-k_2+1$,
(\ref{M/M-1}) still holds because $R_s=R_{s-1}$. Suppose $t_1\ge
R_{n}$. Since $n\ge |S|-k_2$, so $R_n=T_1^{\sigma(|S|)}=\vee_{m=1}^n
T_1^m$. From the discussion at the beginning of the proof of (iv),
we have \BGE M_*(t_1,\bar t_2)=M_*(R_n,\bar t_2),\quad\mbox{for
}t_1\in[R_n,\infty].\label{M=M}\EDE If $\bar t_2=0$,
(\ref{0,R0}$\sim$\ref{M=M}) still hold because all $I_s=0$ and so
$M_*(t_1,\bar t_2)=M(t_1,I_s)=M(t_1,0)=1$ for any $t_1\ge 0$.

Let $R_{-1}=0$. We claim that for each $s\in\N\cup\{0\}$, $R_s$ is
an $(\F^1_t\times\F^2_{\bar t_2})_{t\ge 0}$-stopping time, and
$I_{s}$ is $\F_{R_{s-1},\bar t_2}$-measurable. Recall that
$\F_{R_{s-1},\bar t_2}$ is the $\sigma$-algebra obtained from the
filtration $(\F^1_t\times\F^2_{\bar t_2})_{t\ge 0}$ and its stopping
time $R_{s-1}$. It is clear that $R_{-1}=0$ is an
$(\F^1_t\times\F^2_{\bar t_2})_{t\ge 0}$-stopping time, and
$I_0=\bar t_2$ is $\F_{R_{-1},\bar t_2}$-measurable. Now suppose
$I_{s}$ is $\F_{R_{s-1},\bar t_2}$-measurable. Since $I_s\le \bar
t_2$ and $R_{s-1}\le R_s$, so for any $t\ge 0$, $\{R_s\le
t\}=\{R_{s-1}\le t\}\cap {\cal E}_t$, where
$${\cal E}_t=\bigcap_{m=1}^n(\{T_2^m<I_s\}\cup\{T_1^m\le t\})=
 \bigcap_{m=1}^n(\cup_{q\in\Q}(\{T_2^m<q\le \bar t_2\}\cap\{q<I_s\})\cup \{T_1^m\le t\} ).$$
Thus ${\cal E}_t\in \F_{R_{s-1},\bar t_2}\vee
(\F^1_t\times\F^2_{\bar t_2})$, and so $\{R_s\le t\}\in
\F^1_t\times\F^2_{\bar t_2}$ for any $t\ge 0$. Therefore $R_s$ is an
$(\F^1_t\times\F^2_{\bar t_2})_{t\ge 0}$-stopping time. Next we
consider $I_{s+1}$. For any $h\ge 0$, $$\{I_{s+1}>
h\}=\cup_{m=1}^n(\{h<T_2^m<I_s\}\cap\{T_1^m>R_s\})$$
$$=\cup_{m=1}^n(\cup_{q\in\Q}(\{h<T_2^m<q< \bar t_2\}\cap\{q<I_s\})
\cap\{T_1^m>R_s\})\in \F_{R_{s},\bar t_2}.$$ Thus $I_{s+1}$ is
$\F_{R_{s},\bar t_2}$-measurable. So the claim is proved by
induction.

Since $\bar t_2\le\vee_{m=1}^n T_2^m<T_2$, so from Theorem
\ref{M-I}, for any $s\in\N_n$, $(M(R_{s-1}+t,I_s),0\le
t<T_1(I_s)-R_{s-1})$ is a continuous $(\F_{R_{s-1}+t,\bar
t_2})_{t\ge 0}$-local martingale. For $m\in\N_n$, if $T_2^m\ge I_s$,
then $T_1^m<T_1(T_2^m)\le T_1(I_s)$. So from (\ref{RI}) we have
$R_s<T_1(I_s)$. From (\ref{M/M-1}), we find that
 $(M_*(R_{s-1}+t,\bar t_2),0\le t\le R_s-R_{s-1})$ is
a continuous $(\F_{R_{s-1}+t,\bar t_2})_{t\ge 0}$-local martingale
for any $s\in\N_n$. From Theorem \ref{M} and (\ref{0,R0}),
$(M_*(t,\bar t_2),0\le t\le R_0)$ is a continuous $(\F_{t,\bar
t_2})_{t\ge 0}$-local martingale. From (\ref{M=M}), $(M_*(R_n+t,\bar
t_2),t\ge 0)$ is a continuous $(\F_{R_n+t,\bar t_2})_{t\ge 0}$-local
martingale. Thus $(M_*(t,\bar t_2),t\ge 0)$ is a continuous
$(\F_{t,\bar t_2})_{t\ge 0}$-local martingale. Since  by (iii)
$M_*(t_1,t_2)\in[C_1, C_2]$, so this local martingale is a bounded
martingale. Thus (iv) is satisfied. Finally, (v) follows from the
symmetry in the definition (\ref{M*M})
 and (\ref{M*}) of $M_*$. $\Box$

\subsection{Coupling measures} \label{subsec-coup}
 Let ${\cal C}:=\cup_{T\in(0,\infty]}
C([0,T))$. The map $T:{\cal C}\to (0,\infty]$ is such that
$[0,T(\xi))$ is the definition domain of $\xi$. For
$t\in[0,\infty)$, let $\F_t$ be the $\sigma$-algebra on $\cal C$
generated by $\{T>s,\xi(s)\in A\}$, where $A$ is a Borel set on $\R$
and $s\in[0,t]$. Then $(\F_t)$ is a
filtration on $\cal C$, and $T$ is an $(\F_t)$-stopping time. Let $\F_\infty=\vee_t \F_t$.

For $\xi\in\cal C$, let $K_\xi(t)$, $0\le t<T(\xi)$, denote the
chordal Loewner hulls driven by $\xi$.
Let $H$ be a hull in $\HH$ w.r.t.\ $\infty$. Let
$T_H(\xi)\in[0,T(\xi)]$ be the maximal number such that
$K_\xi(t)\cap\lin{\HH\sem H}=\emptyset$ for $0\le t<T_H$. Then $T_H$
is an $(\F_t)$-stopping time. Let ${\cal C}_H=\{T_H>0\}$. Then
$\xi\in{\cal C}_H$ iff $H$ contains some neighborhood of $\xi(0)$ in
$\HH$. Define $P_H:{\cal C}_H\to {\cal C}$ such that $P_H(\xi)$ is
the restriction of $\xi$ to $[0,T_H(\xi))$. Then $P_H({\cal
C}_H)=\{T_H=T\}$, and $P_H\circ P_H=P_H$. Let ${\cal C}_{H,\pa}$
denote the set of $\xi\in\{T_H=T\}$ such that $\lin{\cup_{0\le
t<T(\xi)} K_\xi(t)}\cap\lin{(\HH\sem H)}\ne\emptyset$. If
$\xi\in{\cal C}_H\cap\{T_H<T\}$, then $P_H(\xi)\in{\cal C}_{H,\pa}$.
If $A$ is a Borel set on $\R$ and $s\in[0,\infty)$, then
$$P_H^{-1}(\{\xi\in{\cal C}:T(\xi)>s,\xi(s)\in A\})=\{\xi\in{\cal
C}_H:T_H(\xi)>s,\xi(s)\in A\}\in\F_{T_H^-}.$$ Thus $P_H$ is
$(\F_{T_H^-},\F_{\infty})$-measurable on ${\cal C}_H$. On the other
hand, the restriction of $\F_{T_H^-}$ to ${\cal C}_H$ is the
$\sigma$-algebra generated by $\{\xi\in{\cal
C}_H:T_H(\xi)>s,\xi(s)\in A\}$, where $s\in[0,\infty)$ and $A$ is a
Borel set on $\R$. Thus $P_H^{-1}(\F_\infty)$ agrees with the
restriction of $\F_{T_H^-}$ to ${\cal C}_H$.

Let $\ha\C=\C\cup\{\infty\}$ be the Riemann sphere with spherical
metric. Let $\Gamma_{\ha\C}$ denote the space of nonempty compact
subsets of $\ha\C$ endowed with Hausdorff metric. Then
$\Gamma_{\ha\C}$ is a compact metric space. Define $G:{\cal C}\to
\Gamma_{\ha\C}$ such that $G(\xi)$ is the spherical closure of
$\{t+i\xi(t):0\le t<T(\xi)\}$. Then $G$ is a one-to-one map. Let
$I_G=G({\cal C})$. Let $\F^H_{I_G}$ denote the $\sigma$-algebra on
$I_G$ generated by Hausdorff metric. Let $${\cal
R}=\{\{z\in\C:a<\Ree z<b,c<\Imm z< d\}:a,b,c,d\in\R\}.$$ Then
$\F^H_{I_G}$ agrees with the $\sigma$-algebra on $I_G$ generated by
$\{\{F\in I_G:F\cap R\ne\emptyset\}:R\in{\cal R}\}$. Using this
result, one may check that $G$ and $G^{-1}$ (defined on $I_G$) are
both measurable with respect to $\F_\infty$ and $\F^H_{I_G}$.

For $j=1,2$, let $\xi_j(t)$, $0\le t<T_j$, be the maximal solution
to (\ref{SDE-xi}). Then $\xi_j$ is a $\cal C$-valued random
variable, and $T(\xi_j)=T_j$. Since $B_1(t)$ and $B_2(t)$ are
independent, so are $\xi_1(t)$ and $\xi_2(t)$. Now we write $K_j(t)$
for $K_{\xi_j}(t)$, $0\le t<T_j$, $j=1,2$. For $j=1,2$, let $\mu_j$
denote the distribution of $\xi_j$, which is a probability measure
on $\cal C$. Let $\mu=\mu_1\times\mu_2$ be a probability measure on
${\cal C}^2$. Then $\mu$ is the joint distribution of $\xi_1$ and
$\xi_2$. Let $(H_1,H_2)\in\HP$. For $j=1,2$, $H_j$ contains some
neighborhood of $x_j=\xi_j(0)$ in $\HH$, so $\xi_j\in {\cal
C}_{H_j}$. Since $\cup_{0\le t<T_j} K_j(t)$ disconnects some force
point from $\infty$, so we do not have $\cup_{0\le t<T_j}
K_j(t)\subset H_j$, which implies that $T_{H_j}(\xi_j)<T_j$,
$j=1,2$. Thus $P_{H_j}(\xi_j)\in {\cal C}_{H_j,\pa}$, and so
$(P_{H_1}\times P_{H_2})_*(\mu)$ is supported by ${\cal
C}_{H_1,\pa}\times {\cal C}_{H_2,\pa}$.

Let $\HP_*$ be the set of $(H_1,H_2)\in\HP$ such that for $j=1,2$,
$H_j$ is a polygon whose vertices have rational coordinates. Then
$\HP_*$ is countable. Let $(H_1^m,H_2^m)$, $k\in\N$, be an
enumeration of $\HP_*$. For each $n\in\N$, let $M_*^n(t_1,t_2)$ be
the $M_*(t_1,t_2)$ given by Theorem \ref{martg} for $(H_1^m,H_2^m)$,
$1\le m\le n$, in the above enumeration. For each $n\in\N$ define
$\nu^n=(\nu^n_1,\nu^n_2)$ such that
${d\nu^n}/{d\mu}=M_*^n(\infty,\infty)$. From Theorem \ref{martg},
$M_*^n(\infty,\infty)>0$ and $\int
M_*^n(\infty,\infty)d\mu=\EE_\mu[M_*^n(\infty,\infty)]=1$, so
$\nu^n$ is a probability measure on ${\cal C}^2$. Since
$d\nu^n_1/d\mu_1=\EE_\mu[M_*^n(\infty,\infty)|\F^2_\infty]=M_*^n(\infty,0)=1$,
so $\nu^n_1=\mu_1$. Similarly, $\nu^n_2=\mu_2$. So each $\nu^n$ is a
coupling of $\mu_1$ and $\mu_2$.

Let $\bar\nu^n=(G\times G)_*(\nu^n)$ be a probability measure on
$\Gamma_{\ha \C}^2$. Since $\Gamma_{\ha\C}^2$ is compact, so
$(\bar\nu^n)$ has a subsequence $(\bar\nu^{n_k})$ that converges
weakly to some probability measure $\bar\nu=(\bar\nu_1,\bar\nu_2)$
on $\Gamma_{\ha\C}\times \Gamma_{\ha\C}$. Then for $j=1,2$,
$\bar\nu^{n_k}_j\to\bar\nu_j$ weakly. For $n\in\N$ and $j=1,2$,
since $\nu^n_j=\mu_j$, so $\bar \nu^n_j=G_*(\mu_j)$. Thus
$\bar\nu_j=G_*(\mu_j)$, $j=1,2$. So $\bar\nu$ is supported by
$I_G^2$. Let $\nu=(\nu_1,\nu_2)=(G^{-1}\times G^{-1})_*(\bar\nu)$ be
a probability measure on ${\cal C}^2$. Here we use the fact that
$G^{-1}$ is $(\F^H_{I_G},\F^j_\infty)$-measurable. For $j=1,2$, we
have $\nu_j=(G^{-1})_*(\bar\nu_j)=\mu_j$. So $\nu$ is also a
coupling measure of $\mu_1$ and $\mu_2$.

\begin{Lemma} For any $r\in\N$, the restriction of $\nu$ to
$\F^1_{T_{H_1^r}}\times \F^2_{T_{H_2^r}}$ is absolutely continuous
w.r.t.\ $\mu$, and the Radon-Nikodym derivative is
$M(T_{H_1^r}(\xi_1),T_{H_2^r}(\xi_2))$. \label{Radon}
\end{Lemma}
{\bf Proof.} We may choose $s\in\N$ such that
$\lin{H_j^r}\cap\lin{\HH\sem H_j^s}=\emptyset$, $j=1,2$. Since $M$
is continuous, so $M(T_{H_1^s}(\xi_1),T_{H_2^s}(\xi_2))$ is
$\F^1_{T_{H_1^s}^-}\times \F^2_{T_{H_2^s}^-}$-measurable. Let
$\nu_{(s)}$ be defined on $\F^1_{T_{H_1^s}^-}\times
\F^2_{T_{H_2^s}^-}$ such that
$d\nu_{(s)}/d\mu=M(T_{H_1^s}(\xi_1),T_{H_2^s}(\xi_2))$. From Theorem
\ref{M} and Theorem \ref{bound}, $\nu_{(s)}$ is a probability
measure on $\F^1_{T_{H_1^s}^-}\times \F^2_{T_{H_2^s}^-}$. Let
$\mr\nu_{(s)}=(P_{H_1^s}\times P_{H_2^s})_*(\nu_{(s)})$. Since
$P_{H_j^s}$ is $(\F^j_{T_{H_j^s}^-},\F^j_\infty)$-measurable,
$j=1,2$, so $\mr\nu_{(s)}$ is a probability measure on ${\cal C}^2$,
and is absolute continuous w.r.t.\ $(P_{H_1^s}\times
P_{H_2^s})_*(\mu)$. Let
 $\bar\nu_{(s)}=(G\times G)_*(\mr
\nu_{(s)})$. Since $G$ is $(\F^j_\infty,\F^H_{I_G})$-measurable,
$j=1,2$, so $\bar\nu_{(s)}$ is a probability measure on $I_G^2$.
Since $d\nu^{n_k}/d\mu=M^{n_k}_*(\infty,\infty)$, and
$M^{n_k}_*(\cdot,\cdot)$ satisfies the martingale properties, so the
Radon-Nikodym derivative of the restriction of $\nu^{n_k}$ to
$\F^1_{T_{H_1^s}^-}\times \F^2_{T_{H_2^s}^-}$ w.r.t.\ $\mu$ is
$M^{n_k}_*(T_{H_1^s}(\xi_1),T_{H_2^s}(\xi_2))$. If $n_k\ge s$ then
$M^{n_k}_*(T_{H_1^s}(\xi_1),T_{H_2^s}(\xi_2))=
M(T_{H_1^s}(\xi_1),T_{H_2^s}(\xi_2))$. Thus the restriction of
$\nu^{n_k}$ to $\F^1_{T_{H_1^s}^-}\times \F^2_{T_{H_2^s}^-}$ equals
to $\nu_{(s)}$, which implies that
$$(G\times G)_*\circ (P_{H_1^s}\times
P_{H_2^s})_*(\nu^{n_k})=(G\times G)_*\circ (P_{H_1^s}\times
P_{H_2^s})_*(\nu_{(s)})=\bar\nu_{(s)}.$$

For $n\in\N$, let a ${\cal C}^2$-valued random variable
$(\zeta^{n}_1,\zeta^{n}_2)$ have the distribution $\nu^{n}$, and
$\eta^{n}_j=P_{H_j^s}(\zeta^n_j)$, $j=1,2$. Let $\bar\tau^n_{(s)}$
denote the distribution of
$(G(\zeta^n_1),G(\zeta^n_2),G(\eta^n_1),G(\eta^n_2))$. Then
$\bar\tau^n_{(s)}$ is supported by $\Xi$, which is the set of
$(L_1,L_2, F_1,F_2)\in\Gamma_{\ha C}^4$ such that $F_j\subset L_j$,
$j=1,2$. It is easy to check that $\Xi$ is a closed subset of
$\Gamma_{\ha\C}^4$. Then $(n_k)$ has a subsequence $(n_k')$ such
that $(\bar\tau^{n_k'}_{(s)})$ converges weakly to some probability
measure $\bar\tau_{(s)}$ on $\Xi$. Since the marginal of
$\bar\tau^{n_k'}_{(s)}$ at the first two variables equals to
$(G\times G)_*(\nu^{n_k'})=\bar\nu^{n_k'}$, and
$\bar\nu^{n_k'}\to\bar\nu$ weakly, so the marginal of
$\bar\tau_{(s)}$ at the first two variables equals to $\bar\nu$.
Since the marginal of $\bar\tau^{n_k'}_{(s)}$ at the last two
variables equals to $(G\times G)_*\circ (P_{H_1^s}\times
P_{H_2^s})_*(\nu^{n_k'})=\bar\nu_{(s)}$ when $n_k'\ge s$, so the
marginal of $\bar\tau_{(s)}$ at the last two variables equals to
$\bar\nu_{(s)}$.

Now $\bar\tau_{(s)}$ is supported by $I_G^4$. Let
$\tau_{(s)}=(G\times G\times G\times G)^{-1}_*(\bar \tau_{(s)})$.
Let a ${\cal C}^4$-valued random variable
$(\zeta_1,\zeta_2,\eta_1,\eta_2)$ have distribution $\tau_{(s)}$.
Since $\bar\nu=(G\times G)_*(\nu)$ and $\bar\nu_{(s)}=(G\times
G)_*(\mr \nu_{(s)})$, so the distribution of $(\zeta_1,\zeta_2)$ is
$\nu$, and the distribution of $(\eta_1,\eta_2)$ is $\mr\nu_{(s)}$.
For $j=1,2$, since $G(\eta_j)\subset G(\zeta_j)$, so $\eta_j$ is
some restriction of $\zeta_j$. Note that for $j=1,2$, $K_j(t)$ does
not always stay in $H_j^s$, so $\mu_j$ is supported by
$\{T_{H_j^s}<T_j\}$, so $(P_{H_j^s})_*(\mu_j)$ is supported by
${\cal C}_{H_j^s,\pa}$. Thus $(P_{H_1^s}\times P_{H_2^s})_*(\mu)$ is
supported by ${\cal C}_{H_1^s,\pa}\times {\cal C}_{H_2^s,\pa} $.
Since $\mr\nu_{(s)}$ is absolutely continuous w.r.t.\
$(P_{H_1^s}\times P_{H_2^s})_*(\mu)$, so $\mr\nu_{(s)}$ is also
supported by ${\cal C}_{H_1^s,\pa}\times{\cal C}_{H_2^s,\pa}$. Thus
for $j=1,2$, $K_{\eta_j}(t)\cap \lin{\HH\sem H_j^s}=\emptyset$ for
$0\le t<T(\eta_j)$, and $\lin{\cup_{0\le
t<T(\eta_j)}K_{\eta_j}(t)}\cap\lin{(\HH\sem H_j^s)}\ne\emptyset$.
Since $\eta_j$ is a restriction of $\zeta_j$, so from the above
observation, we have $\eta_j=P_{H_j^s}(\zeta_j)$, $j=1,2$. Thus
$\mr\nu_{(s)}=(P_{H_1^s}\times P_{H_2^s})_*(\nu)$.

We now have $(P_{H_1^s}\times
P_{H_2^s})_*(\nu)=\mr\nu_{(s)}=(P_{H_1^s}\times
P_{H_2^s})_*(\nu_{(s)})$. So $\nu({\cal E})=\nu_{(s)}({\cal E})$ for
any ${\cal E}\in P_{H_1^s}^{-1}(\F^1_\infty)\times
P_{H_2^s}^{-1}(F^2_\infty)$. Since $P_{H_j^s}^{-1}(\F^j_\infty)$
agrees with the restriction of $\F^j_{T_{H_j^s}^-}$ to ${\cal
C}_{H_j^s}$, $j=1,2$, and both $\nu$ and $\nu_{(s)}$ are supported
by ${\cal C}_{H_1^s}\times {\cal C}_{H_2^s}$, so the restriction of
$\nu$ to $\F^1_{T_{H_1^s}^-}\times \F^2_{T_{H_2^s}^-}$ equals to
$\nu_{(s)}$. From the definition of $\nu_{(s)}$, the Radon-Nikodym
derivative of the restriction of $\nu$ to $\F^1_{T_{H_1^s}^-}\times
\F^2_{T_{H_2^s}^-}$ w.r.t.\ $\mu$ is
$M(T_{H_1^s}(\xi_1),T_{H_2^s}(\xi_2))$.

For $j=1,2$, since $\lin{H_j^r}\cap\lin{\HH\sem H_j^s}=\emptyset$,
so $\mu_j$ and $\nu_j$ are supported by $\{T_{H_j^r}<T_{H_j^s}\}$.
Since  $\F^j_{T_{H_j^r}}\subset \F^j_{T_{H_j^s}^-}$ on
$\{T_{H_j^r}<T_{H_j^s}\}$, $j=1,2$, so the restriction of $\nu$ to
$\F^1_{T_{H_1^r}}\times \F^2_{T_{H_2^r}}$ is absolutely continuous
w.r.t.\ $\mu$, and the Radon-Nikodym derivative equals to
$$\EE_\mu[M(T_{H_1^s}(\xi_1),T_{H_2^s}(\xi_2))|\F^1_{T_{H_1^r}}\times
\F^2_{T_{H_2^r}}]=M(T_{H_1^r}(\xi_1),T_{H_2^s}(\xi_r)).\quad\Box$$

\vskip 4mm

\no {\bf Proof of Theorem \ref{bartk}.} Now let the ${\cal
C}^2$-valued random variable $(\xi_1,\xi_2)$ have distribution $\nu$
in the above theorem. Let $K_j(t)$ and $\vphi_j(t,\cdot)$, $0\le
t<T_j$, be the chordal Loewner hulls and maps, respectively, driven
by $\xi_j$, $j=1,2$. For $j\ne k\in\{1,2\}$, since $\nu_j=\mu_j$, so $K_j(t)$,
$0\le t<T_j$, is a chordal
SLE$(\kappa_j;-\frac{\kappa_j}2,\vec{\rho}_{j})$ process started
from $(x_j;x_{k},\vec{p})$.

Fix $j\ne k\in\{1,2\}$. Suppose $\bar t_k$ is an $(\F^k_t)$-stopping
time with $\bar t_k<T_k$. For $n\in\N$, define
$$R_n=\sup\{T_j(H^m_j):1\le m\le n, T_k(H^m_k)\ge \bar t_k\}.$$
Here we set $\sup(\emptyset)=0$. Then for any $t\ge 0$,
$$\{R_n\le t\}=\cap_{m=1}^n(\{\bar t_k>
T_k(H^m_k)\}\cup\{T_j(H^m_j)\le t\}\in\F^j_t\times\F^k_{\bar t_k}.$$
So $R_n$ is an $(\F^j_t\times\F^k_{\bar t_k})_{t\ge 0}$-stopping
time for each $n\in\N$. For $1\le m\le n$, let $\bar t_k^m=\bar
t_k\wedge T_k(H^m_2)$.  Then $\bar t_k^m$ is an $(\F^k_t)$-stopping
time, and $\bar t_k^m\le T_k(H^m_k)$. Let $(L(t))$ be a chordal
SLE$(\kappa_j;-\frac{\kappa_j}2,\vec{\rho}_{j})$ process started
from $(\vphi_k(\bar t_2,x_j);\xi_k(\bar t_k),\vphi_k(\bar
t_k,\vec{p}))$. From Lemma \ref{Radon} and the discussion after
Theorem \ref{bound}, $\vphi_k(\bar t_k^m,K_j(t))$, $0\le t\le
T_j(H_j^m)$, has the distribution of a time-change of a partial
$(L(t))$, i.e., $(L(t))$ stopped at some stopping time. Let ${\cal E}_{n,m}=\{\bar t_k\le
T_k(H^m_k)\}\cap\{R_n=T_j(H^m_j)\}$. Since
$\{R_n>0\}=\cup_{m=1}^n{\cal E}_{n,m}$, and on each ${\cal
E}_{n,m}$, $\bar t_k=\bar t_k^m$ and $R_n=T_j(H^m_j)$, so
$\vphi_k(\bar t_k,K_j(t))$, $0\le t\le R_n$, has the distribution of
a time-change of a partial $(L(t))$. Since $T_j(\bar
t_k)=\sup\{T_j(H_j^m):m\in\N,T_k(H_k^m)\ge \bar
t_k\}=\vee_{n=1}^\infty R_n$, so $\vphi_k(\bar t_k,K_j(t))$, $0\le
t<T_j(\bar t_k)$, has the distribution of a time-change of a partial
$(L(t))$. Thus after a time-change,
$\vphi_k(\bar t_k,K_j(t))$, $0\le t<T_j(\bar t_k)$, is a partial
chordal SLE$(\kappa_j;-\frac{\kappa_j}2,\vec{\rho}_{j})$ process started
from $(\vphi_k(\bar t_2,x_j);\xi_k(\bar t_k),\vphi_k(\bar
t_k,\vec{p}))$. $\Box$

\subsection{Coupling in degenerate cases}
Now we will prove that Theorem \ref{bartk} still holds if one or
more than one force points $p_m$ are degenerate, i.e., $p_m$ equals to
$x_1^\pm$ or $x_2^\pm$. The results do not immediately follow from
Theorem \ref{bartk} in the generic case. However, we may modify the
proof of Theorem \ref{bartk} to deal with the degenerate cases. We
need to find some suitable two-dimensional local martingales, and
obtain some boundedness.

We use the following simplest example to illustrate the idea.
Suppose there is only one degenerate force point, which is
$p_1=x_1^+$. Then the $(K_1(t))$ and $(K_2(t))$ in Theorem
\ref{bartk} should be understood as follows: $(K_1(t))$ is a chordal
SLE$(\kappa_1;-\frac{\kappa_1}2,\vec{\rho}_{1})$ process started
from $(x_1;x_2,x_1^+,p_2,\dots,p_N)$, and $(K_2(t))$ is a chordal
SLE$(\kappa_2;-\frac{\kappa_2}2+\rho_{2,1},\rho_{2,2},\dots,\rho_{2,N})$
process started from $(x_2;x_1,p_2,\dots,p_N)$. Here the force
points $x_1$ and $p_1=x_1^+$ for $(K_2(t))$ are combined to be a
single force point $x_1$. And in Theorem \ref{bartk},
$\vphi_2(t_2,p_1)=\vphi_2(t_2,x_1^+)$ should be understood as
$\vphi_2(t_2,x_1)$; and $\vphi_1(t_1,p_1)=\vphi_1(t_1,x_1^+)$ should
be understood as $p_1(t_1)$, which is a component of the solution to
the equation that generates $(K_1(t))$.

We want to define $M(t_1,t_2)$ by  (\ref{def-M}) and
(\ref{def-tilM}). However, for the case we study here, some factors
in (\ref{def-tilM}) does not make sense, and some factors become
zero, which will cause trouble in (\ref{def-M}). Let's check the
factors in (\ref{def-tilM}) one by one. Let $j\ne k\in\{1,2\}$.
First, $A_{j,h}(t_1,t_2)=\pa_z^h \vphi_{k,t_j}(t_k,\xi_j(t_j))$ is
well defined for $h=0,1$, and $A_{j,1}$ is a positive number; and
$E(t_1,t_2)=|A_{1,0}(t_1,t_2)-A_{2,0}(t_1,t_2)|>0$ is well defined.
Then $F(t_1,t_2)>0$ is well defined by (\ref{def-E}). Now
$B_{m,0}(t_1,t_2)=\vphi_{K_1(t_1)\cup K_2(t_2)}(p_m)$ is well
defined for each $1\le m\le n$. For the degenerate force point
$p_1=x_1^+$, the formula $\vphi_{K_1(t_1)\cup K_2(t_2)}(x_1^+)$ is
understood as
$\vphi_{2,t_1}(t_2,\vphi_1(t_1,x_1^+))=\vphi_{2,t_1}(p_1(t_1))$. So
$E_{j,m}=A_{j,0}-B_{m,0}$ and $C_{m_1,m_2}=B_{m_1,0}-B_{m_2,0}$ are
all well defined. Among these numbers, $|C_{m_1,m_2}(t_1,t_2)|$ is
positive if $m_1\ne m_2$, and $|E_{j,m}(t_1,t_2)|$ is positive
except when $j=1$, $m=1$ and $t_1=0$. The factor
$B_{m,1}(t_1,t_2)=\pa_z\vphi_{K_1(t_1)\cup K_2(t_2)}(p_m)$ is well
defined and positive except when $m=1$. Now for $(t_1,t_2)\in{\cal
D}$, define $$\til
N(t_1,t_2)=\frac{A_{1,1}^{\alpha_1}A_{2,1}^{\alpha_2}}{E^{1/2}}\,F^\lambda
|E_{2,1}|^{\rho^*_{2,1}}\,\prod_{m=2}^N
\Big(B_{m,1}^{\gamma_m}\prod_{j=1}^2|E_{j,m}|^{\rho^*_{j,m}}\Big)
\prod_{1\le m_1<m_2\le N} |C_{m_1,m_2}|^{\delta_{m_1,m_2}};$$ and
\BGE N(t_1,t_2)=(\til N(t_1,t_2)\til N(0,0))/(\til N(t_1,0)\til
N(0,t_2)).\label{def-N}\EDE Then in the generic case, we have
$M(t_1,t_2)/N(t_1,t_2)=L_1(t_1,t_2)/L_2(t_1,t_2)$, where
$$L_1(t_1,t_2)=\frac{\pa_z \vphi_{K_1(t_1)\cup
K_2(t_2)}(p_1)^{\gamma_1}|\vphi_{2,t_1}(t_2,\xi_1(t_1))-\vphi_{2,t_1}(t_2,\vphi_1(t_1,p_1))|^{\rho^*_{1,1}}}
{\pa_z \vphi_{K_1(t_1)\cup
K_2(0)}(p_1)^{\gamma_1}|\vphi_{2,t_1}(0,\xi_1(t_1))-\vphi_{2,t_1}(0,\vphi_1(t_1,p_1))|^{\rho^*_{1,1}}}$$
$$=\pa_z\vphi_{2,t_1}(t_2,\vphi_1(t_1,p_1))^{\gamma_1}\,
\frac{|\vphi_{2,t_1}(t_2,\xi_1(t_1))-\vphi_{2,t_1}(t_2,\vphi_1(t_1,p_1))|^{\rho^*_{1,1}}}
{|\xi_1(t_1)-\vphi_1(t_1,p_1)|^{\rho^*_{1,1}}},$$
$$L_2(t_1,t_2)=\frac{\pa_z \vphi_{K_1(0)\cup
K_2(t_2)}(p_1)^{\gamma_1}|\vphi_{2,0}(t_2,\xi_1(0))-\vphi_{2,0}(t_2,\vphi_1(0,p_1))|^{\rho^*_{1,1}}}
{\pa_z \vphi_{K_1(0)\cup
K_2(0)}(p_1)^{\gamma_1}|\vphi_{2,0}(0,\xi_1(0))-\vphi_{2,0}(0,\vphi_1(0,p_1))|^{\rho^*_{1,1}}}$$
$$= \pa_z\vphi_{2}(t_2,p_1)^{\gamma_1}\,
\frac{|\vphi_{2}(t_2,x_1)-\vphi_{2}(t_2,p_1)|^{\rho^*_{1,1}}}
{|x_1-p_1|^{\rho^*_{1,1}}}.$$ In the above equalities,
(\ref{circ=1}) is used. Thus in the generic case,
$$\frac{M(t_1,t_2)}{N(t_1,t_2)}=\left(\frac{\pa_z\vphi_{2,t_1}(t_2,\vphi_1(t_1,p_1))}
{\pa_z\vphi_{2}(t_2,p_1)}\right)^{\gamma_1}\cdot
\left(\frac{|x_1-p_1|}{|\vphi_{2}(t_2,x_1)-\vphi_{2}(t_2,p_1)|}\right)^{\rho^*_{1,1}}$$
$$\cdot\left( \frac{
|\vphi_{2,t_1}(t_2,\xi_1(t_1))-\vphi_{2,t_1}(t_2,\vphi_1(t_1,p_1))|
}{ |\xi_1(t_1)-\vphi_1(t_1,p_1)|}\right)^{\rho^*_{1,1}}.$$

Now come back to the degenerate case $p_1=x_1^+$ we are studying
here. Then
$${\pa_z\vphi_{2,t_1}(t_2,\vphi_1(t_1,p_1))}={\pa_z\vphi_{2,t_1}(t_2,\vphi_1(t_1,x_1^+))}
\quad\mbox{and} \quad
{\pa_z\vphi_{2}(t_2,p_1)}={\pa_z\vphi_{2}(t_2,x_1)}$$ still make
sense and are both positive. If $t_1>0$, then
$|\vphi_{2,t_1}(t_2,\xi_1(t_1))-\vphi_{2,t_1}(t_2,\vphi_1(t_1,p_1))|$
and $|\xi_1(t_1)-\vphi_1(t_1,p_1)|$ both make sense and are
positive. And we have
$$\lim_{t_1\to 0^+}|\vphi_{2,t_1}(t_2,\xi_1(t_1))-\vphi_{2,t_1}(t_2,\vphi_1(t_1,p_1))|/
|\xi_1(t_1)-\vphi_1(t_1,p_1)|=\pa_z \vphi_{2,t_1}(t_2,\xi_1(t_1)).$$
Since $p_1=x_1^+$, we may view
$|x_1-p_1|/|\vphi_{2}(t_2,x_1)-\vphi_{2}(t_2,p_1)|$ as
$$\lim_{p\to x_1^+}|x_1-p|/|\vphi_{2}(t_2,x_1)-\vphi_{2}(t_2,p)|=
1/\pa_z\vphi_2(t_2,x_1).$$

These observations suggest us to define $M(t_1,t_2)$ in the case
$p_1=x_1^+$ as follows. For $(t_1,t_2)\in\cal D$, define
$U(t_1,t_2)$ such that $U(0,t_2)=\pa_z
\vphi_{2,t_1}(t_2,\xi_1(t_1))$; and if $t_1>0$, then
$$U(t_1,t_2)=|\vphi_{2,t_1}(t_2,\xi_1(t_1))-\vphi_{2,t_1}(t_2,\vphi_1(t_1,p_1))|/
|\xi_1(t_1)-\vphi_1(t_1,p_1)|.$$ Then $U$ is continuous on ${\cal
D}$. Now for $(t_1,t_2)\in\cal M$, define \BGE
M(t_1,t_2)=N(t_1,t_2)\cdot\frac{U(t_1,t_2)^{\rho^*_{1,1}}}
{U(0,t_2)^{\rho^*_{1,1}}}\cdot
\frac{\pa_z\vphi_{2,t_1}(t_2,\vphi_1(t_1,x_1^+))^{\gamma_1}}
{\pa_z\vphi_{2}(t_2,x_1)^{\gamma_1}}.\label{def-M*}\EDE Then $M$ is
continuous on ${\cal D}$. It is direct to check that
$M(t_1,0)=M(0,t_2)=1$ for any $t_1\in[0,T_1)$ and $t_2\in[0,T_2)$.

Suppose $(\xi_1(t),0\le t<T_1)$ and $(\xi_2(t),0\le t<T_2)$ are
independent. Let $\mu_j$ denote the distribution of $(\xi_j(t))$,
$j=1,2$, and $\mu=\mu_1\times\mu_2$. Let $(\F^j_t)$ be the
filtration generated by $(\xi_j(t))$, $j=1,2$. Let $j\ne
k\in\{1,2\}$. Then for any fixed $(\F^k_t)$-stopping time $\bar t_k$
with $\bar t_k<T_k$, the process $M|_{t_j=t,t_k=\bar t_k}$, $0\le
t<T_j(\bar t_k)$, is an $(\F^j_t\times \F^k_{\bar t_k})$-adapted
local martingale, under the probability measure $\mu$. The argument
is similar to that used in Section \ref{subsec-mart}.

Let $\HP$ denote the set of $(H_1,H_2)$ such that $H_j$ is a hull in
$\HH$ w.r.t.\ $\infty$ that contains some neighborhood of $x_j$ in
$\HH$, $j=1,2$, $\lin{H_1}\cap\lin{H_2}=\emptyset$, and $p_m\not\in
\lin{H_1}\cup\lin{H_2}$, $2\le m\le N$. Here we only require that
the non-degenerate force points are bounded away from $H_1$ and
$H_2$. Then Theorem \ref{bound} still holds here. For the proof, one
may check that Theorem \ref{bound} holds with $M(t_1,t_2)$ replaced
by $N(t_1,t_2)$, $U(t_1,t_2)$,
$\pa_z\vphi_{2,t_1}(t_2,\vphi_1(t_1,x_1^+))$, and
$\pa_z\vphi_{2}(t_2,x_1)$, respectively. So for any
$(H_1,H_2)\in\HP$, $\EE_\mu[M(T_1(H_1),T_2(H_2))]=1$. Suppose $\nu$
is a measure on $\F^1_{T_1(H_1)}\times\F^2_{T_2(H_2)}$ such that
$d\nu/d\mu=M(T_1(H_1),T_2(H_2))$. Then $\nu$ is a probability
measure. Now suppose the joint distribution of $(\xi_1(t),0\le t\le
T_1(H_1))$ and $(\xi_2(t),0\le t\le T_2(H_2))$ is $\nu$ instead of
$\mu$. Let $j\ne k\in\{1,2\}$. Using Girsanov Theorem, one may
check that for any fixed $(\F^k_t)$-stopping time $\bar t_k$ with
$\bar t_k\le T_k(H_k)$. conditioned on $\F^k_{\bar t_k}$,
 $(\vphi_k(\bar t_k,K_j(t)))$,
$0\le t<T_j(H_j)$, is a time-change of a partial chordal
SLE$(\kappa_j;-\frac{\kappa_j}2,\vec{\rho}_{j})$ process started
from $(\vphi_k(\bar t_k,x_j);\xi_k(\bar t_k),\vphi_k(\bar
t_k,\vec{p}))$. We now can use the argument in Section
\ref{subsec-coup} to derive Theorem \ref{bartk} in this degenerate
case.

\section{Applications} \label{appli}
\subsection{Duality}
We say $\alpha$ is a crosscut in $\HH$ on $\R$ if $\alpha$ is a
simple curve that lies inside $\HH$ except for the two ends of
$\alpha$, which lie on $\R$. If $\alpha$ is a crosscut, then
$\HH\sem\alpha$ has two connected components: one is bounded, the
other is unbounded. Let $D(\alpha)$ denote the bounded component. We
say that such $\alpha$ strictly encloses some $S\subset\lin{\HH}$ if
$\lin{S}\subset\lin{D(\alpha)}$ and $\lin{S}\cap\alpha=\emptyset$.

In Theorem \ref{bartk}, let $\kappa_1<4<\kappa_2$; $x_1<x_2$; $N=3$;
$p_1\in(-\infty,x_1)$, $p_2\in(x_2,\infty)$, $p_3\in (x_1,x_2)$; for
$j=1,2$, $\rho_{j,1}=C_1(\kappa_j-4)$, $\rho_{j,2}=C_2(\kappa_j-4)$,
and $\rho_{j,3}=\frac 12(\kappa_j-4)$ for some $C_1\le 1/2$ and
$C_2=1-C_1$. Let $K_j(t)$, $0\le t<T_j$, $j=1,2$, be given by
Theorem \ref{bartk}. Let $\vphi_j(t,\cdot)$ and $\beta_j(t)$, $0\le
t<T_j$, $j=1,2$, be the corresponding Loewner maps and traces.

Let $K_2(T^-_2)=\cup_{0\le t<T_2} K_2(t)$. Since $\kappa_1\in(0,4)$,
so $\beta_1(t)$, $0\le t<T_j$, is a simple curve, and $\beta_1(t)\in
\HH$ for $0<t<T_j$. From Theorem \ref{lim-1p3} and Lemma
\ref{coordinate2}, a.s.\ $\beta_1(T_1):=\lim_{t\to T_1} \beta_1(t)$
exists and lies on $(x_2,p_2)$. For simplicity, we use $\beta_1$ to
denote the image $\{\beta_1(t):0\le t\le T_1\}$. Thus $\beta_1$ is a
crosscut in $\HH$ on $\R$.

Suppose $S\subset \HH$ is bounded. Then there is a unique unbounded
component of $\HH\sem\lin{S}$, which is denoted by $D_\infty$. Then
we call $\pa D_\infty\cap\HH$ the outer boundary of $S$ in $\HH$.
Let it be denoted by $\pa^{\out}_\HH S$.

\begin{Lemma} Almost surely
$\beta_1=\pa^{\out}_{\HH}K_2(T^-_2)$. \label{K-infty}
\end{Lemma}
{\bf Proof.} For $j=1,2$, let ${\cal P}_j$ denote the set of
polygonal crosscuts in $\HH$ on $\R$ whose vertices have rational
coordinates, which strictly enclose $x_j$, and which do not contain
or enclose $x_{3-j}$ or $p_m$, $m=1,2,3$. For each $\gamma\in{\cal
P}_j$, let $T_j(\gamma)$ denote the first time that $\beta_j$ hits
$\gamma$. Then $T_j(\gamma)$ is an $(\F^j_t)$-stopping time, and
$T_j(\gamma)<T_j$. Moreover, we have $T_j=\vee_{\gamma\in{\cal P}_j}
T_j(\gamma)$. Let ${\cal P}^*_2$ denote the set of polygonal
crosscuts in $\HH$ on $\R$ whose vertices have rational coordinates,
and which strictly enclose $x_2$.

We first show that $K_2(T^-_2)\subset D(\beta_1)\cup \beta_1$ a.s..
Let $\cal E$ denote the event that $\beta_2$ intersects
$\HH\sem(D(\beta_1)\cup\beta_1)$. We need to show that $\PP[{\cal
E}]=0$.  For $\alpha\in{\cal P}^*_2$ and $\gamma\in{\cal P}_2$, let
${\cal E}_{\alpha;\gamma}$ denote the event that $\alpha$ strictly
encloses $\beta_1$, and $\beta_2$ hits $\alpha$ before $\gamma$.
Then ${\cal E}=\cup_{\alpha\in{\cal P}^*_2;\gamma\in{\cal P}_2}{\cal
E}_{\alpha;\gamma}$. Since ${\cal P}^*_2$ and ${\cal P}_2$ are
countable, so we suffice to show $\PP[{\cal E}_{\alpha;\gamma}]=0$
for any $\alpha\in{\cal P}^*_2$ and $\gamma\in{\cal P}_2$.

Now fix $\alpha\in{\cal P}^*_2$ and $\gamma\in{\cal P}_2$. Let $\bar
t_2$ denote the first time that $\beta_2$ hits $\alpha\cup\gamma$.
Then $\bar t_2$ is an $(\F^2_t)$-stopping time, and $\bar t_2\le
T_2(\gamma)<T_2$. From Theorem \ref{bartk}, after a time-change,
$\vphi_2(\bar t_2,\beta_1(t))$, $0\le t<T_1(\bar t_2)$, has the same
distribution as a full chordal SLE$(\kappa_1;-\frac{\kappa_1}2,
C_1(\kappa_1-4),C_2(\kappa_1-4),\frac 12(\kappa_1-4))$ trace started
from $(\vphi_2(\bar t_2,x_1);\xi_2(\bar t_2),\\\vphi_2(\bar
t_2,p_1),\vphi_2(\bar t_2,p_2),\vphi_2(\bar t_2,p_3))$. Here we have
$$\vphi_2(\bar t_2,p_1)<\vphi_2(\bar t_2,x_1)<\vphi_2(\bar
t_2,p_3)<\xi_2(\bar t_2)<\vphi_2(\bar t_2,p_2).$$ Since
$C_1(\kappa_1-4)\ge \kappa_1/2-2$, $\frac 12(\kappa_1-4))\ge
\kappa_1/2-2$, and
$$|(C_1(\kappa_1-4)+C_2(\kappa_1-4))-(\frac12(\kappa_1-4)+(-\frac{\kappa_1}2))|=|\kappa_1-2|<2,$$
so from Theorem \ref{lim-1p3} and Lemma \ref{coordinate2}, a.s.\
$\lim_{t\to T_1(\bar t_2)} \vphi_2(\bar
t_2,\beta_1(t))\in(\xi_2(\bar t_2),\vphi_2(\bar t_2,p_2))$. Thus
a.s.\ $\{\vphi_2(\bar t_2,\beta_1(t)):0\le t<T_1(\bar t_2)\}$
disconnects $\xi_2(\bar t_2)$ from $\infty$ in $\HH$. So a.s.\
$\beta_1$ disconnects $\beta_2(\bar t_2)$ from $\infty$ in $\HH\sem
K_2(\bar t_2)$.

Assume that the event ${\cal E}_{\alpha;\gamma}$ occurs. Since $\beta_2$
starts from $x_2$, which is strictly enclosed by $\alpha$, so
$\beta_2(t)\in\lin{D(\alpha)}$ for $0\le t\le \bar t_2$, which
implies that $K_2(\bar t_2)\subset \lin{D(\alpha)}$. On the other
hand, $\beta_1$ is strictly enclosed by $\alpha$, and $\beta_2(\bar
t_2)\in\alpha$. Thus $\beta_1$ does not disconnect $\beta_2(\bar
t_2)$ from $\infty$ in $\HH\sem K_2(\bar t_2)$. So we have
$\PP[{\cal E}_{\alpha;\gamma}]=0$. Thus  $K_2(T^-_2)\subset
D(\beta_1)\cup \beta_1$ a.s..

Next we show that a.s.\ $\beta_1\subset\lin{K_2(T^-_2)}$. Fix
$\gamma\in{\cal P}_1$ and $q\in\Q_{\ge 0}$. Let $\bar t_1=q\wedge
T_1(\gamma)$. Then $\bar t_1$ is an $(\F^1_t)$-stopping time, and
$\bar t_1\le T_1(\gamma)<T_1$. From Theorem \ref{bartk}, after a
time-change, $\vphi_1(\bar t_1,\beta_2(t))$, $0\le t<T_2(\bar t_1)$,
has the same distribution as a full chordal
SLE$(\kappa_2;-\frac{\kappa_2}2,
C_1(\kappa_2-4),C_2(\kappa_2-4),\frac 12(\kappa_2-4))$ trace started
from $(\vphi_1(\bar t_1,x_2);\xi_1(\bar t_1),\vphi_1(\bar
t_1,p_1),\vphi_1(\bar t_1,p_2),\vphi_1(\bar t_1,p_3))$. Here we have
$$\vphi_1(\bar
t_1,p_1)<\xi_1(\bar t_1)<\vphi_1(\bar t_1,p_3)<\vphi_1(\bar
t_1,x_2)<\vphi_1(\bar t_1,p_2).$$ Since $\frac12(\kappa_2-4)\ge
\kappa_2/2-2$, $C_2(\kappa_2-4)\ge \kappa/2-2$, and
$C_2(\kappa_2-4)+C_1(\kappa_2-4)=\kappa_2-4\ge\kappa_2/2-2$, so from
Theorem \ref{more force} and Lemma \ref{coordinate}, a.s.\
$\xi_1(\bar t_1)$ is a subsequential limit point of $\vphi_1(\bar
t_1,\beta_2(t))$ as $t\to T_2(\bar t_1)$. Thus a.s.\ $\beta_1(\bar
t_1)$ is a subsequential limit point of $\beta_2(t)$ as $t\to
T_2(\bar t_1)$. So $\beta_1(\bar t_1)\in\lin{K_2(T_2(\bar
t_1)^-)}\subset \lin{K_2(T_2^-)}$ a.s.. Since $\Q_{\ge 0}$ is
countable, so a.s.\ $\beta_1(q\wedge T_1(\gamma))\in
\lin{K_2(T_2^-)}$ for any $q\in\Q_{\ge 0}$. Since $\Q_{\ge 0}$ is
dense in $\R_{\ge 0}$, so a.s.\ $\beta_1(t)\in \lin{K_2(T_2^-)}$ for
any $t\in[0,T_1(\gamma)]$. Since ${\cal P}_1$ is countable and
$T_1=\vee_{\gamma\in{\cal P}_1} T_1(\gamma)$, so almost surely
$\beta_1(t)\in \lin{K_2(T_2^-)}$ for any $t\in[0,T_1)$, i.e.,
$\beta_1\subset\lin{K_2(T^-_2)}$ a.s.. Finally, $K_2(T^-_2)\subset
D(\beta_1)\cup \beta_1$ and $\beta_1\subset\lin{K_2(T^-_2)}$ imply
that $\beta_1=\pa^{\out}_{\HH}K_2(T^-_2)$. $\Box$

\begin{Theorem} Suppose $\kappa>4$;
 $p_1<x_1<p_3<x_2<p_2$; $C_1\le 1/2\le C_2$ and $C_1+C_2=1$. Let
 $K(t)$, $0\le t<T$, be chordal
 SLE$(\kappa;-\frac{\kappa}2,C_1(\kappa-4),C_2(\kappa-4),\frac12(\kappa-4))$
 process
 started from $(x_2;x_1,p_1,p_2,p_3)$. Let $K(T^-)=\cup_{0\le t<T} K(t)$.
 Then a.s.\ $K(T^-)$ is bounded, and $\pa^{\out}_\HH K(T^-)$
  has the distribution of the image of a chordal
 SLE$(\kappa';-\frac{\kappa'}2,C_1(\kappa'-4),C_2(\kappa'-4),\frac12(\kappa'-4))$
 trace started from $(x_1;x_2,p_1,p_2,p_3)$, where
 $\kappa'=16/\kappa$. \label{K-infty*}
\end{Theorem}

Theorem \ref{K-infty*} still holds if we let $p_1\in(-\infty,x_1)$,
or $=-\infty$, or $=x_1^-$; let $p_2\in(x_2,\infty)$, or $=\infty$,
or $=x_2^+$; and let $p_3\in(x_1,x_2)$, or $=x_1^+$, or $=x_2^-$. In
some cases we may use Theorem \ref{lim-1p}, Theorem \ref{lim-1p*},
or Theorem \ref{3pt**} instead of Theorem \ref{4pt*} to prove that
$\beta_1$ is a crosscut. We may derive some interesting theorems
from some  cases.

\begin{Theorem} Suppose $\kappa\ge 8$, and $K(t)$, $0\le t<\infty$,
is a standard chordal SLE$(\kappa)$ process, i.e., the chordal
Loewner chain driven by $\xi(t)=\sqrt\kappa B(t)$. Let $x\in\R\sem\{0\}$ and
$T_x$ be the first $t$ such that $x\in\lin{K(t)}$. Then $\pa
K(T_x)\cap\HH$ has the same distribution as the image of a chordal
SLE$(\kappa';-\frac{\kappa'}2,-\frac{\kappa'}2,\frac{\kappa'}2-2)$
trace started from $(x;0,x^a,x^b)$, where $\kappa'=16/\kappa$, $a=\sign(x)$
and $b=\sign(-x)$, and so $\pa K(T_x)\cap\HH$ is a crosscut in $\HH$ on $\R$ connecting $x$
with some $y\in\R\sem\{0\}$ with $\sign(y)=\sign(-x)$. 
\label{kappa>8}
\end{Theorem}
{\bf Proof.} $K(t)$, $0\le t<T_x$, is a full chordal SLE$(\kappa;0)$
process started from $(0;x)$. Since $\kappa\ge 8$, so
$K(T_x)=\cup_{0\le t<T_x}K_t$ and $\pa
K(T_x)\cap\HH=\pa^{\out}_{\HH}K(T_x)$. If $x<0$, this follows from
Theorem \ref{K-infty*} with
 $x_1=x$, $x_2=0$, $p_1=x_1^-$, $p_2=\infty$, $p_3=x_1^+$;
 $C_1=2/(\kappa-4)$ and
$C_2=1-C_1$. If $x>0$, this follows from symmetry. $\Box$

\vskip 3mm One may expect that after reasonable modifications, the
above theorem also holds for $\kappa\in(4,8)$. In this case, for the
SLE$(\kappa';-\frac{\kappa'}2,-\frac{\kappa'}2,\frac{\kappa'}2-2)$
trace started from $(x;0,x^a,x^b)$, the force $-\frac{\kappa'}2$
that corresponds to the degenerate force point $x^a$ does not
satisfy $-\frac{\kappa'}2\ge \kappa'/2-2$. So we must allow that the
process continue growing after the degenerate force point is
swallowed. This will make sense because $-\frac{\kappa'}2>-2$.

\begin{Corollary}
For $\kappa>8$, chordal SLE$(\kappa)$ trace is not reversible.
\end{Corollary}
Proof. Let $\beta(t)$, $0\le t < \infty$, be a standard chordal SLE$(\kappa)$ trace.
Let $W(z) = 1/\lin{z}$ and
$\gamma(t) = W(\beta(1/t))$. Suppose chordal SLE$(\kappa)$ trace is reversible, then after a time-change,
$(\gamma(t), 0 < t < \infty)$ has the same distribution as $(\beta(t), 0 < t < \infty)$. Let T be the first
time such that $1\in \beta(t)$. Since $\kappa > 8$, $1$ is visited by $\beta$ exactly once a.s..
Thus $1/T$ is the first time such that −$1\in \gamma(t)$. From the above theorem, $\pa(\beta((0, T]))\cap\HH$
and $\pa (\gamma((0, 1/T]))\cap \HH$ both have the distribution of the image of a chordal SLE$(\kappa';-\frac{\kappa'}2,-\frac{\kappa'}2,\frac{\kappa'}2-2)$ trace
started from $(1; 0, 1^+−, 1^-)$, where $\kappa'′ = 16/\kappa$. From Lemma 2.1
 and the definition of $\gamma$, we find that
$\pa(\beta([T,\infty)))\cap\HH$ has the distribution of the image of a chordal
SLE$(\kappa';\frac{3\kappa'}2-4,\frac{\kappa'}2-2,-\frac{\kappa'}2)$ trace
started from $(1; 0, 1^+−, 1^-)$. Since $\kappa'<2$, so $-\frac{\kappa'}2\ne \frac{3\kappa'}2-4$.
 Thus
$\pa(\beta((0, T]))\cap\HH$ and $\pa(\beta([T,\infty)))\cap \HH$ have different distributions. However, since
$\beta$ is space-filling and never crosses its past, the two boundary curves coincide, which gives
a contradiction. $\Box$

\vskip 3mm Suppose $S\subset\HH$ and
$\lin{S}\cap[a,\infty)=\emptyset$ for some $a\in\R$. Then there is a
unique component of $\HH\sem\lin{S}$, which has $[a,\infty)$ as part
of its boundary. Let $D_+$ denote this component. Then $\pa D_+\cap
\HH$ is called the right boundary of $S$ in $\HH$. Let it be denoted
by $\pa_\HH^+ S$.

\begin{Theorem} Let $\kappa>4$, $C\ge 1/2$,
 and $K(t)$, $0\le t<\infty$, be a
chordal SLE$(\kappa;C(\kappa-4),\frac 12(\kappa-4))$ process started
from $(0;0^+,0^-)$. Let $K(\infty)=\cup_{t<\infty}K(t)$. Let
$W(z)=1/\lin{z}$. Then $W(\pa_\HH^+ K(\infty))$ has the same distribution
as the image of a chordal SLE$(\kappa';C'(\kappa'-4))$ trace started
from $(0;0^+)$, where $\kappa'=16/\kappa$ and $C'=1-C$.
\label{duality1}
\end{Theorem}
{\bf Proof.} Let $W_0(z)=1/(1-z)$. Then $W_0$ is a conformal
automorphism of $\HH$, and $W_0(0)=1$, $W_0(\infty)=0$, $W_0(0^\pm)=1^\pm$.
>From Lemma \ref{coordinate}, after a time-change, $(W_0(K(t)))$ has the same
distribution as a chordal SLE$(\kappa;C'(\kappa-4)-\frac\kappa
2,C(\kappa-4),\frac 12(\kappa-4))$ process started from
$(1;0,1^+,1^-)$. Applying Theorem \ref{K-infty*} with $x_1=0$,
$x_2=1$, $p_1=0^-$, $p_2=1^+$, $p_3=1^-$, $C_1=C'$ and $C_2=C$, we
find that $\pa_\HH^{\out} W_0(K_\infty)$ has the same distribution
as the image of a chordal SLE$(\kappa';
C(\kappa'-4)-2,C'(\kappa'-4))$ trace started from $(0;1,0^-)$. Let
$\beta$ denote this trace. From Lemma \ref{coordinate2} and Theorem
\ref{lim-1p*}, $\beta$ is a crosscut in $\HH$ from $0$ to $1$. Thus
$\pa_{\HH}^{+} K_\infty= W_0^{-1}(\beta)$, and so $W(\pa_{\HH}^{+}
K_\infty)=W\circ W_0^{-1}(\beta)$. Let $W_1=W\circ W_0^{-1}$. Then
$W_1(z)=\lin{z}/(\lin{z}-1)$. So $W_1(0)=0$,
$W_1(1)=\infty$, $W_1(0^-)=0^+$. From Lemma \ref{coordinate}, after a time-change,
$W_1(\beta)$ has the same distribution as the image of a chordal
SLE$(\kappa';C'(\kappa'-4))$ trace started from $(0;0^+)$. $\Box$

\begin{Theorem} Let $\kappa>4$, $C\ge 1/2$,
 and $K(t)$, $0\le t<\infty$, be a
chordal SLE$(\kappa;C(\kappa-4))$ process started from $(0;0^+)$.
Let $K(\infty)=\cup_{t<\infty}K(t)$. Let $W(z)=1/\lin z$. Then
$W(\pa_\HH^+ K(\infty))$ has the same distribution as the image of a
chordal SLE$(\kappa';C'(\kappa'-4),\frac 12(\kappa'-4))$ trace
started from $(0;0^+,0^-)$, where $\kappa'=16/\kappa$ and $C'=1-C$.
\label{duality2}
\end{Theorem}
{\bf Proof.} This proof is similar to the previous one. We use the
same $W_0$, $W_1$, $x_1$, $x_2$, $p_1$, and $p_2$, except that now
$p_3=0^+$ instead of $1^-$. And the $\beta$ here is a chordal
SLE$(\kappa'; C(\kappa'-4)-\frac{\kappa'}2,\frac
12(\kappa'-4),C'(\kappa'-4))$ trace started from $(0;1,0^+,0^-)$.
$\Box$

\begin{Corollary} Let $\kappa>4$,
 and $K(t)$, $0\le t<\infty$, be a
chordal SLE$(\kappa;\kappa-4,\frac 12(\kappa-4))$ process started
from $(0;0^+,0^-)$. Let $K(\infty)=\cup_{t<\infty}K(t)$. Then
$\pa_\HH^+ K(\infty)$ has the same distribution as the image of a
standard chordal SLE$(\kappa')$ trace, where $\kappa'=16/\kappa$.
\end{Corollary}
{\bf Proof.} This follows from Theorem \ref{duality1} and the
reversibility of chordal SLE$(\kappa')$ trace when $\kappa'\in(0,4]$
(see \cite{reversibility}). $\Box$

\vskip 3mm If we assume that Conjecture \ref{conjec} is true, then
in Theorem \ref{duality1} we conclude that $\pa_\HH^+ K(\infty)$ has
the same distribution as a chordal SLE$(\kappa';C'(\kappa'-4))$
trace started from $(0;0^+)$; and in Theorem \ref{duality2} we
conclude that $\pa_\HH^+ K(\infty)$ has the same distribution as the
image of a chordal SLE$(\kappa';C'(\kappa'-4),\frac 12(\kappa'-4))$
trace started from $(0;0^+,0^-)$, where $\kappa'=16/\kappa$ and
$C'=1-C$. Moreover, assuming Conjecture \ref{conjec}, and letting
$C=1$ in Theorem \ref{duality2}, we conclude that the right boundary
of the final hull of a chordal SLE$(\kappa;\kappa-4)$ process
started from $(0;0^+)$ has the same distribution as the image of a
chordal SLE$(\kappa';\frac 12(\kappa-4))$ trace started from
$(0;0^-)$, which is Conjecture 2 in \cite{Julien-Duality}. Moreover,
we conjecture that for $C_r,C_l\ge 1/2$, if $(K(t))$ is a chordal
SLE$(\kappa;C_r(\kappa-4),C_l(\kappa-4))$ started from
$(0;0^+,0^-)$, then $\pa_\HH^+ K(\infty)$ has the same distribution
as the image of a chordal
SLE$(\kappa';C_r'(\kappa'-4),C_l'(\kappa'-4))$ trace started from
$(0;0^+,0^-)$, where  $C_r'=1-C_r$ and $C_l'=1/2-C_l$.

\subsection{Reversibility}
\begin{Theorem} Let $\vec{p}_\pm=(p_{\pm
1},\dots,p_{\pm N_\pm})$ and $\vec{\rho}_\pm=(\rho_{\pm
1},\dots,\rho_{\pm N_\pm})$, where $0<\pm p_{\pm m}<\pm p_{\pm n}$ for $1\le m<n\le
N_\pm$; $\sum_{m=1}^n \rho_{\pm m}\ge 0$ for $1\le n\le N_\pm$, and
$\sum_{m=1}^{N_\pm}\rho_{\pm m}=0$.  Let $\beta(t)$, $0\le t<\infty$, be a
chordal SLE$(4;\vec{\rho}_+,\vec{\rho}_{-})$ trace started from $(0;
\vec{p}_+,\vec{p}_-)$. Let $W(z)=1/\lin{z}$. Then a.s.\ $\lim_{t\to
\infty}\beta(t)=\infty$, and after a time-change, the reversal of $(W(\beta(t)))$ has the
same distribution as a chordal
SLE$(4;-\vec{\rho}_+,-\vec{\rho}_{-})$ trace started from $(0;
W(\vec{p}_+),W(\vec{p}_-)$, where
$W(\vec{p}_\pm)=(W(p_{\pm 1}),\dots,W(p_{\pm N_\pm}))$.
 \label{kappa=4}
\end{Theorem}
{\bf Proof.} Choose $x_0>p_{N_+}$. Let $W_0(z)=x_0/(x_0-z)$. Then
$W_0$ maps $\HH$ conformally onto $\HH$, and $W_0(0)=1$,
$W_0(\infty)=0$. Let $q_{\pm j}=W_0(p_{\pm j})$, $1\le j\le N_\pm$.
Then $0<q_{-N_-}<\dots <q_{-1}<
1<q_1 <\dots<q_{N_+}$. Let $x_1=1$, $x_2=0$,
$\vec{\rho}_{1,\pm}=\vec{\rho}_\pm$, and
$\vec{\rho}_{2,\pm}=-\vec{\rho}_\pm$. From Theorem \ref{bartk},
there is a coupling of two curves $\beta_j(t)$, $0\le t<T_j$,
$j=1,2$, such that for fixed $j\ne k\in\{1,2\}$, (i) $(\beta_j(t))$
is a chordal SLE$(4;-2,\vec{\rho}_{j,+},\vec{\rho}_{j,-})$ trace
started from $(x_j;x_{k},\vec{p}_+,\vec{p}_-)$; and (ii) for any
$(\F^k_t)$-stopping time $\bar t_k$ with $\bar t_k<T_k$,
$\vphi_k(\bar t_k,\beta_j(t))$, $0\le t<T_j(\bar t_k)$, has the same
distribution as a chordal
SLE$(4;-2,\vec{\rho}_{j,+},\vec{\rho}_{j,-})$ trace started from
$(\vphi_k(\bar t_k,x_j);\xi_k(\bar t_k),\vphi_k(\bar
t_k,\vec{p}_+),\vphi_k(\bar t_k,\vec{p}_-))$, where
$\vphi_j(t,\cdot)$ and $\xi_j(t)$, $0\le t<T_j$, are chordal Loewner
maps and driving function for the trace $\beta_j$, $j=1,2$.
Note the symmetry between $\vec{\rho}_{1,\pm}$ and
$\vec{\rho}_{2,\pm}$: $\sum_{m=1}^n \rho_{1,\pm m}\ge 0$ for all
$1\le n\le N_\pm$, and $\sum_{m=1}^{N_\pm}\rho_{1,\pm m}=0$;
$\sum_{m=n}^{N_\pm} \rho_{2,\pm m}\ge 0$ for all $1\le n\le N_\pm$,
and $\sum_{m=1}^{N_\pm}\rho_{2,\pm m}=0$.

Fix  $j\ne k\in\{1,2\}$. From Lemma \ref{coordinate} and Theorem
\ref{more force}, we have a.s.\ $x_k\in\lin{\beta_j((0,T_j))}$. Now
fix an $(\F^k_t)$-stopping time $\bar t_k\in(0,T_k)$. From Lemma
\ref{coordinate} and Theorem \ref{more force}, we have a.s.\ $
\lin{\vphi_k(\bar t_k,\beta_j((0,T_j(\bar
t_k))))}\cap\R=\{\xi_k(\bar t_k)\}$, which implies that
$\lin{\beta_j((0,T_j(\bar t_k)))}\cap(\R\cup\beta_k((0,\bar
t_k)))=\{\beta_k(\bar t_k)\}$. Since $\bar t_k>0$, so $\beta_k(\bar
t_k)\ne \beta_k(0)=x_k$. If $T_j(\bar t_k)=T_j$, then
$x_k\in\lin{\beta_j((0,T_j(\bar t_k)))}$, which a.s.\ does not
happen.
Thus a.s.\ $T_j(\bar t_k)<T_j$. So we have a.s.\
$\beta_j(T_j(\bar t_k))=\lim_{t\to T_j(\bar t_k)^-}\beta_j(t)\in \lin{\beta_j((0,T_j(\bar t_k)))}$.
 From the definition of
$T_j(\bar t_k)$, we have a.s.\ $\beta_j(T_j(\bar
t_k))\in\beta_k([0,\bar t_k])$. Thus a.s.\ $\beta_j(T_j(\bar
t_k))=\beta_k(\bar t_k)$.

We may choose a sequence of $(\F^k_t)$-stopping times $(\bar
t^{(n)}_k)$ on $(0,T_k)$ such that $\{\bar t^{(n)}_k:n\in\N\}$ is
dense on $[0,T_k]$. Then a.s. $\beta_k(\bar
t_k^{(n)})=\beta_j(T_j(\bar t^{(n)}_k))$ for any $n\in\N$. From the
denseness of $\{\bar t^{(n)}_k:n\in\N\}$ and the continuity of
$\beta_j$ and $\beta_k$, we have a.s.\
$\beta_k((0,T_k))\subset\beta_j((0,T_j))$. Similarly, a.s.\
$\beta_j((0,T_j))\subset\beta_k((0,T_k))$. So a.s.\ $\beta_2$ is a
time-change of the reversal of $\beta_1$.

>From Lemma \ref{coordinate}, $(W_0(\beta(t)))$ has the same
distribution as $(\beta_1(t))$ after a time-change. Thus the
reversal of $(W(\beta(t)))$ has the same distribution as $(W\circ
W_0^{-1}(\beta_2(t)))$ after a time-change. From Lemma
\ref{coordinate}, $(W\circ W_0^{-1}(\beta_2(t)))$ has the same
distribution as a chordal SLE$(4;-\vec{\rho}_+,-\vec{\rho}_{-})$
trace started from $(0; W(\vec{p}_+),W(\vec{p}_-)$. $\Box$

\vskip 3mm

This theorem may also be proved using the convergence of discrete
Gaussian free field on some triangle lattice with suitable boundary
conditions (see \cite{SS}). It also holds in the degenerate cases,
i.e., $p_1=0^+$ and/or $p_{-1}=0^-$ and/or $p_{N_+}=+\infty$ and/or
$p_{-N_-}=-\infty$. For example, let $\rho_+,\rho_-\ge 0$, and apply
Theorem \ref{kappa=4} with $N_+=N_-=2$, $p_1=0^+$, $p_{-1}=0^-$,
$p_2=+\infty$, $p_{-2}=-\infty$, $\rho_1=\rho_+$, $\rho_2=-\rho_+$,
$\rho_{-1}=\rho_-$, and $\rho_{-2}=-\rho_-$. Then we conclude that
if $\beta(t)$, $0\le t<\infty$, is a chordal SLE$(4;\rho_+,\rho_-)$
trace started from $(0;0^+,0^-)$, then after a time-change, the
reversal of $(W(\beta(t)))$ has the same distribution as $(\beta(t))$. This is the
case when $\kappa=4$ in Conjecture \ref{conjec} of this paper.

\vskip 5mm

\no{\bf Acknowledgment}. The author would like to thank Oded Schramm for
his suggestion on applying the technique of coupling SLE processes to this work.
The author also thanks the referee's work in reviewing this paper.


\begin{thebibliography}{99}
\bibitem{Ahl} Lars V.\ Ahlfors. {\it Conformal invariants: topics
in geometric function theory}. McGraw-Hill Book Co., New York, 1973.
\bibitem{Dipolar} M.\ Bauer, D.\ Bernard and J.\ Houdayer.
Dipolar stochastic Loewner evolutions. J. Stat. Mech, P03001, 2005.
\bibitem{dim-SLE} V.\ Beffara. The dimension of the SLE curves.
arXiv:math/0211322.
\bibitem{Julien-Duality} Julien Dub\'edat. SLE$(\kappa,\rho)$
martingales and duality, {\it Ann. Probab.}, 33(1):223-243, 2005.
\bibitem{Julien-Comm} Julien Dub\'edat. Commutation relations for SLE,
{\it Comm. Pure Applied Math.}, arXiv:math/0411299.
\bibitem{LawSLE} Gregory F.\ Lawler. {\it Conformally Invariant Processes in the Plane}.
American Mathematical Society, 2005.
\bibitem{LSW1} Gregory F.\ Lawler, Oded Schramm and Wendelin Werner.
Values of Brownian intersection exponents I: Half-plane exponents.
{\it Acta Mathematica}, 187(2):237-273, 2001.
\bibitem{LSW-8/3} Gregory F.\ Lawler, Oded Schramm and Wendelin Werner.
Conformal restriction: the chordal case, {\it J.\ Amer.\ Math.\
Soc.}, 16(4): 917-955, 2003.
\bibitem{LSW-2} Gregory F.\ Lawler, Oded Schramm and Wendelin Werner.
Conformal invariance of planar loop-erased random walks and uniform
spanning trees. {\it Ann. Probab.}, 32(1B):939-995, 2004.
\bibitem{RY} Daniel Revuz and Marc Yor. {\it Continuous Martingales
and Brownian Motion}. Springer-Verlag, 1991.
\bibitem{RS-basic} Steffen Rohde and Oded Schramm. Basic properties of
SLE. {\it Ann.\ Math.}, 161(2):883-924, 2005.
\bibitem{S-SLE} Oded Schramm. Scaling limits of loop-erased random walks
and uniform spanning trees. {\it Israel J. Math.}, 118:221-288,
2000.
\bibitem{SS} Oded Schramm and Scott Sheffield.
Contour lines of the two-dimensional discrete Gaussian free field,
arXiv:math.PR/0605337.
\bibitem{SW} Oded Schramm and David B.\ Wilson. SLE coordinate
changes. {\it New York Journal of Mathematics}, 11:659--669, 2005.
\bibitem{SS-3} Stanislav Smirnov. Towards conformal invariance of 2D lattice
models. {\it International Congress of Mathematicians.},
II:1421-1451, 2006.
\bibitem{SS-6} Stanislav Smirnov. Critical percolation in the plane:
conformal invariance, Cardy's formula, scaling limits. {\it C.\ R.\
Acad.\ Sci.\ Paris S\'er.\ I Math.}, 333(3):239-244, 2001.
\bibitem{Wilson} D.\ B.\ Wilson. Generating random trees more quickly
than the cover time. {\it Proceedings of the 28th ACM Symposium on
the Theory of Computing} 296-303. Assoc.\ Comp.\ Mach., New York.
\bibitem{reversibility} Dapeng Zhan. Reversibility of chordal SLE,
arXiv:0705.1852. To appear in {\it Ann. Probab.}.
\bibitem{LERW} Dapeng Zhan. The Scaling Limits of
Planar LERW in Finitely Connected Domains. {\it Ann. Probab.}, 36(2):467-529, 2008.
\bibitem{thesis} Dapeng Zhan. Random Loewner Chains in Riemann
Surfaces. PhD dissertation in Caltech.
\end{thebibliography}
\end{document}